# A stochastic analysis of subcritical Euclidean fermionic field theories

Francesco De Vecchi

Department of Mathematics, University of Pavia,
Via Adolfo Ferrata 5 27100, Pavia, Italy
*Email:* francescocarlo.devecchi@unipv.it

Luca Fresta

Institute for Applied Mathematics &
Hausdorff Center for Mathematics, University of Bonn,
Endenicher Allee 60 53115, Bonn, Germany
*Email:* fresta@iam.uni-bonn.de

Massimiliano Gubinelli

Mathematical Institute, University of Oxford,
Woodstock Road OX2 6GG Oxford, United Kingdom
*Email:* gubinelli@maths.ox.ac.uk

**Abstract**

Building on previous work on the stochastic analysis for Grassmann random variables, we introduce a forward-backward stochastic differential equation (FBSDE) which provides a stochastic quantisation of Grassmann measures. Our method is inspired by the so-called continuous renormalisation group, but avoids the technical difficulties encountered in the direct study of the flow equation for the effective potentials. As an application, we construct a family of weakly coupled subcritical Euclidean fermionic field theories and prove exponential decay of correlations.

**Keywords:** stochastic quantisation, renormalisation group, non-commutative probability, Euclidean quantum field theory

**A.M.S. subject classification:** 60H30, 81T16

# Table of contents

# 1 Introduction

Let $(\mathcal{M}, \omega)$ be an algebraic probability space, i.e. a unital $C^*$-algebra $\mathcal{M}$ endowed with a normalised, positive linear functional $\omega: \mathcal{M} \to \mathbb{C}$ and let $(\psi(f))_{f \in \mathfrak{h}}$ be a centred Grassmann Gaussian field indexed by the Hilbert space $\mathfrak{h} := L^2(\mathbb{R}^d; \mathbb{C}^2) \oplus L^2(\mathbb{R}^d; \mathbb{C}^2)$ with covariance

$$\omega(\psi(f)\psi(g)) = \langle \Theta f, G g \rangle, \qquad f, g \in \mathfrak{h},$$

where $\Theta$ is the anti-unitary involution on $\mathfrak{h}$ given by $\Theta(f_1 \oplus f_2) := (\bar{f}_2 \oplus \bar{f}_1)$ and $G$ the bounded operator

$$G := (\mathbb{1} - \Delta)^{\gamma - d/2} \oplus -(\mathbb{1} - \Delta)^{\gamma - d/2} \tag{1.1}$$

$\Delta$ being the Laplacian on $L^2(\mathbb{R}^d; \mathbb{C}^2)$ and $0 < \gamma < d/2$. The field $\psi$ is said to be Grassmann because it satisfies anticommutation relations

$$\psi(f)\psi(g) + \psi(g)\psi(f) = 0, \qquad \forall f, g \in \mathfrak{h},$$

hence the need for the algebraic, i.e. non-commutative, setting $(\mathcal{M}, \omega)$ [Mey93, Par92, Bia95]. One can think of $\psi$ as a massive Grassmann Gaussian free field: the fractional exponent in the covariance is there only to model the scaling behaviour one would expect for arbitrary (non-integer) dimension.

We want to prove the existence and cluster properties of a class of Grassmann measures that are constructed as Gibbsian perturbations of the Grassmann Gaussian field $\psi$ on $\mathbb{R}^d$: namely, we want to make sense of the continuous linear functional $\omega^V$ on the algebra $\mathcal{A} \subseteq \mathcal{M}$ of polynomials of $\psi$, and given by

$$\omega^V(A) := \frac{\omega(A e^{V(\psi)})}{\omega(e^{V(\psi)})}, \qquad A \in \mathcal{A}, \tag{1.2}$$

which we call the Gibbsian Grassmann measure with covariance $G$ and *interaction potential* $V(\psi)$. In particular, we will take $V(\psi)$ as a Gross–Neveu-type interaction of the form [GN74]

$$V(\psi) := \int_{\mathbb{R}^d} \left[ \frac{\lambda}{4} (\psi_x)^4 + \mu (\psi_x)^2 \right] \mathrm{d}x,$$

for some constants $\lambda, \mu \in \mathbb{R}$ and with $\psi_x = \psi(\delta_x)$ and $(\psi_x)^4, (\psi_x)^2$ suitable quartic and quadratic homogeneous polynomial of the field $\psi$ evaluated at a point $x \in \mathbb{R}^d$. Due to singularities in the covariance $G$ the random variable $\psi_x$ cannot really be defined and the expression of the interaction potential $V$ is at best a suggestive notation for a more concrete family of approximations which avoid both the small-scale (or ultraviolet, or UV) singularities implicit in considering local monomials and the large-scale (or infrared, or IR) singularities relative to the computation of the integral over all $\mathbb{R}^d$.

Gibbsian Grassmann measures appear naturally in the study of interacting Euclidean fermionic quantum field theories, see, e.g., [GJ87] for a review, but also in the study of condensed-matter systems [GM10, GMP17] and as an emergent description in statistical mechanics, see [Mas08]. Although the Grassmann measures in (1.2) do not describe any relevant quantum field theory nor statistical mechanics system, they are a class of Gross–Neveu-type toy models that exhibit the same mathematical challenges of subcritical (that is, super-renormalisable) Euclidean field theories. They are therefore important from a methodological perspective and should be considered in view of the stochastic analysis in the broader context of super-renormalisable (possibly supersymmetric) field theories with bosons and fermions.

For a rigorous definition of the measure $\omega^V$ we need to regularise the IR and UV singularities and then perform a limit to remove the regularisation parameters. A popular choice is to restrict the randomness in the problem to a finite number of degrees of freedom by working on a finite set of space points. Let $L \in \mathbb{N}$, $\varepsilon \in 2^{-\mathbb{N}}$ and introduce the $d$-dimensional toroidal lattice $\mathbb{T}^d_{L,\varepsilon} := ((\varepsilon\mathbb{Z})/(L\mathbb{Z}))^d$ and a suitable $\mathbb{T}^d_{L,\varepsilon}$-lattice regularisation[1.1] $G^{L,\varepsilon}$ of the continuum covariance $G$ in eq. (1.1). Moreover, let $(\psi^{L,\varepsilon}_x)_{x \in \mathbb{T}^d_{L,\varepsilon}}$ be the Grassmann Gaussian field with covariance $G^{L,\varepsilon}$ (and with the same involution $\Theta$) and discretise the interaction potential as

$$V^{L,\varepsilon}(\psi^{L,\varepsilon}) := \int_{\mathbb{T}^d_{L,\varepsilon}} \left[\frac{\lambda}{4}(\psi^{L,\varepsilon}_x)^4 + \mu^\varepsilon(\lambda)(\psi^{L,\varepsilon}_x)^2\right]dx, \tag{1.3}$$

where $\int_{\mathbb{T}^d_{L,\varepsilon}} \cdot\, dx := \varepsilon^d \sum_{x \in \mathbb{T}^d_{L,\varepsilon}}$ is a normalised counting measure on $\mathbb{T}^d_{L,\varepsilon}$ and where we choose $\mu^\varepsilon(\lambda)$ depending on $\varepsilon, \lambda$ in order to achieve a non-trivial limit as $\varepsilon \to 0$. Due to the singularities in the covariance $G$ we expect to need to consider a diverging family $(\mu^\varepsilon(\lambda))_{\varepsilon > 0}$ in order to compensate for the divergences introduced by the non-linear factor $e^{V^{L,\varepsilon}(\psi)}$ in the averages. To the regularised potential there correspond regularised measures $\omega^{V^{L,\varepsilon}}$, so that our problem becomes the analysis of the weak limit of the families $(\omega^{V^{L,\varepsilon}})_{\varepsilon > 0, L > 0}$:

$$\begin{aligned}\omega^V(A) &:= \lim_{L \to \infty} \lim_{\varepsilon \to 0} \omega^{V^{L,\varepsilon}}(A) \\ &:= \lim_{L \to \infty} \lim_{\varepsilon \to 0} \frac{\omega(A e^{V^{L,\varepsilon}(\psi^{L,\varepsilon})})}{\omega(e^{V^{L,\varepsilon}(\psi^{L,\varepsilon})})}, \quad A \in \mathcal{A}.\end{aligned} \tag{1.4}$$

## 1.1 Methodology and main results

The point of our work is mainly methodological: to construct the Grassmann measure $\omega^V$ we use a novel approach based on a non-commutative forward–backward stochastic differential equation (FBSDE), which provides a *stochastic quantisation* of the measure of interest. Our approach applies to any subcritical Grassmann measure, including the equivalents of the well-known $\varphi^4_2$ and $\varphi^4_3$ bosonic (or commutative) theories.

To describe in detail our strategy, we introduce a differentiable interpolation $(G^{L,\varepsilon}_t)_{t \geqslant 0}$ such that

$$G^{L,\varepsilon}_0 = 0, \qquad G^{L,\varepsilon}_\infty = G^{L,\varepsilon}, \qquad \Theta(G^{L,\varepsilon}_t)^*\Theta = -G^{L,\varepsilon}_t.$$

This interpolation should also suppress momenta larger than $2^s$, so that $G^{L,\varepsilon}_t$ is a bounded operator, uniformly in $L$ and $\varepsilon$, if $t < \infty$. In $(\mathcal{M}, \omega)$ we introduce a Grassmann Brownian martingale (GBM) as the family $(X^{L,\varepsilon}_t)_{t \geqslant 0}$ of centred Grassmann Gaussian random fields indexed by the Hilbert space[1.2] $\mathfrak{h}_{L,\varepsilon} := L^2(\mathbb{T}^d_{L,\varepsilon}; \mathbb{C}^2) \oplus L^2(\mathbb{T}^d_{L,\varepsilon}; \mathbb{C}^2)$ and with covariance $G^{L,\varepsilon}_t$,

$$\omega(X^{L,\varepsilon}_t(f) X^{L,\varepsilon}_s(g)) = \langle \Theta f, G^{L,\varepsilon}_{t \wedge s} g\rangle_{\mathfrak{h}_{L,\varepsilon}}.$$

Note that the $X^{L,\varepsilon}_\infty$ is distributed in law as the Gaussian free field $\psi^{L,\varepsilon}$. See Section 2 for the detailed definitions and properties of these objects. From a stochastic perspective, the crucial point to note is that the parameter $t$ is not a physical nor a stochastic time. Rather, following the core idea of RG techniques, it is a continuous flow parameter associated with the scale decomposition of the problem, which we choose in such a way that ultraviolet scales correspond to $t \to \infty$.

---

1.1. Because of technical reasons, we will avoid simply replacing $\Delta$ by the lattice Laplacian $\Delta_{L,\varepsilon}$ on $L^2(\mathbb{T}^d_{L,\varepsilon}; \mathbb{C}^2)$, see Definition 3.16 in Section 3.3.

1.2. The space $L^2(\mathbb{T}^d_{L,\varepsilon}; \mathbb{C}^2)$ is associated with $\varepsilon^d$ times the counting measure: if $f, g \in L^2(\mathbb{T}^d_{L,\varepsilon}; \mathbb{C}^2)$ we write $\langle f, g\rangle_{L^2(\mathbb{T}^d_{L,\varepsilon}; \mathbb{C}^2)} = \sum_i \sum_{x \in \mathbb{T}^d_{L,\varepsilon}} \varepsilon^d \overline{f_i(x)} g_i(x)$.

As we have already announced, we would like to study the interacting measure (1.4) via stochastic quantisation, that is, by identifying it as the marginal law of a suitable stochastic process coupled with the Grassmann Gaussian field. In particular, we will prove that for any nice enough function $P$ the following identity holds true

$$\omega^V(P(X_\infty^{L,\varepsilon})) = \omega(P(\Psi_\infty^{L,\varepsilon})), \tag{1.5}$$

provided that the process $(\Psi_s^{L,\varepsilon})_s$ solves the following forward-backward stochastic differential equation (FBSDE) on $[0,\infty]$

$$d\Psi_s^{L,\varepsilon} = -\dot{G}_s^{L,\varepsilon} \omega_s(DV^{L,\varepsilon}(\Psi_\infty^{L,\varepsilon}))ds + dX_s^{L,\varepsilon}, \qquad \Psi_0^{L,\varepsilon} = 0, \tag{1.6}$$

where $DV^{L,\varepsilon}$ denotes the functional derivative of $V^{L,\varepsilon}$ and where $\omega_t(\cdot)$ denotes the conditional expectation with respect to a filtration such that $(X_t^{L,\varepsilon})_t$ is adapted, see Section 2.

It is important to note that the drift term in (1.6) satisfies the identity

$$\omega_s(DV^{L,\varepsilon}(\Psi_\infty^{L,\varepsilon})) = DV_s^{L,\varepsilon}(\Psi_s^{L,\varepsilon}),$$

where $V_t^{L,\varepsilon}$ is the solution of the Hamilton–Jacobi–Bellman (HJB) flow equation

$$\partial_t V_t^{L,\varepsilon} + \frac{1}{2} D_{\dot{G}_t^{L,\varepsilon}}^2 V_t^{L,\varepsilon} + \frac{1}{2}(DV_t)_{\dot{G}_t^{L,\varepsilon}}^2 = 0, \tag{1.7}$$

where $D_{\dot{G}_t^{L,\varepsilon}}^2$ is the functional Laplacian and where $(\cdot)_{\dot{G}_t^{L,\varepsilon}}^2$ is a suitable quadratic form, associated with the operator $\dot{G}_t^{L,\varepsilon}$. The typical approach in the physics and mathematical physics literature is based on the study of (1.7), see below. However, we do not want to make this substitution: the FBSDE is a powerful tool that allow us to truncate the flow equation (1.7) for the effective force $\omega_s(DV^{L,\varepsilon}(\Psi_\infty^{L,\varepsilon}))$ in a suitable way, and to control the remainder due to the truncation for any subcritical theory. In fact, we decompose the drift term as

$$\omega_s(DV^{L,\varepsilon}(\Psi_\infty^{L,\varepsilon})) = F_s^{L,\varepsilon}(\Psi_s^{L,\varepsilon}) + R_s^{L,\varepsilon},$$

where the "remainder process" $R_s^{L,\varepsilon}$ solves the following self-consistent equation, for $s > 0$

$$R_s^{L,\varepsilon} = \int_s^\infty \omega_s(\mathcal{H}_r[F_r^{L,\varepsilon}](\Psi_r^{L,\varepsilon}))dr + \int_s^\infty \omega_s(DF_r^{L,\varepsilon}(\Psi_r^{L,\varepsilon}) \cdot \dot{G}_r^{L,\varepsilon} R_r^{L,\varepsilon})dr, \tag{1.8}$$

with

$$\mathcal{H}_r[F_r^{L,\varepsilon}] := \partial_r F_r^{L,\varepsilon} + \frac{1}{2} D_{\dot{G}_r^{L,\varepsilon}}^2 F_r^{L,\varepsilon} + DF_r^{L,\varepsilon} \cdot \dot{G}_r^{L,\varepsilon} F_r^{L,\varepsilon}.$$

We have now the freedom to choose $F_s^{L,\varepsilon}$ in any convenient way. If we let $F_s^{L,\varepsilon}$ be the solution of the HJB flow equation, then $\mathcal{H}_r[F_r^{L,\varepsilon}]$ vanishes and so does $R_s^{L,\varepsilon}$. More generally, we can take $F_s^{L,\varepsilon}$ to be the solution of a simpler flow equation, provided that we can control eq. (1.8). In fact, once we have good estimates for $F_s^{L,\varepsilon}$, we can prove the existence and uniqueness of a pair $(\Psi^{L,\varepsilon}, R^{L,\varepsilon})$ solving eqs. (1.6) and (1.8) by a fixed-point argument, for which $\lambda$ small is needed.

To provide a precise statement, for $a, \eta \geqslant 0$ we introduce the following topology on fields[1.3] on $\mathbb{R}^+ \times \mathbb{R}^d$, that is, maps $\psi \colon \mathbb{R}^+ \times \mathbb{R}^d \to \mathcal{M}^4$, $(t,x) \mapsto (\psi_{t,x,\mu})_{\mu \in \{\downarrow,\uparrow\} \times \{\pm\}}$:

$$\|\psi\|_{a,\eta} := \sup_\mu \sup_{(t,x) \in \mathbb{R}^+ \times \mathbb{R}^d} \varrho_\eta(x) 2^{-at} \|\psi_{t,x,\mu}\|. \tag{1.9}$$

---

1.3. If fields are defined on $\mathbb{R}^+ \times \mathbb{T}_{L,\varepsilon}^d$, they can be suitably extended to $\mathbb{R}^+ \times \mathbb{R}^d$, see Section 3.3.

where $\varrho_\eta(x) := (1+|x|^2)^{-\eta/2}$. Note that when $t \in \mathbb{R}^+$ is the scale parameter described above, $-a$ measures the Hölder regularity in space of the field, in a suitable weighted $L^\infty$ topology. We also remark that, in order to carry out the fixed-point argument described above, we also need to control the derivatives of $\Psi^{L,\varepsilon}$, see Section 4.3 for details; for the moment (1.9) suffices to state the theorem below. On fields on $\mathbb{R}^d$, we consider the weighted Besov spaces $B^a_{p,q,\eta}(\mathbb{R}^d;\mathcal{M}^4)$, see, e.g., [Ama19, BCD11] for the definition of Besov spaces and [ABDG20] for their use in the analysis of fermionic systems, where $\eta \in \mathbb{R}$ indicates integration with weight $\varrho_\eta$; we denote the Besov norms simply by $\|\cdot\|_{B^a_{p,q,\eta}}$. If $\eta = 0$, as usual, we abridge the notation to $B^a_{p,q}$ and $\|\cdot\|_{B^a_{p,q}}$.

**Theorem 1.1.** *Let $d \in \mathbb{N}$, $\gamma < \min\{d/4, 1\}$. Then, there exist $\lambda_0 = \lambda_0(d,\lambda)$ and a function $\mu^\varepsilon_\infty$: $\mathbb{R}^+ \to \mathbb{R}$ such that if $\lambda \leq \lambda_0$ and $\mu^\varepsilon(\lambda) = \mu^\varepsilon_\infty(\lambda)$ then for any $L \in \mathbb{N}$ and $\varepsilon \in 2^{-\mathbb{N}}$ if $V^{L,\varepsilon}$ is chosen as in (1.3), the FBSDE in (1.6) has unique global solution $\Psi^{L,\varepsilon}$. The sequence $(\Psi^{L,\varepsilon})_{L,\varepsilon}$ is bounded in the $\|\cdot\|_{\gamma,0}$ topology and Cauchy in the $\|\cdot\|_{\gamma+\theta,\eta}$ topology for $\eta,\theta > 0$ small enough. Denoting by $\Psi$ the limit we have*

$$\|\Psi - \Psi^{L,\varepsilon}\|_{\gamma+\theta,\eta} \lesssim_{\theta,\eta} \varepsilon^\theta + L^{-\eta}. \tag{1.10}$$

*Furthermore, the sequence $(\Psi^{L,\varepsilon}_\infty)_{L,\varepsilon}$ is bounded in $B^{-\gamma-\rho}_{\infty,\infty}(\mathbb{R}^d;\mathcal{M}^4)$ for any $\rho > 0$ and Cauchy in $B^{-\gamma-\rho}_{\infty,\infty,\eta}(\mathbb{R}^d;\mathcal{M}^4)$ for $\rho,\eta > 0$ small enough. Denoting the limit $\Psi_\infty$, for $\theta > \rho > 0$ small enough we have*

$$\|\Psi_\infty - \Psi^{L,\varepsilon}_\infty\|_{B^{-\gamma-\theta}_{\infty,\infty,\eta}} \lesssim_{\theta,\eta} (\varepsilon^{\theta-\rho} + L^{-\eta}). \tag{1.11}$$

**Remark 1.2.** The constraint $\gamma < 1$ is for technical reasons and can be removed provided that further "counter-terms" are included in the potential, e.g., of the form $\int_{\mathbb{T}^d_{L,\varepsilon}} v^\varepsilon(\lambda)(\partial_x X^{L,\varepsilon}_{t,x})^2 dx$, where $\partial_x$ is the lattice derivative.

Theorem 1.1 allows us to identify the algebra of observables $\mathcal{A}$ on which the measures $\omega^{V^{L,\varepsilon}}$ can be defined via the stochastic quantisation equation (1.5). This is the polynomial algebra generated by the Grassmann Gaussian free field $\psi$

$$\mathcal{A} := \operatorname{span}\{\mathbb{1}, \psi(f_1) \cdots \psi(f_n) | n \in \mathbb{N}, (f_i)_{i=1}^n \subseteq B^{\gamma^+}_{1,1,0^-}\}, \tag{1.12}$$

see Section 4.3, where $B^{\gamma^+}_{1,1,0^-} := \bigcup_{a > \gamma, \eta < 0} B^a_{1,1,\eta}$. Since $B^a_{1,1,\eta} \subset B^{a'}_{1,1,\eta'}$ for $a \geq a'$ and $\eta \leq \eta'$, if $f \in B^{\gamma^+}_{1,1,0^-}$, with abuse of notation we write $\|f\|_{B^{\gamma^+}_{1,1,0^-}} := \|f\|_{B^{\bar{a}}_{1,1,\bar{\eta}}}$, where $\bar{a}$ and $\bar{\eta}$ are such that $f \in B^{\bar{a}}_{1,1,\bar{\eta}}$ but $f \notin B^{\bar{a}^+}_{1,1,\bar{\eta}^-}$. We prove existence, cluster properties and short-distance divergence of the weak limit of $\omega^{V^{L,\varepsilon}}$ on $\mathcal{A}$.

**Theorem 1.3.** *Under the same assumption of Theorem 1.1, for any $L \in \mathbb{N}$, $\varepsilon \in 2^{-\mathbb{N}}$, $\lambda \leq \lambda_0$, the measures $\omega^{V^{L,\varepsilon}}$ exist and have a weak limit $\omega^V$ on $\mathcal{A}$, as $\varepsilon \to 0$ and $L \to \infty$. Furthermore:*

1. *Consider $m_1, m_2 \in \mathbb{N}$, and let $((f^{(i,k)})_{k=1}^{m_i})_{i=1,2} \subseteq B^{\gamma^+}_{1,1,0^-}(\mathbb{R}^d;\mathbb{C}^4)$. Let $\operatorname{Cov}_{\omega^V}(A;B) := \omega^V(AB) - \omega^V(A)\omega^V(B)$ and let $\psi$ be the Grassmann Gaussian field with covariance (1.1). Then, letting $D_i := \bigcup_{k=1}^{m_i} \operatorname{supp}(f^{(i,k)})$, for some universal $c > 0$ we have*

$$\left|\operatorname{Cov}_{\omega^V}\left(\prod_{k=1}^{m_1} \psi(f^{(1,k)}); \prod_{k=1}^{m_2} \psi(f^{(2,k)})\right)\right| \lesssim_{d,\gamma,\lambda} \left(\prod_{i=1}^2 \prod_{k=1}^{m_i} \|f^{(i,k)}\|_{B^{\gamma^+}_{1,1,0^-}}\right) e^{-c\operatorname{dist}(D_1,D_2)}.$$

2. *Fix some $\chi \in C^\infty_c(\mathbb{R}^d)$ such that $\|\chi\|_{L^1} = 1$ and for any $\varepsilon > 0$ and $x \in \mathbb{R}^d$ define $\chi^\varepsilon_x := \varepsilon^{-d}\chi(\varepsilon^{-1}(\cdot-x))$. Then, for any $x \neq y \in \mathbb{R}^d$ and $\varepsilon > 0$ small enough,*

$$|\omega^V(\psi(\chi^\varepsilon_x)\psi(\chi^\varepsilon_y)) - \omega(\psi(\chi^\varepsilon_x)\psi(\chi^\varepsilon_y))| \lesssim \lambda (\varepsilon \vee |x-y|)^{-2\gamma}.$$

**Remark 1.4.** In particular, Theorem 1.3 implies the exponential decay of the correlation function: for $m_1 = m_2 = 1$, writing $f_i = f^{(i,1)}$, $i = 1, 2$, then

$$|\omega^V(\psi(f_1)\psi(f_2))| \lesssim_{d,\gamma,\lambda,f_1,f_2} e^{-c\,\mathrm{dist}(D_1,D_2)}.$$

Furthermore, since the non-interacting correlation function diverges as $\omega(\psi(\chi_x^\varepsilon)\psi(\chi_y^\varepsilon)) \sim (\varepsilon \vee |x-y|)^{-2\gamma}$, Theorem 1.3 implies that for $\lambda$ small enough

$$\omega^V(\psi(\chi_x^\varepsilon)\psi(\chi_y^\varepsilon)) \sim (\varepsilon \vee |x-y|)^{-2\gamma}.$$

## 1.2 Comparison with the literature

The problem of taking the IR and UV limits of the renormalised model with interaction (1.3) is a well-known mathematical question in constructive quantum field theory. Many techniques have been developed both for bosonic and fermionic theories. We will review below some of the relevant literature.

**Stochastic quantisation.** Stochastic quantisation is understood in the broad sense of having a map which transports a measure of our choice into a target measure of interest, in this case a Gibbsian measure of the form $\omega^V$. In our particular case, the measure of choice is given by a non-commutative probability space $(\mathcal{M}, \tilde{\omega})$ endowed with a suitable Gaussian Grassmann field $X$ and another Grassmann field $\Psi$ (non-Gaussian, in general) which will belong to the sub-algebra $\mathcal{M}_X \subseteq \mathcal{M}$ generated by $X$ and such that

$$\mathrm{Law}_{\tilde{\omega}}(\Psi)(P) := \tilde{\omega}(P(\Psi)) = \omega^V(P(\psi)) = \mathrm{Law}_{\omega^V}(\psi)(P),$$

for a sufficiently large class of functions $P$ on Grassmann fields. The usefulness of this construction is that one can identify sufficiently nice maps $F: X \to \Psi$ which perform the push-forward of the Gaussian law of $X$ under $\tilde{\omega}$ to the law of $\psi$ under $\omega$, i.e. $\mathrm{Law}_{\omega^V}(\psi) = F_\# \mathrm{Law}_{\tilde{\omega}}(X)$. The nice features of the map $F$ together with the Gaussianity of $\mathrm{Law}_{\tilde{\omega}}(X)$ are the key to the construction and analysis of $\mathrm{Law}_{\omega^V}(\psi)$.

From a more conceptual viewpoint, the stochastic quantisation, as understood above, is a tool to perform the stochastic analysis of a given measure or process and this paper is one among a series of recent results which try to push forward a programme of systematically developing the stochastic analysis of Euclidean quantum field theories or more generally of irregular random fields parametrised by multidimensional continuous index sets like $\mathbb{R}^d$ and possessing particular properties of locality and symmetry. In particular, let us mention the following works:

- The stochastic analysis of Grassmann random fields, in the language of the present paper was initiated in [ABDG20] where the infinite volume limit of some weakly interacting Grassmann Gibbsian measures was established via parabolic stochastic quantisation, i.e. stochastic quantisation based on a Langevin-type equation.

- The stochastic quantisation of the $\Phi_3^4$ bosonic model has been pushed forward in [GH21]. While this work has been able to establish most of the Osterwalder–Schrader [GJ87] axioms, the uniqueness and rotation invariance of the limit is still an open problem.

- A variational method for stochastic quantisation in the bosonic setting has been introduced in [BG20, BG21b, BG21a, Bar22] and used to prove estimates for $\Phi_3^4$ in [CGW20] towards establishing the existence of a phase transition. Our contribution in this paper can be understood as the fermionic version of the variational approach. In the bosonic setting the FBSDE indeed appears as the Euler–Lagrange equation of the optimal control problem for the interacting measure, see also [GM24] for recent work on the FBSDE in this setting.

- An elliptic variant of stochastic quantisation for bosonic theories has been introduced and studied in [AVG20, ADG21, BD21], including the full construction of an EQFT in two space dimensions and the verification of the Osterwalder–Schrader axioms.

- Stochastic quantisation is a very active field of modern stochastic analysis, with recent advances in the analysis of low-dimensional quantum (Euclidean) Yang–Mills theories [CCHS20, CCHS22], large-$N$ limit of $O(N)$ vector models in $d = 3$ [SZZ21b] and study of the perturbation expansion of correlations in $\Phi_2^4$ [SZZ21a].

- The point of view on stochastic quantisation as the search for good representation as push-forwards of the target measures is parallel to a recent analysis from the point of view of functional inequalities and optimal transport, see e.g. [She22], inspired by a renormalisation group analysis of [BB21, BD22, BBD23] which is very much related to our use of the flow equations and the renormalisation group ideas as discussed below.

- One common feature of stochastic quantisation methods is that they naturally provide a coupling between the Gaussian free field and the interacting field. This is also true in our Grassmann framework, which is a novelty with respect to other constructive fermionic technology. This coupling has been recently used to study the rigorous properties of the measure [BH20, DGT22].

Our analysis builds on the Grassmann probabilistic setting developed in [ABDG20] and pioneered in [OS72]. See also [Lot87] for a similar approach to Fermi fields, in the context of supersymmetric Wiener integrals. In addition to what was done in [ABDG20, OS72], we are here able to remove the ultraviolet cut-off and thus to describe ultraviolet singular super-renormalisable fermionic field theories. Let us point out right now that this is not enough to study some of the well-known fermionic field theories, e.g., the Gross–Neveu model in $d = 2, 3$ [GK85a, FMRS86] and the Yukawa model in $d = 2$ [Les87], which, within the formalism of the present paper, behave as critical theories which require a more detailed handling of the small-scale contributions to the drift of the stochastic equation. We leave to future work the analysis of these more delicate situations.

**The continuous renormalisation group.** The stochastic quantisation that we use is introduced via a scale decomposition of the Gaussian process according to the basic ideas of the renormalisation group (RG). In the constructive QFT literature, RG is one of the most successful methods to study Euclidean fermionic theories (i.e. Grassmann measures). Indeed, Grassmann measures have early on been understood to be a nice setting in which to develop constructive tools and RG techniques thanks to the intrinsic *bounded* nature of the relevant Gaussian variables, in striking contrast with the usual (bosonic) Gaussians which have unbounded support.

More generally, RG is a crucial tool in the construction of Euclidean quantum field theories and in the study of statistical mechanics and critical phenomena. It was pioneered in the work of Wilson [Wil71a, Wil71b], building on the work of Kadanoff [Kad66], where the idea of progressive integration of scales was first introduced for studying critical phenomena. Since then, RG techniques have witnessed a very intense and diverse development, depending on whether they are used to study, e.g., Euclidean scalar measures [BCG+80, GK85b, FMRS87, BY90, Abd07, BBS19], Hamiltonian systems [ED81, MKS82, FT92, GGM95], Grassmann measures [GK85a, FMRS86, Les87, BG90, BGM92, FT04], lattice gauge theories [Bal87, Bal88, Dim20], supersymmetric measures [BBS15b, BBS15a, AFP20], or stochastic differential equations [Kup16, KM17, Duc21, Duc22] and depending on whether the separation of scales is done in a continuous or discrete fashion. Our list of references is by no means exhaustive and should only be understood as a sample for the interested reader.

A well-known approach in the case of a continuum of scales is the Polchinski flow equation [Pol84], that is, the infinite-dimensional version of the HJB eq. (1.7). While this HJB equation has been an effective tool to estimate the perturbative expansion of bosonic Euclidean field theories [Pol84, Sal99], the rigorous study of eq. (1.7) at the non-perturbative level is very challenging. An important contribution in this sense was made in [BK87] and later applied to a slightly different equation [BY90, BBS15b, BBS15a], where the basic idea is to show that the solution of the HJB equation is majorised, in a suitable topology, by a quantity solving a simpler first-order equation, without a functional Laplacian.

In the context of Grassmann measures the direct study of the flow equation (1.7) is quite challenging (and essentially an open problem to our knowledge, see also [GMR21]). The validity of the flow equation in relevant examples and without ultraviolet cutoffs has been established by Disertori–Rivasseau [DR00] using a careful analysis of the tree expansion, see also the recent work [Gre24]. In most of the literature on the RG of Euclidean Fermionic theories, the technical difficulties posed by the continuous equation are avoided by studying instead a discrete flow, that is, a sequence $(V_{t_n})_{n \geqslant 0}$ with the points $(t_n)_{n \geqslant 0}$ sufficiently separated. The discrete flow equation, despite being more cumbersome, turns out to be controllable in a surprisingly wide range of regimes because of the Pauli principle [GK85a], which later assumed the form of the celebrated determinant bounds [Les87], see also [GMR21] and [FMRS86]. More recently, however, the continuous equation for fermions has received more attention in the spirit of [BK87], see [KM22].

**A synergy.** The key contribution of the present paper, apart from the explicit constructions and bounds, is the idea that the synergy of RG and stochastic quantisation, in the form of a FBSDE, provides a powerful tool to analyse Gibbsian measures. Indeed, here we control the solution of the FBSDE by means of a suitable flow equation for the change of the drift with respect to the scale of the decomposition. However, in contrast with other works using a RG approach, and thanks to the presence of the FBSDE, we would only need to solve Polchinski's equation in an *approximate way*, which results in great simplification of the analysis. The approximate flow equation is studied in a space of polynomial Grassmann functionals with controlled locality properties within an analytic setup largely inspired by the work [GMR21]. Note that flow equations have been recently applied to the study of super-renormalisable stochastic partial differential equations [Duc21, Duc22] as well, but without taking full advantage of the presence of the SPDE. By transposing the ideas developed in the present paper, and in particular the use of the approximate flow equation in conjunction with the analysis of an associated equation for the fields, the paper [DGR24] obtains the stochastic quantisation of the bosonic fractional $\Phi_3^4$ theory in the full subcritical regime and in the infinite-volume limit. Finally, we also mention the recent application of flow equations to the generalised KPZ equation [CF24].

**Comparison with Berezin integration.** In the standard approach, Grassmann measures on a finite dimensional Grassmann algebra are described via the Berezin integration [Sal99, KTF02, Mas08]. By analogy with the commutative setting, one can think of the Berezin integral as a "flat" measure for Grassmann variables, i.e. the equivalent of the Lebesgue measure. More in detail, consider the Grassmann algebra generated by $(\psi_{x,\sigma}^\rho)_{x \in \mathbb{T}_{L,\varepsilon}^d, \sigma=\uparrow,\downarrow}^{\rho=\pm}$, that is, the unital complex algebra whose generators satisfy for any $x, x' \in \mathbb{T}_{L,\varepsilon}^d$, $\rho, \rho' \in \{\pm\}$, $\sigma \in \{\uparrow, \downarrow\}$

$$\psi_{x,\sigma}^\rho \psi_{x',\sigma'}^{\rho'} + \psi_{x',\sigma'}^{\rho'} \psi_{x,\sigma}^\rho = 0. \tag{1.13}$$

Let us denote by $\mathcal{G}_{L,\varepsilon}^d$ said Grassmann algebra and note that $\mathcal{G}_{L,\varepsilon}^d$ is a finite-dimensional linear space thanks to (1.13).

By abuse of language we call Grassmann measure any linear functional on $\mathcal{G}_{L,\varepsilon}^d$. The Berezin integral $\omega_{\text{Ber}}$ is the Grassmann measure defined on the linear generators of the algebra as

$$\omega_{\text{Ber}}\left(\prod_{x\in\mathbb{T}_{L,\varepsilon}^d,\sigma=\uparrow,\downarrow}\psi_{x,\sigma}^-\psi_{x,\sigma}^+\right)=1,$$

and equal to zero on any other linearly independent element of $\mathcal{G}_{L,\varepsilon}^d$. Note that the order in the product above is unimportant because $\psi_{x,\sigma}^-\psi_{x,\sigma}^+$ are commuting elements of $\mathcal{G}_{L,\varepsilon}^d$. It is customary to formally write the Berezin integration as

$$\int\left[\prod_{x\in\mathbb{T}_{L,\varepsilon}^d,\sigma=\uparrow,\downarrow}\mathrm{d}\psi_{x,\sigma}^+\mathrm{d}\psi_{x,\sigma}^-\right]P(\psi):=\omega_{\text{Ber}}(P(\psi)).$$

This is consistent with interpreting the symbol $\int\mathrm{d}\psi_{x,\sigma}^\varepsilon\cdot$ as the anticommuting derivative $\frac{\partial}{\partial\psi_{x,\sigma}^\varepsilon}$ on $\mathcal{G}_{L,\varepsilon}^d$; in fact, the Berezin integration is also written as $\omega_{\text{Ber}}(\cdot)=\prod_{x,\sigma}\frac{\partial}{\partial\psi_{x,\sigma}^+}\frac{\partial}{\partial\psi_{x,\sigma}^-}\cdot\,$. As long as $\varepsilon>0$ and $L\in\mathbb{N}$, we can always write the law $\omega^{\rho,\Psi}:=\text{Law}_\rho(\Psi)$ in the state $\rho$ of any Grassmann random variable $\Psi$ with values in $\mathcal{G}_{L,\varepsilon}^d$ in terms of the Berezin measure as $\omega^{\rho,\Psi}(\cdot)=\omega_{\text{Ber}}(\nu\cdot):=\int\mathrm{d}\nu(\psi)\cdot$ for some "density" $\nu$ taking values in $\mathcal{G}_{L,\varepsilon}^d$, and, by abuse of notation, refer also to $\mathrm{d}\nu$ as a Grassmann measure. Therefore, as long as $L\in\mathbb{N}$ and $\varepsilon>0$, the law $\mathrm{d}\nu$ of $\psi$ under the measure $\omega^V$ as described in (1.4) can be written in terms of the Berezin integral as follows:

$$\mathrm{d}\nu(\psi)\propto\left[\prod_{x\in\mathbb{T}_{L,\varepsilon}^d,\sigma=\uparrow,\downarrow}\mathrm{d}\psi_{x,\sigma}^+\mathrm{d}\psi_{x,\sigma}^-\right]e^{-(\psi,(G^{L,\varepsilon})^{-1}\psi)+\int_{\mathbb{T}_{L,\varepsilon}^d}\left[\frac{\lambda}{4}(\psi_x)^4+\mu^\varepsilon(\lambda)(\psi_x)^2\right]\mathrm{d}x},\qquad(1.14)$$

where we used the notation $(\psi,A\psi):=\sum_{\sigma,\sigma'}\int_{\mathbb{T}_{L,\varepsilon}^d\times\mathbb{T}_{L,\varepsilon}^d}\psi_{x,\sigma}^+A_{\substack{x,y\\\sigma,\sigma'}}\psi_{y,\sigma'}^-\mathrm{d}x\mathrm{d}y$, and set $(\psi_x)^2:=\sum_\sigma\psi_{x,\sigma}^+\psi_{x,\sigma}^-$, $(\psi_x)^4:=((\psi_x)^2)^2$. In other words, (1.14) makes the Gibbsian Grassmann measure (1.2) quite concrete in terms of functionals on an abstract finite dimensional Grassmann algebra and it is the customary object of investigation in the mathematical physics literature.

The language of Berezin integration becomes cumbersome when dealing with Grassmann stochastic processes which cannot be faithfully realised on finite dimensional Grassmann algebras. Furthermore, although using an infinite-dimensional abstract Grassmann algebra is possible [KTF02], the useful topology generated by a Gaussian covariance essentially reproduce the analytic setting we realised quite compactly here. See also [ABDG20] for further discussion on similar matters and on a review of previous works which tried to make sense of Grassmann stochastic analysis without using our $C^*$-algebraic topological setting.

We look at the $C^*$-algebra $\mathcal{M}$ as a convenient place where to realise all the Grassmann random variables we need; in this sense, it is a suitable substitute for the finite-dimensional Berezin integral. An analogy drawn from the standard commutative probabilistic setting can help the reader get our point of view: while in finite dimensions all Gaussian measures are absolutely continuous with respect to the Lebesgue measure, in infinite dimensions this is not true anymore and there is not a privileged "base" measure. Similarly, when dealing with infinite-dimensional Grassmann Gaussian fields and processes, the appeal of the language of Berezin integration is quite limited. This observation was clear to some researchers in the field, in particular Osterwalder–Schrader [OS72, Ost73, OS73] and Lott [Lot87].

**Structure of the paper.** In Section 2 we introduce the general formalism of filtered non-commutative probability spaces and characterise in detail Grassmann Brownian martingales. With respect to [ABDG20], we focus on the conditional expectation, which is a crucial ingredient for the derivation of the FBSDE. The latter is the content of Section 3, in which we also describe the model in more detail, the topology on the Grassmann fields and prove some bounds on the covariance of the GBM used in our FBSDE. Finally, in Section 4 we solve the FBSDE by using a truncated flow equation for the drift and, as an application, we prove exponential clustering of the correlation functions.

**Acknowledgements.** The work of L.F. has been supported by the Swiss National Science Foundation (SNSF), grant number 200160. This work has been partially funded by the German Research Foundation (DFG) under Germany's Excellence Strategy - GZ 2047/1, project-id 390685813. This paper has been written with GNU T<sub>E</sub>X<sub>MACS</sub> (www.texmacs.org).

## 2 Grassmann stochastic analysis

In this section, we recall some aspects of non-commutative (or algebraic) probability spaces that are necessary for the formulation of the Grassmann FBSDE. In particular, we introduce the concept of a Grassmann process, a filtration and a conditional expectation associated with the latter and describe Grassmann Brownian martingales. We refer to [ABDG20] for further details. The problem of showing that Grassmann random variables can actually be represented as operators acting on suitable Hilbert spaces, see [Ost73, OS73, ABDG20], is postponed to Appendix A.

### 2.1 Grassmann probability

A non-commutative probability space (NPS) is the pair $(\mathcal{M}, \omega)$, consisting of a unital $C^*$-algebra $\mathcal{M}$ of operators acting on a separable Hilbert space and $\omega$ a state on $\mathcal{M}$, that is, a positive (hence continuous) normalised linear functional. We are interested in generalised processes or fields indexed by a complex Hilbert space and satisfying exact anticommutation relations.

**Definition 2.1.** *(Grassmann field). Let $(\mathcal{M}, \omega)$ be a NPS and $\mathcal{H}$ a complex separable Hilbert space. A Grassmann field indexed by $\mathcal{H}$ is a linear map $\psi: \mathcal{H} \to \mathcal{M}$ such that*

$$\{\psi(f), \psi(g)\} := \psi(f)\psi(g) + \psi(g)\psi(f) = 0, \qquad \forall f, g \in \mathcal{H}.$$

**Remark 2.2.** Denote by $\bigwedge \mathcal{H}$ the Grassmann algebra (also known as exterior algebra) associated with the vector space $\mathcal{H}$, that is, the quotient of the tensor algebra $\bigoplus_{n \geqslant 0} \mathcal{H}^{\otimes n}$ by the two-sided ideal generated by $\{f \otimes f | f \in \mathcal{H}\}$. It is implied that a Grassmann field extends to an algebra homomorphism of the Grassmann algebra $\bigwedge \mathcal{H}$ into $\mathcal{M}$, by setting $\psi(f \wedge g) := \psi(f)\psi(g)$. Equivalently, the set of Grassmann fields on $(\mathcal{M}, \omega)$ indexed by $\mathcal{H}$ can be identified with $\mathrm{Hom}(\bigwedge \mathcal{H}, \mathcal{M})$. Note that the algebra of observable $\mathcal{A}$ introduced in (1.12) is simply $\mathcal{A} = \psi(\bigwedge B^{\gamma^+}_{1,1,0^-})$, $\psi$ being the Grassmann Gaussian free field. Note, however, that $\psi(\bigwedge \mathcal{H})$ is not a Grassmann algebra in general, since $\psi(\bigwedge \mathcal{H})$ could have additional relations.

Given two Grassmann fields $\psi_1$ and $\psi_2$, in general one has that $\{\psi_1(f), \psi_2(g)\} \neq 0$, for some $f, g \in \mathcal{H}$, that is, the generalised random process $\zeta$ indexed by $\mathcal{H} \oplus \mathcal{H}$ and defined by $\zeta(f \oplus g) := \psi_1(f) + \psi_2(g)$, is not a Grassmann field. For example, the adjoint field $\psi^*$ is still a Grassmann field but does not satisfy anticommutation relations with $\psi$. This motivates the following definition, see [ABDG20].

**Definition 2.3.** *(Compatibility). Let $\psi_1$ and $\psi_2$ be two Grassmann fields indexed by $\mathcal{H}_1$ and $\mathcal{H}_2$ respectively. We say that they are compatible or jointly Grassmann fields iff the random field defined by $\zeta(f \oplus g) := \psi_1(f) + \psi_2(g)$, for any $(f,g) \in (\mathcal{H}_1, \mathcal{H}_2)$ is a Grassmann field indexed by $\mathcal{H}_1 \oplus \mathcal{H}_2$.*

A crucial class of Grassmann fields is given by Gaussian Grassmann fields. More specifically, consider a complex separable Hilbert space $\mathcal{H}$ with scalar product $\langle \cdot, \cdot \rangle_\mathcal{H}$ and with conjugation $\Theta$, that is, an anti-linear map such that $\langle f, g \rangle_\mathcal{H} = \langle \Theta g, \Theta f \rangle_\mathcal{H}$ for any $f, g \in \mathcal{H}$, and such that $\Theta^2 = \text{Id}$. We introduce the bilinear, symmetric and non-degenerate form

$$\langle\!\langle \cdot, \cdot \rangle\!\rangle_\Theta := \langle \Theta \cdot, \cdot \rangle_\mathcal{H},$$

and for any bounded linear operator $A \in \mathcal{B}(\mathcal{H})$ define the $\Theta$-transpose $A^\Theta := \Theta A^* \Theta$, which satisfies $\langle\!\langle f, Ag \rangle\!\rangle_\Theta = \langle\!\langle A^\Theta f, g \rangle\!\rangle_\Theta$. Because of the antisymmetric nature of the Grassmann fields, Gaussian Grassmann fields are associated with a covariance $G \in \mathcal{B}(\mathcal{H})$ that is $\Theta$-antisymmetric:

$$G^\Theta = -G. \tag{2.1}$$

Let us define the truncated expectation, or cumulant, of a Grassmann field.

**Definition 2.4.** *(Cumulant). Let $(\mathcal{M}, \omega)$ be a NPS and let $\psi$ be a Grassmann field in it indexed by a complex separable Hilbert space $\mathcal{H}$. We define the cumulants $(\mathcal{K}_n[\psi])_{n \in \mathbb{N}}$ of $\psi$, as the family of multilinear maps given by $\mathcal{K}_1[\psi](f) := \omega(\psi(f))$ and for $n \geq 2$ by recursively solving*

$$\omega(\psi(f_1) \cdots \psi(f_n)) =: \sum_{\Pi \in \text{Partitions}} (-)^{\pi(\Pi)} \prod_{I \in \Pi} \mathcal{K}_{|I|}[\psi](f_{i_1}, \ldots, f_{i_{|I|}}),$$

*where $i_1 < \cdots < i_{|I|}$ are the elements of $I$, and where the $\pi(\Pi)$ is the parity of the permutation associated with the partition.*

**Definition 2.5.** *(Gaussian Grassmann field). Let $(\mathcal{M}, \omega)$ be a NPS and let $\psi$ be a Grassmann field in it indexed by a complex separable Hilbert space $\mathcal{H}$, with conjugation $\Theta$. We say that $\psi$ is a (centred) Gaussian field, if*

$$\mathcal{K}_2[\psi](f_1, f_2) = \langle\!\langle f_1, G f_2 \rangle\!\rangle_\Theta, \qquad \mathcal{K}_n[\psi](f_1, \ldots, f_n) = 0, \quad \forall n \neq 2.$$

*for some $\Theta$-antisymmetric $G \in \mathcal{B}(\mathcal{H})$, see (2.1).*

**Remark 2.6.** Equivalently, one can define (centred) Gaussianity via the anticommuting Wick rule:

$$\omega(\psi(f_1) \cdots \psi(f_{2n})) = \text{Pf}((\langle\!\langle f_i, G f_j \rangle\!\rangle_\Theta)_{1 \leq i, j \leq 2n}), \qquad \omega(\psi(f_1) \cdots \psi(f_{2n+1})) = 0.$$

## 2.2 Processes and conditional expectations

In this paper, we are deeply inspired by some techniques coming from stochastic analysis. For this reason we need to consider Grassmann fields labelled by a real variable $t \in \mathbb{R}^+$ (which takes the place of Brownian motion or more general stochastic processes), a non-commutative conditional expectation and the related notion of martingality.

Let us begin with defining a filtration in a NPS.

**Definition 2.7.** *(Filtration). Let $(\mathcal{M}, \omega)$ be a NPS and let $I \subseteq \mathbb{R}$. A filtration over $I$ is a family $(\mathcal{M}_t)_{t \in I}$ of sub-$C^*$-algebras of $\mathcal{M}$ such that $\mathcal{M}_t \subseteq \mathcal{M}_{t'} \subseteq \mathcal{M}$ for any $t \leq t'$. The triple $(\mathcal{M}, \omega, (\mathcal{M}_t)_{t \in I})$ is called a filtered non-commutative probability space (FNPS).*

We then define processes and adapted processes as follows.

**Definition 2.8.** *(Process). Let $(\mathcal{M},\omega)$ be a NPS. A process is a family $(\psi_t)_{t\in I}$, of elements of $\mathcal{M}$, for some interval $I\subseteq\mathbb{R}$. Let $(\mathcal{M},\omega,(\mathcal{M}_t)_t)$ be a FNPS. A process $(\psi_t)_{t\in I}$ is called adapted (to the filtration $(\mathcal{M}_t)_t$) if $\psi_t\in\mathcal{M}_t$ for any $t\in I$.*

In particular, we are interested in Grassmann processes.

**Definition 2.9.** *(Grassmann process). Let $(\mathcal{M},\omega)$ be a NPS and $(\psi_t)_{t\in I}$ a process in it. We say that $(\psi_t)_{t\in I}$ is a Grassmann process (indexed by $\mathcal{H}$) if for any $s,t\in I$ $\psi_s$ and $\psi_t$ are compatible Grassmann fields (indexed by $\mathcal{H}$). Let $(\mathcal{M},\Omega,(\mathcal{M}_t)_t)$ be a FNPS. A Grassmann process $(\psi_t)_{t\in I}$ on it is called adapted if $\psi_t\in\mathcal{M}_t$ for any $t\in I$.*

Let us now introduce the conditional expectation with respect to a smaller $C^*$-algebra [Tom57, Tom58]. See also [Hia20] for a review and note that more general definitions have been provided, e.g., in [AC82].

**Definition 2.10.** *(Conditional expectation). Let $(\mathcal{M},\omega)$ be a NPS and $\mathcal{N}\subset\mathcal{M}$ a $C^*$-subalgebra. A linear mapping $\varphi\colon\mathcal{M}\to\mathcal{N}$ is a conditional expectation of $\mathcal{M}$ onto $\mathcal{N}$ with respect to $\omega$ if it is positive and satisfies*

$$\varphi(\mathbb{1}_\mathcal{M})=\mathbb{1}_\mathcal{N}, \qquad \varphi(N_1 M N_2)=N_1\varphi(M)N_2, \qquad \omega(M)=\omega(\varphi(M)),$$

*for any $N_1,N_2\in\mathcal{N}$ and $M\in\mathcal{M}$.*

*Let $(\mathcal{M},\omega,(\mathcal{M}_t)_{t\in I})$ be a FNPS. We say that $(\omega_t)_{t\in I}$ is a conditional expectation for $(\mathcal{M},\omega,(\mathcal{M}_t)_{t\in I})$ if for any $t\in I$, $\omega_t$ is a conditional expectation of $\mathcal{M}$ onto $\mathcal{M}_t$ with respect to $\omega$ and if the tower property is satisfied, that is, for any $A\in\mathcal{M}$ and for any $s\leqslant t$ we have*

$$\omega_s(\omega_t(A))=\omega_s(A). \tag{2.2}$$

**Remark 2.11.** The existence and uniqueness of a conditional expectation of a $C^*$-algebra onto a (generic) $C^*$-subalgebra with respect to a (generic) positive state are not guaranteed without additional assumption on the $C^*$-algebra, the state, and the $C^*$-subalgebra (see Appendix A for more details). For this reason, we need to explicitly include the tower law (2.2) in the definition of a conditional expectation for a FNPS. In Appendix A, we show that the FNPS we work on admit a conditional expectation as per Definition 2.10.

The following properties hold true, see [Kad04] for a proof.

**Proposition 2.12.** *Let $\mathcal{M}$ and $\mathcal{N}\subset\mathcal{M}$ be $C^*$-algebras and let $\varphi$ be a conditional expectation of $\mathcal{M}$ onto $\mathcal{N}$. Then, for any $M\in\mathcal{M}$:*

$$\varphi(M^*)=\varphi(M)^*, \qquad |\varphi(M)|^2\leqslant\varphi(|M|^2), \qquad \|\varphi(M)\|\leqslant\|M\|.$$

We conclude with the definition of a martingale.

**Definition 2.13.** *(Martingale). Let $(\mathcal{M},\omega,(\mathcal{M}_t)_{t\in I})$ be a FNPS and let $(\omega_t)_{t\in I}$ be the conditional expectation. We say that an adapted process $(\psi_t)_{t\in I}$ is a martingale iff*

$$\omega_s(\psi_t)=\psi_s, \qquad \forall s\leqslant t.$$

## 2.3 Grassmann Brownian martingales

We now discuss in detail Grassmann Brownian martingales (GBM). Our exposition is based on abstract definitions, yet one can prove the existence of a FNPS where these objects exist, see Appendix A.

We first of all define independence and independent increments.

**Definition 2.14.** *Let $(\mathcal{M},\omega)$ be a NPS. We say that the subalgebras $\mathcal{N}_1, \mathcal{N}_2 \subseteq \mathcal{M}$ are independent if for any $A_1 \in \mathcal{N}_1$ and $A_2 \in \mathcal{N}_2$ we have $\omega(A_1 A_2) = \omega(A_2 A_1) = \omega(A_1)\omega(A_2)$. Furthermore, let $(\mathcal{M},\omega,(\mathcal{M}_t)_t)$ be a FNPS and let $(\omega_t)_t$ be a conditional expectation for it. We say that an adapted process $(Y_t)_t$ has independent increments if for any $s \leqslant t$, the subalgebra $\mathcal{M}_t$ is independent of the $C^*$-algebra*

$$\mathcal{M}^Y_{>t} = \overline{\mathrm{span}\{Y_\tau - Y_t | \tau > t\}}$$

*with respect to $\omega_s$, namely for any $A_1 \in \mathcal{M}_t$ and $A_2 \in \mathcal{M}^Y_{>t}$ we have*

$$\omega_s(A_1 A_2) = \omega_s(A_2 A_1) = \omega_s(A_1)\omega_s(A_2) = \omega_s(A_2)\omega_s(A_1).$$

**Definition 2.15.** *(Grassmann Brownian martingale). Let $\mathfrak{h}$ be a Hilbert space with conjugation $\Theta$ and let $(G_t)_t \subset \mathcal{B}(\mathfrak{h})$ be $\Theta$-antisymmetric. A Grassmann Brownian martingale on $(\mathcal{M},\omega,(\mathcal{M}_t)_t)$ indexed by $\mathfrak{h}$ and with covariance $(G_t)_t$ is the adapted centred Gaussian Grassmann process $(X_t)_{t \in I}$ with covariance $(G_t)_t$ and with independent increments. In particular, $\omega_s((X_t - X_s)(f)) = 0$ for any $s \leqslant t$ and $f \in \mathfrak{h}$, and*

$$\omega(X_t(f) X_s(g)) = \langle\!\langle f, G_{t \wedge s} g \rangle\!\rangle_\Theta \qquad \forall s,t \in I, \forall f,g \in \mathfrak{h}. \tag{2.3}$$

Grassmann Brownian martingales are indeed martingales, as follows from independent increments and the fact that $(X_t)_t$ is adapted. On the other hand, polynomials of a GBM are not martingales in general. This is however the case if one consider Wick's products, defined as follows.

**Definition 2.16.** *(Wick's products). Wick's products of a GBM $(X_t)_t$ are defined recursively by setting $[\![\mathbb{1}]\!] := \mathbb{1}$ and*

$$[\![X_t(v_1) \cdots X_t(v_{n+1})]\!] := X_t(v_1) [\![X_t(v_2) \cdots X_t(v_{n+1})]\!]$$
$$- \sum_{j=1}^n (-)^j \omega(X_t(v_1) X_t(v_j)) [\![X_t(v_2) \cdots X_t(v_{j-1}) X_t(v_{j+1}) \cdots X_t(v_{n+1})]\!].$$

**Proposition 2.17.** *Wick's products of a GBM $(X_t)_t$ are martingales, that is, for $s \leqslant t$*

$$\omega_s([\![X_t(v_1) \cdots X_t(v_n)]\!]) = [\![X_s(v_1) \cdots X_s(v_n)]\!].$$

**Proof.** This can be directly checked via the inductive definition and by using that the GBM has independent increments. □

We henceforth write $X_{s,t} := X_t - X_s$ and $G_{s,t} := G_t - G_s$ for $s \leqslant t$. Clearly, $X_{s,t}$ is a Gaussian Grassmann field with covariance $G_{s,t}$. As it turns out, one can realise a GBM indexed by $\mathfrak{h}$ as a suitable Grassmann Gaussian field on $L^2(\mathbb{R};\mathfrak{h})$ for which the norm is controlled in a natural way. We encode this information in the following definition.

**Definition 2.18.** *(Norm-compatibility).* Let $\mathfrak{h}$ be a Hilbert space with conjugation $\Theta$ and let $(G_t)_t \subset \mathcal{B}(\mathfrak{h})$ be $\Theta$-antisymmetric and differentiable in t. Let $(X_t)_t$ be a GBM on $(\mathcal{M}, \omega, (\mathcal{M}_t)_t)$ indexed by $\mathfrak{h}$ and with covariance $(G_t)_t$. We say that $(X_t)_t$ is norm-compatible if, for any decomposition $\dot{G}_t = C_t^* C_t U_t$, with $U_t$ unitary, $(X_t)_t$ can be realised in such a way that

$$\|X_t(f)\|^2 \leq c_X \int_0^t \left( \|C_s U_s f\|_{\mathfrak{h}}^2 + \|C_s \Theta f\|_{\mathfrak{h}}^2 \right) ds, \qquad \forall t \in I, \forall f \in \mathfrak{h}, \tag{2.4}$$

*for some constant $c_X > 0$.*

**Remark 2.19.** In our setting, we will have $\dot{G}_t = C_t^* C_t U$ with $[C_t, U] = [C_t, \Theta] = 0$. In this case, (2.4) can be rewritten as $\|X_t(f)\|^2 \lesssim \int_0^t \|C_s f\|_{\mathfrak{h}}^2 ds$.

**Theorem 2.20.** *Let $\mathfrak{h}$ be a separable Hilbert space with conjugation $\Theta$ and let $(G_t)_t \subset \mathcal{B}(\mathfrak{h})$ be $\Theta$-antisymmetric differentiable operators. Then, there exists a FNPS $(\mathcal{M}, \omega, (\mathcal{M}_t)_t)$ where we can realise a norm-compatible GBM indexed by $\mathfrak{h}$ and with covariance $(G_t)_t$.*

The proof, provided in Appendix A, is by construction and follows the same ideas presented in [OS73, ABDG20].

# 3 The FBSDE approach

In this section, we derive the forward-backward stochastic differential equation (FBSDE) which realises the stochastic quantisation in (1.5) for processes indexed by a finite-dimensional Hilbert space – which is our setting with the cutoffs $L \in \mathbb{N}$, $\varepsilon \in 2^{-\mathbb{N}}$. The derivation of the FBSDE is ultimately a consequence of Itô's lemma for Grassmann processes.

## 3.1 Finite-dimensional SDEs

Let $\mathfrak{h}$ be a finite-dimensional complex Hilbert space with scalar product $\langle \cdot, \cdot \rangle_{\mathfrak{h}}$ and with conjugation $\Theta$. Because $\mathfrak{h}$ is finite dimensional, all elements of the Grassmann algebra $\bigwedge \mathfrak{h}$ are polynomials. If $(v_j)_{j \in S}$ is an orthonormal basis for $\mathfrak{h}$, where $S$ is some ordered set, we define a linear basis for $\bigwedge \mathfrak{h}$ of monomials $(v_\alpha)_{\alpha \subseteq S}$, with $v_\alpha := \bigwedge_{j \in \alpha} v_j$ for all ordered subsets $\alpha \subseteq S$. We also set $\bigwedge_{\text{odd}} \mathfrak{h} := \bigoplus_{n \geq 0} \mathfrak{h}^{\wedge(2n+1)}$ and $\bigwedge_{\text{even}} \mathfrak{h} := \bigoplus_{n \geq 0} \mathfrak{h}^{\wedge 2n}$ and say that $F \in \bigwedge \mathfrak{h}$ is odd (resp. even) if in particular $F \in \bigwedge_{\text{odd}} \mathfrak{h}$ (resp. $F \in \bigwedge_{\text{odd}} \mathfrak{h}$). Finally, we say that $F$ is purely Grassmann (or nilpotent) if $F \in \bigoplus_{n \geq 1} \mathfrak{h}^{\wedge n}$.

Let $F$ be purely Grassmann and let $f$ be sufficiently smooth around $F_0 \in \mathbb{C}$. We can always define $f(F_0 + F) \in \bigwedge \mathfrak{h}$ via the power series expansion, $f(F_0 + F) := \sum_{n \geq 0} \frac{f^{(n)}(F_0)}{n!} F^n$, which is well-defined because it is a finite sum due to the nilpotency of $F$. Typical examples that will be used in the paper are the (principal value) logarithm and the exponential function: for any $F$ purely Grassmann and any $F_0 \in \mathbb{C} \setminus (-\infty, 0]$ we have $\log(F_0 + F) := \log(F_0) + \sum_{n \geq 1} (-1)^{n+1} \frac{F_0^{-n} F^n}{n}$ and, furthermore, for any $F \in \bigwedge \mathfrak{h}$ we have $\exp(F) := \sum_{n \geq 0} \frac{F^n}{n!}$. If $\psi$ is a Grassmann field indexed by $\mathfrak{h}$ and $F \in \bigwedge \mathfrak{h}$, we will write $F(\psi) := \psi(F)$, see Remark 2.2, and interpret the left-hand side as a function of $\psi$.

**Definition 3.1.** *(Functional derivative).* If $U \in \bigwedge \mathfrak{h}$ then $DU: \mathfrak{h} \to \bigwedge \mathfrak{h}$

$$(DU)(v) := \sum_i \langle\!\langle v, v_i \rangle\!\rangle_\Theta \partial_{v_i} U,$$

where $\partial_{v_j}$ is the anticommuting derivative, that is, the linear operator satisfying $\partial_{v_j} 1 = 0$ and $\partial_{v_j} v_i F + v_i \partial_{v_j} F = \delta_{i,j} F$ for any $F \in \bigwedge \mathfrak{h}$. We can also define the functional Laplacian $\mathrm{D}^2_C \colon \bigwedge \mathfrak{h} \to \bigwedge \mathfrak{h}$ associated with $C \in \mathcal{B}(\mathfrak{h})$ $\Theta$-antisymmetric

$$\mathrm{D}^2_C := \sum_{i,j} \langle\!\langle v_i, C v_j \rangle\!\rangle_\Theta \, \partial_{v_i} \partial_{v_j},$$

and, if $U \in \bigwedge_{\mathrm{even}} \mathfrak{h}$ the quadratic form

$$(\mathrm{D} U)^2_C := \sum_{i,j} \langle\!\langle v_i, C v_j \rangle\!\rangle_\Theta \, (\partial_{v_i} U)(\partial_{v_j} U).$$

**Remark 3.2.** Note that our definitions do not depend on the choice of orthonormal basis $(v_j)_{j \in S}$. The functional derivative $\mathrm{D} U$ is linear in its argument and obviously $(\mathrm{D} U)(\Theta v_j) = \partial_{v_j} U$. Finally, note that we require $\langle\!\langle v_i, C v_j \rangle\!\rangle_\Theta$ to be antisymmetric because $\partial_{v_i} \partial_{v_j} = -\partial_{v_j} \partial_{v_i}$.

If $F, F' \colon \mathfrak{h} \to \bigwedge \mathfrak{h}$ are both linear or anti-linear we also use the notation:

$$\langle F, F' \rangle := \sum_j F(\Theta v_j) F'(v_j), \tag{3.1}$$

whereas, if one is linear and the other anti-linear we write

$$F \cdot F' := \sum_j F(v_j) F'(v_j). \tag{3.2}$$

The expression above do not depend on the choice of the orthonormal basis $(v_j)_{j \in S}$ and, of course, coincide for $\Theta$-real orthonormal bases. By abuse of notation, we extend (3.2) to $\langle \psi, \mathrm{D} U \rangle(\varphi) := \sum_j \psi(\Theta v_j)((\mathrm{D} U)(v_j))(\varphi)$, with $\psi$ and $\varphi$ Grassmann fields indexed by $\mathfrak{h}$.

Let us now discuss in more detail the type of SDEs we are going to consider.

**Definition 3.3.** *A drift $F$ on $\bigwedge \mathfrak{h}$ is a linear mapping $F \colon \mathfrak{h} \to \bigwedge_{\mathrm{odd}} \mathfrak{h}$. We call a family of drifts $(F_t)_{t \in I}$ with $I \subseteq \mathbb{R}$ admissible if it is $L^\infty_{\mathrm{loc}}(I)$, that is, if $F_t = \sum_\alpha (F_t)_\alpha v_\alpha$, then the coefficients $(F_t)_\alpha$ are in $L^\infty_{\mathrm{loc}}(I)$.*

**Remark 3.4.** Requiring that $F$ takes values in $\bigwedge_{\mathrm{odd}} \mathfrak{h}$ is natural in the Grassmann context because we want the solution of the SDE to be anticommuting, see below.

Let $(F_t)_t$ be an admissible family of drifts. We consider the additive SDE

$$\begin{cases} \mathrm{d} \Psi_t = F_t(\Psi_t) \mathrm{d} t + \mathrm{d} X_t, & t \in [t_0, T], \\ \Psi_{t_0} = \varphi, \end{cases} \tag{3.3}$$

where $(X_t)_t$ is a norm-compatible GBM with continuously differentiable covariance $(G_t)_t$ and where the Grassmann field $\varphi \colon \mathfrak{h} \to \mathcal{M}_{t_0}$ is compatible with $X_t$ for any $t \in [t_0, T]$, recall Definition 2.3. More precisely, (3.3) is a shorthand for

$$\Psi_t(f) = \varphi(f) + \int_{t_0}^t (F_s(\Psi_s))(f) \mathrm{d} s + X_{t_0, t}(f), \qquad t \in [t_0, T], \ \forall f \in \mathfrak{h}. \tag{3.4}$$

A similar equation with time-independent drift and Brownian motion in place of the GBM was already studied in [ABDG20], where the local existence and uniqueness were proven. The same result holds for (3.3), under mild assumptions on the time-dependent drift.

**Proposition 3.5.** *(Local existence and uniqueness). Let $(X_t)_t$ be a norm-compatible GBM indexed by a finite-dimensional Hilbert space $\mathfrak{h}$ with covariance $(G_t)_t$, let $(F_t)_t$ be an admissible family of drifts and let $\varphi: \mathfrak{h} \to \mathcal{M}_{t_0}$ be a Grassmann field compatible with $X_t$ for any $t \geq t_0$. Then, for some $T > t_0$, depending on $F_t$ and $G_t$ there exists a unique adapted Grassmann process $\Psi_t$ solving (3.4) for $t \in [t_0, T]$ and compatible with $\varphi$ and with $X_t$ for any $t \geq t_0$.*

**Proof.** Local existence and uniqueness is proven by a Banach fixed-point argument in the space $C([t_0, T], \mathcal{M})$ for suitable $T$, see [ABDG20] for details. In particular, the solution can be constructed by a Picard iteration. Because the $n$-th Picard iteration is an adapted Grassmann process compatible with $\varphi$ and with $X_t$ for any $t \geq t_0$, the thesis follows. □

**Remark 3.6.** To obtain global existence, one needs more information on $F_t$ and $G_t$. This will be done in the next section. In the case of $F_t$ constant in $t$, a simple argument for global existence can be found in [ABDG20].

The next lemma presents an averaged version of Itô's formula for Grassmann variables. A similar statement was proven in [ABDG20] in the special case of SDE with additive Grassmann Brownian motion.

**Lemma 3.7.** *(Itô's lemma). Let $(P_t)_t \subset \bigwedge \mathfrak{h}$ be a continuously differentiable family and let $\Psi$ solve (3.4) (same setting) over some interval $[t_0, T]$. Then, for $s, t \in [t_0, T]$ we have:*

$$\omega_s(P_t(\Psi_t) - P_s(\Psi_s)) = \int_s^t \omega_s \left[ (\partial_r P_r)(\Psi_r) + \frac{1}{2}(D^2_{G_r} P_r)(\Psi_r) + \langle F_r, D P_r \rangle (\Psi_r) \right] dr \qquad s \leq t.$$

**Proof.** Here we follow a standard argument, see also the proof of [ABDG20, Theorem 32]. The fundamental step in the proof is to show that for any $t_1, t_2 \in [s, t]$ $t_1 \leq t_2$ we have

$$\left\| \omega_s \left[ P_{t_2}(\Psi_{t_2}) - P_{t_1}(\Psi_{t_1}) - (t_2 - t_1) \left( (\partial_t P_{t_1})(\Psi_{t_1}) + \frac{1}{2}(D^2_{G_{t_1}} P_{t_1})(\Psi_{t_1}) \right) \right. \right.$$
$$\left. \left. - (t_2 - t_1) \langle F_{t_1}(\Psi_{t_1}), D P_{t_1}(\Psi_{t_1}) \rangle \right] \right\| = o(|t_2 - t_1|), \tag{3.5}$$

where we recall that $\Psi_{s,t} := \Psi_t - \Psi_s$, and the estimate is uniform in $|t_2 - t_1|$ small enough. We write

$$\begin{aligned} P_{t_2}(\Psi_{t_2}) - P_{t_1}(\Psi_{t_1}) &= P_{t_2}(\Psi_{t_2}) - P_{t_2}(\Psi_{t_1}) + P_{t_2}(\Psi_{t_1}) - P_{t_1}(\Psi_{t_1}) \\ &= P_{t_2}(\Psi_{t_2}) - P_{t_2}(\Psi_{t_1}) + (t_2 - t_1)(\partial_t P_{t_1})(\Psi_{t_1}) + o(|t_2 - t_1|), \end{aligned} \tag{3.6}$$

where we used the differentiability of $P_t$ with respect to $t$. Then, by using that $\Psi_s$ is a Grassmann process, that is, that it satisfies exact anticommuting relations at all times, we write

$$P_{t_2}(\Psi_{t_2}) - P_{t_2}(\Psi_{t_1}) = ((e^{\langle \Psi_{t_1, t_2}, D \rangle} - \mathbb{1}) P_{t_2})(\Psi_{t_1}), \tag{3.7}$$

where the exponential is intended as a (truncated) power series, see below (3.1) for the notation $\langle \Psi_{t_1, t_2}, D \rangle$. To bound the right-hand side, we first of all control $\Psi_{s,t}$:

$$\begin{aligned} &\sup_j \| \Psi_{t_1, t_2}(v_j) \| \\ &\leq (t_2 - t_1) \sup_j \sup_{r \in [s,t]} \|(F_r(\Psi_r))(v_j)\| + \sup_j \|X_{t_1, t_2}(v_j)\| \\ &\lesssim (t_2 - t_1) \left( 1 + \sup_j \sup_{r \in [s,t]} \|\Psi_r(v_j)\| \right)^{\deg(F)} \sum_\alpha \|F_\alpha\|_{L^\infty(I)} + \sup_j \left( \int_{t_1}^{t_2} \|C_r v_j\|_{\mathfrak{h}}^2 dr \right)^{\frac{1}{2}} \\ &\lesssim |t_2 - t_1|^{1/2}, \end{aligned} \tag{3.8}$$

where we used the local existence of $\Psi_r$ and that $X_t$ is norm-compatible ($\dot{G}_r = C_r^* C_r U_r$). By the Lagrange remainder formula, this implies

$$\sup_j \left\| \left( \left[ e^{\langle \Psi_{t_1,t_2}, D \rangle} - \mathbb{1} - \langle \Psi_{t_1,t_2}, D \rangle - \frac{1}{2} \langle \Psi_{t_1,t_2}, D \rangle^2 \right] P_{t_2} \right) (\Psi_{t_1}) \right) (v_j) \right\| \\ \lesssim \left( \sup_j \| \Psi_{t_1,t_2}(v_j) \| \right)^3 := O(|t_2 - t_1|^{3/2}). \tag{3.9}$$

Furthermore, we have that

$$\sup_j \| \Psi_{t_1,t_2}(v_j) - (t_2 - t_1)(F(\Psi_{t_1}))(v_j) - X_{t_1,t_2}(v_j) \| \lesssim |t_2 - t_1|^{\frac{3}{2}},$$

from which we control the first term in the series in (3.7) as follows

$$\begin{aligned} \omega_s((\langle \Psi_{t_1,t_2}, DP_{t_2} \rangle)(\Psi_{t_1})) &= (t_2 - t_1)\, \omega_s(\langle F_{t_1}(\Psi_{t_1}), DP_{t_2} \rangle(\Psi_{t_1})) \\ &\quad + \omega_s((\langle X_{t_1,t_2}, DP_{t_2} \rangle)(\Psi_{t_1})) + O\left(|t_2 - t_1|^{\frac{3}{2}}\right) \\ &= (t_2 - t_1)\, \omega_s(\langle F_{t_1}(\Psi_{t_1}), DP_{t_1} \rangle)(\Psi_{t_1}) \\ &\quad + \langle \omega_s(X_{t_1,t_2}), \omega_s(DP_{t_1}(\Psi_{t_1})) \rangle + O\left(|t_2 - t_1|^{\frac{3}{2}}\right) \\ &= (t_2 - t_1)\, \omega_s(\langle F_{t_1}(\Psi_{t_1}), DP_{t_1} \rangle(\Psi_{t_1})) + O\left(|t_2 - t_1|^{\frac{3}{2}}\right) \end{aligned} \tag{3.10}$$

where we use that $P_t$ is continuously differentiable (in $t$) and $X_t$ has independent increments (i.e. for any $K \in \mathcal{M}_{t_1}$ and $H \in \bigwedge \mathfrak{h}$ we have $\omega_s(KH(X_{t_1,t_2})) = \omega_s(K)\, \omega_s(H(X_{t_1,t_2}))$, see Definition 2.14).

Regarding the second term in (3.7), we have instead:

$$\begin{aligned} \omega_s((\langle \Psi_{t_1,t_2}, D \rangle^2 P_{t_2})(\Psi_{t_1})) &= \omega_s((\langle X_{t_1,t_2}, D \rangle^2 P_{t_2})(\Psi_{t_1})) + o(|t_2 - t_1|) \\ &= \sum_{i,j} \omega_s((X_{t_1,t_2}(v_i) X_{t_1,t_2}(v_j)) \partial_{v_i} \partial_{v_j} P_{t_2}(\Psi_{t_1})) + o(|t_1 - t_2|) \\ &= (t_2 - t_1)\, \omega_s((D^2_{\dot{G}_s} P_{t_2})(\Psi_{t_1})) + o(|t_2 - t_1|) \\ &= (t_2 - t_1)\, \omega_s((D^2_{\dot{G}_s} P_{t_1})(\Psi_{t_1})) + o(|t_2 - t_1|), \end{aligned} \tag{3.11}$$

where we used that $P_t$ is continuously differentiable (in $t$), $X_t$ has independent increments, and that

$$\sup_j \| \Psi_{t_1,t_2}(v_j) - X_{t_1,t_2}(v_j) \| = O(|t_1 - t_2|).$$

Putting together (3.6), (3.9), (3.10) and (3.11) we get (3.5). Once we have (3.5), by observing that

$$\omega_s(P_t(\Psi_t) - P_s(\Psi_s)) = \lim_{|\pi_n| \to 0} \sum_{t_j \in \pi_n} \omega_s(P_{t_j}(\Psi_{t_j}) - P_{t_{j-1}}(\Psi_{t_{j-1}})),$$

where $\pi_n$ is a sequence of partition of $[s, t]$ whose diameter is going to 0, we can proceed as in the proof of [ABDG20, Theorem 32] obtaining the thesis. □

## 3.2 Flow under the conditional expectation

As a simple consequence of Itô's lemma, we establish the following flow equations.

**Proposition 3.8.** *Let $U \in \bigwedge_{\text{even}} \mathfrak{h}$. Let $(X_t)_t$ be a GBM with covariance $(G_t)_t \subset \mathcal{B}(\mathfrak{h})$. Define*

$$U_t(\psi) := \omega_t(U(\psi + X_{t,T})) \qquad \text{for } t \in [0, T], \tag{3.12}$$

where $\psi: \mathfrak{h} \to \mathcal{M}_t$ is a Grassmann field compatible with $X_{t,T}$, see Definition 2.3. Then, $U_t \in \bigwedge_{\text{even}} \mathfrak{h}$ solves Kolmogorov's equation

$$\partial_t U_t + \frac{1}{2} D^2_{\dot{G}_t} U_t = 0. \tag{3.13}$$

Furthermore, as long as $U_t(0) \neq 0$, set $F_t(\psi) := D(\log U_t(\psi))$. Then, $F_t$ solves a HJB-type equation

$$\partial_t F_t + \frac{1}{2} D^2_{\dot{G}_t} F_t + DF_t \cdot \dot{G}_t F_t = 0. \tag{3.14}$$

**Proof.** Eq. (3.13) is a straightforward consequence of Itô's lemma. Because the equation is linear and the boundary datum $U$ even, it follows that the solution $U_t$ is likewise even. As long as $U_t(0) \neq 0$, $V_t := \log U_t$ is well-defined. Therefore, plugging $U_t = e^{V_t}$ into (3.13), using that $U_t$ is even and that $U_t(0) \neq 0$ proves that $V_t$ satisfies the following HJB-type equation:

$$\partial_t V_t + \frac{1}{2} D^2_{\dot{G}_t} V_t + \frac{1}{2} (DV_t)^2_{\dot{G}_t} = 0.$$

Taking the functional derivative proves that $F_t$ satisfies (3.14). $\square$

**Remark 3.9.** Note that since Wick's products are martingales, they satisfy Kolmogorov's equation, that is, if one defines $H^n_t \in \bigwedge \mathfrak{h}$ by $H^n_t(X_t) := [\![ X_t(f_1) \cdots X_t(f_n) ]\!]$, then $H^n_t$ solves eq. (3.13).

We can finally prove a stochastic quantisation formula. Our first result connects the interacting measure with the solution of an SDE where the drift $(F_t)_t$ is given by the solution of a HJB-type equation.

**Proposition 3.10.** Let $0 \leq t \leq T$ and let $V_T \in \bigwedge_{\text{even}} \mathfrak{h}$, $\mathfrak{h}$ a finite-dimensional Hilbert space with conjugation $\Theta$. Let $(X_t)_t$ be a norm-compatible GBM with differentiable covariance $(G_t)_t$, assume that

$$\omega(e^{V_T(X_{s,T})}) \neq 0 \qquad \forall s \in [t, T],$$

and set $V_t(\varphi) := \log \omega_t(e^{V_T(\varphi + X_{t,T})})$, compare with Theorem 3.8. Assume that $\Psi_s$ is the solution of the following SDE on $[t, T]$

$$d\Psi_s = -\dot{G}_s DV_s(\Psi_s) ds + dX_s, \qquad \Psi_t = 0. \tag{3.15}$$

Then, for any $P \in \bigwedge \mathfrak{h}$ we have

$$\omega_t(P(\Psi_s)) = \frac{\omega_t(P(X_{t,s}) e^{V_T(X_{t,T})})}{\omega_t(e^{V_T(X_{t,T})})}, \qquad \forall s \in [t, T]. \tag{3.16}$$

**Remark 3.11.** In (3.15) and henceforth in such SDEs $\dot{G}_s DV_s: \mathfrak{h} \to \bigwedge \mathfrak{h}$ is an abuse of notation for $-DV_s \cdot \dot{G}_s \cdot$, that is, it is linear and intended as

$$(\dot{G}_s DV_s)(v) = \sum_i \langle\!\langle v, \dot{G}_r v_i \rangle\!\rangle_\Theta \partial_{v_i} V_r. \tag{3.17}$$

Therefore, the SDE in (3.15) is the shorthand for

$$\Psi_s(v) = -\sum_i \int_t^s \langle\!\langle v, \dot{G}_r v_i \rangle\!\rangle_\Theta (\partial_{v_i} V_r)(\Psi_r) dr + X_s(v), \qquad \forall v \in \mathfrak{h}.$$

**Proof.** Let $U := \exp(V_T)$ and set $U_s(\varphi) := \omega_s(U(\varphi + X_{s,T}))$ for some Grassmann field $\varphi: \mathfrak{h} \to \mathcal{M}_s$, see (3.12). Define the following functions for $\alpha \subseteq S$ and for $s \in [t, T]$:

$$P_\alpha(s) := \omega_t(v_\alpha(\Psi_s)), \qquad \tilde{P}_\alpha(s) := \frac{\omega_t(v_\alpha(X_{t,s})\, e^{V_T(X_{t,T})})}{U_t(0)}$$

where $(v_\alpha)_{\alpha \subseteq S}$ is our choice of linear basis of the Grassmann algebra $\bigwedge \mathfrak{h}$.

By Itô's lemma and by the fact that $\Psi$ solves the SDE (3.15), we find that the family $(P_\alpha(s))_\alpha$ satisfies the following system of ODEs:

$$\begin{aligned}
dP_\alpha(s) &= \omega_t\Big(\frac{1}{2}(D^2_{\dot{G}_s} v_\alpha)(\Psi_s) - \langle \dot{G}_s DV_s, Dv_\alpha \rangle(\Psi_s)\Big) ds \\
&= \omega_t\Big(\frac{1}{2}(D^2_{\dot{G}_s} v_\alpha)(\Psi_s) + \sum_{i,j} \langle\!\langle v_i, \dot{G}_s v_j\rangle\!\rangle_\Theta ((\partial_{v_i} V_s)(\partial_{v_j} v_\alpha))(\Psi_s)\Big) ds \\
&= \Big(\frac{1}{2} D^2_{\dot{G}_s} P_\alpha(s) + \sum_\beta \sum_{i,j} \langle\!\langle v_i, \dot{G}_s v_j\rangle\!\rangle_\Theta (V_s)_\beta \operatorname{sign}(i,j,\alpha,\beta)\, P_{\gamma(i,j,\alpha,\beta)}(s)\Big) ds,
\end{aligned}$$

with suitable $\operatorname{sign}(i,j,\alpha,\beta) \in \{\pm 1\}$ and $\gamma(i,j,\alpha,\beta) \subset S$, where in the first line $\dot{G}_s DV_s$ is intended as in (3.17) and where in the second line we used that $\dot{G}_s$ is $\Theta$-antisymmetric. Concerning $\tilde{P}_\alpha(s)$ we note that by the tower property, by adaptedness, and by the definition of $U$, we have

$$\omega_t(v_\alpha(X_{t,s})\, e^{V_T(X_{t,T})}) = \omega_t(v_\alpha(X_{t,s})\, \omega_s(e^{V_T(X_{t,T})})) = \omega_t((v_\alpha U_s)(X_{t,s})).$$

Therefore, by Itô's lemma and by the fact that $U$ solves Kolmogorov's equation, we find that $(\tilde{P}_\alpha(s))_\alpha$ satisfies the same ODEs as $(P_\alpha(s))_\alpha$:

$$\begin{aligned}
d\tilde{P}_\alpha(s) &= [U_t(0)]^{-1} \omega_t(d(v_\alpha U_s)(X_{t,s})) \\
&= [U_t(0)]^{-1} \omega_t\Big((v_\alpha \partial_s U_s)(X_{t,s}) + \frac{1}{2}(D^2_{\dot{G}_s}(v_\alpha U_s))(X_{t,s})\Big) ds \\
&= [U_t(0)]^{-1} \omega_t\Big(\frac{1}{2}((D^2_{\dot{G}_s} v_\alpha) U_s)(X_{t,s}) + \sum_{i,j} \langle\!\langle v_i, \dot{G}_s v_j\rangle\!\rangle_\Theta ((\partial_{v_i} V_s)(\partial_{v_j} v_\alpha) U_s)(X_{t,s})\Big) ds \\
&= \Big(\frac{1}{2} D^2_{\dot{G}_s} \tilde{P}_\alpha(s) + \sum_\beta \sum_{i,j} \langle\!\langle v_i, \dot{G}_s v_j\rangle\!\rangle_\Theta (V_s)_\beta \operatorname{sign}(i,j,\alpha,\beta)\, \tilde{P}_{\gamma(i,j,\alpha,\beta)}(s)\Big) ds,
\end{aligned}$$

where we used that $U_s$ is even ($V_T$ being even) and where $\operatorname{sign}(i,j,\alpha,\beta)$ and $\gamma(i,j,\alpha,\beta)$ are as above. The existence of $\Psi_s$ over $[t, T]$ implies the existence of $P_\alpha(s)$ over the same interval and thus, since $P_\alpha(s) = \tilde{P}_\alpha(s)$ for $s = t$, we actually have the equality for any $s \in [t, T]$. $\square$

Our ultimate goal is the characterisation of the interacting measure by means of the solution of an SDE that does not require us to solve the HJB equation. It suffices to "undo" the HJB flow of the drift $(DV_t)_t$ and thus study the problem by means of a FBSDE.

**Theorem 3.12.** *(FBSDE). Let $0 \leq t \leq T$ and let $V_T \in \bigwedge_{\text{even}} \mathfrak{h}$, $\mathfrak{h}$ a finite-dimensional Hilbert space with conjugation $\Theta$. Let $(X_t)_t$ be a norm-compatible GBM with differentiable covariance $(G_t)_t$. Assume that $\Psi_s$ solves the following equation for $s \in [t, T]$:*

$$d\Psi_s = -\dot{G}_s\, \omega_s(DV_T(\Psi_T))\, ds + dX_s, \qquad \Psi_t = 0. \tag{3.18}$$

*Then, the following equation holds true for any $P \in \bigwedge \mathfrak{h}$:*

$$\omega_t(P(\Psi_s)) = \frac{\omega_t(P(X_{t,s}) e^{V_T(X_{t,T})})}{\omega_t(e^{V_T(X_{t,T})})}, \qquad \forall s \in [t, T]. \tag{3.19}$$

**Remark 3.13.** As should be clear, the crucial difference between the statements in Proposition 3.10 and Theorem 3.12 is that the latter relies only on the FBSDE and its solution, whereas the former is formulated in terms of $DV_t$, $V_t$ solving the HJB equation. This shift in perspective is one of the key new ideas in this paper, and will allow us to truncate the flow equation and rely on the solution theory for the FBSDE to estimate the truncation error.

**Proof.** Let $U$ be as in the proof of Proposition 3.10. Note that $U_T(0) = e^{V_T(0)} \neq 0$, because $V_T(0) \in \mathbb{C}$. Since $U_s(0)$ is continuous in $s$, we therefore find a $\bar{t}$ such that $U_s(0) \neq 0$ for any $s \in (\bar{t}, T]$ and that if $\bar{t} > -\infty$ then $U_{\bar{t}}(0) = 0$. Accordingly, $F_s := DV_s$ is well-defined and solves the HJB-type equation in (3.14) with $F_T = DV_T$ by Proposition 3.8. Then, by Itô's lemma for $s \in (\bar{t} \vee t, T]$ we see that $(F_s(\Psi_s))_s$ is a martingale, and in particular

$$\omega_s(DV_T(\Psi_T)) = F_s(\Psi_s).$$

If $\bar{t} \leq t$, the solution of (3.18) solves also the SDE in (3.15) and this concludes the proof. Suppose otherwise that $\bar{t} > t$ and let $\tau \in (\bar{t}, T]$. Denoting by $\Psi^{(\tau)}$ the solution of (3.18) on $[\tau, T]$, with initial datum $\Psi^{(\tau)}_\tau = 0$ we clearly have that $\Psi^{(\tau)}$ likewise solves (3.15) on $[\tau, T]$ and therefore by Proposition 3.10 we have that

$$\omega_\tau(P(\Psi^{(\tau)}_s)) = \frac{\omega_\tau(P(X_{\tau,s}) e^{V_T(X_{\tau,T})})}{\omega_\tau(e^{V_T(X_{\tau,T})})}, \qquad \forall s \in [\tau, T]. \tag{3.20}$$

Let us plug $P = e^{-V_T}$ into eq. (3.20) with $s = T$, obtaining

$$\omega_\tau\big(e^{-V_T(\Psi^{(\tau)}_T)}\big) \omega_\tau(e^{V_T(X_{\tau,T})}) = 1, \qquad \bar{t} < \tau \leq T. \tag{3.21}$$

Because $\bar{t} > t$, we have that $\sup_{\tau \in (\bar{t}, T]} \|\Psi^{(\tau)}_T\|$ is bounded and thus by eq. (3.21) and by Proposition 2.12 we have

$$\inf_{\tau \in (\bar{t}, T]} |U_\tau(0)| = \inf_{\tau \in (\bar{t}, T]} |\omega_\tau(e^{V_T(X_{\tau,T})})| \geq e^{-(\sup_{\tau \in (\bar{t}, T]} \|\Psi^{(\tau)}_T\|)^{\deg(V_T)}} > 0.$$

This contradicts the definition of $\bar{t}$, therefore $\bar{t} \leq t$, and the claim is proven. □

**Lemma 3.14.** *Under the same assumptions of Theorem 3.12, $\Psi$ is a solution of the FBSDE in (3.18) iff it is a solution of*

$$d\Psi_s = -\dot{G}_s(F_s(\Psi_s) + R_s) ds + dX_s, \qquad \Psi_0 = 0,$$

*where $(F_t)_{t \in [0,T]} \subset \bigwedge \mathfrak{h}$ is any continuously differentiable interpolating family such that $F_T = DV_T$ and where the remainder process $R_s$ solves the following self-consistent equation:*

$$R_s = \int_s^T \omega_s(\mathscr{H}_r[F_r](\Psi_r)) dr + \int_s^T \omega_s(DF_r(\Psi_r) \cdot \dot{G}_r R_r) dr, \tag{3.22}$$

*with*

$$\mathscr{H}_r[F_r] := \partial_r F_r + \frac{1}{2} D^2_{\dot{G}_r} F_r + DF_r \cdot \dot{G}_r F_r$$

**Proof.** Set $R_t := \omega_t(DV_T(\Psi_T)) - F_t(\Psi_t) = \omega_t(F_T(\Psi_T)) - F_t(\Psi_t)$. Then, by Itô's lemma, and by using that $\dot{G}_s$ is $\Theta$-antisymmetric we can write:

$$R_t = \int_t^T \omega_t\left(\partial_s F_s(\Psi_s) + \frac{1}{2}D^2_{\dot{G}_s} F_s(\Psi_s) + DF_s(\Psi_s)\cdot \dot{G}_s \omega_s(F_T(\Psi_T))\right)ds \qquad (3.23) \quad \square$$
$$= \int_t^T \omega_t(\mathcal{H}_s[F_s])ds + \int_t^T \omega_t(DF_s(\Psi_s)\cdot \dot{G}_s R_s)ds.$$

**Remark 3.15.** Wick-ordering renormalisation corresponds to choosing $F_s$ as the solution of

$$\partial_s F_s(\Psi_s) + \frac{1}{2}D^2_{\dot{G}_s} F_s(\Psi_s) = 0,$$

since in fact $F_s$ is a martingale, see Proposition 2.17 and Remark 3.9. In this case, controlling $F_s$ is simple but the non-linear term $DF_r \cdot \dot{G}_r F_r$ appears in its entirety in the equation for $R_t$. To go beyond the Wick-ordering renormalisation, we will partially include the non-linearity into the flow equation for $F_s$, see Section 4.

## 3.3 FBSDE for the model

We let $\mathbb{N}_\infty := \mathbb{N} \cup \{\infty\}$ and $\mathbb{T}^d_{L,\varepsilon} := ((\varepsilon\mathbb{Z})/(L\mathbb{Z}))^d$ be the toroidal lattice of size $L \in \mathbb{N}_\infty$ and lattice spacing $\varepsilon \in 2^{-\mathbb{N}_\infty}$ with the understanding that $\mathbb{T}^d_{L,0} := \mathbb{T}^d_L$ and $\mathbb{T}^d_{\infty,0} := \mathbb{R}^d$. Because we deal with spin-$1/2$ fermions, for $L \in \mathbb{N}_\infty$ and $\varepsilon \in 2^{-\mathbb{N}_\infty}$ we introduce the Hilbert space

$$\mathfrak{h}_{L,\varepsilon} := L^2(\mathbb{T}^d_{L,\varepsilon}; \mathbb{C}^2) \oplus L^2(\mathbb{T}^d_{L,\varepsilon}; \mathbb{C}^2),$$

where for $\varepsilon > 0$ the scalar product on $L^2(\mathbb{T}^d_{L,\varepsilon}; \mathbb{C}^2)$ is given by

$$\langle f, g\rangle_{L^2(\mathbb{T}^d_{L,\varepsilon}; \mathbb{C}^2)} = \sum_{\sigma=\uparrow,\downarrow}\sum_{x\in\mathbb{T}^d_{L,\varepsilon}} \varepsilon^d \overline{f_\sigma(x)} g_\sigma(x) =: \int_{\mathbb{T}^d_{L,\varepsilon}} \overline{f(x)} \cdot g(x) dx.$$

Clearly, $\mathfrak{h}_{L,\varepsilon}$ are finite-dimensional Hilbert spaces as long as $L \in \mathbb{N}$ and $\varepsilon > 0$. We abridge $\mathfrak{h}_L := \mathfrak{h}_{L,0}$, $\mathfrak{h} := \mathfrak{h}_\infty$ and note that for any finite $L$ we can embed $\mathfrak{h}_{L,\varepsilon} \hookrightarrow \mathfrak{h}_L$ via truncation in the Fourier series. In more detail, we let $(\mathbb{T}^d_{L,\varepsilon})^* := \{k \in 2\pi L^{-1}\mathbb{Z}^d \mid \|k\|_\infty < \pi\varepsilon^{-1}\}$ and $\mathcal{F}_{L,\varepsilon}: \mathfrak{h}_{L,\varepsilon} \to \ell^2((\mathbb{T}^d_{L,\varepsilon})^*; \mathbb{C}^4)$ be the Fourier transform $(\mathcal{F}_{L,\varepsilon}f)(k) := \int_{\mathbb{T}^d_{L,\varepsilon}} e^{ikx} f(x) dx$. To any element $f \in \mathfrak{h}_{L,\varepsilon}$ we associate the following element of $\mathfrak{h}_L$,

$$(\mathcal{F}^{-1}_{L,0}\mathcal{F}_{L,\varepsilon}f)(x) := L^{-d}\sum_{k\in(\mathbb{T}^d_{L,\varepsilon})^*} e^{ikx}(\mathcal{F}_{L,\varepsilon}f)(k).$$

Furthermore, if we identify $\mathbb{T}^d_L$ with the box $(-L/2, L/2]^d \subset \mathbb{R}^d$, we obtain the embedding

$$\mathfrak{h}_{L,\varepsilon} \hookrightarrow \mathfrak{h}_L \hookrightarrow \mathfrak{h},$$

that is, an element of $\mathfrak{h}_L$ is identified with an element of $\mathfrak{h}$ with compact support on $\mathbb{T}^d_L$. This is important to keep in mind because, even though we need to index the GBMs using the finite-dimensional Hilbert spaces $\mathfrak{h}_{L,\varepsilon}$ with $L \in \mathbb{N}$ and $\varepsilon > 0$, we want to represent such GBMs in a sufficiently large FNPS that is directly linked to the infinite-dimensional space $\mathfrak{h}$.

On $\mathfrak{h}$ we define the conjugation $\Theta(v \oplus w) := (\bar{w} \oplus \bar{v})$ for any $v \oplus w \in \mathfrak{h}$. The operator $G$ in (1.1) is $\Theta$-antisymmetric and bounded for $\gamma \leq d/2$ but has an ultraviolet divergence when $\gamma \geq 0$, in the sense that the kernel is singular on the diagonal. We introduce the interpolation, see also [Ste70, Riv91]:

$$G_t := \frac{U}{\Gamma(d/2-\gamma)}\left(t\mathbf{1}_{0\leq t<1}\int_{2^{-2}}^\infty + \mathbf{1}_{t\geq 1}\int_{2^{-2t}}^\infty\right)\zeta^{-\gamma+d/2-1}e^{-\zeta(\mathbb{1}-\Delta)}d\zeta, \qquad (3.24)$$

where $U \coloneqq \mathbb{1} \oplus -\mathbb{1}$ and where $\Gamma$ is the Gamma function. It is easy to see that $G_\infty = G$ and that $G_t$ is $\Theta$-antisymmetric and bounded.

Although it is tempting to introduce a regularisation by replacing $\Delta$ with the lattice Laplacian, for technical reasons, we instead make the following choice. We let $\min_\varepsilon$ be a smooth non-decreasing function on $\mathbb{R}^+$ such that $\min_\varepsilon(a) = a$ on $[0, \varepsilon^{-1}]$ and $\min_\varepsilon(a) = \varepsilon^{-1}$ on $[\varepsilon^{-1}+1, \infty)$. In other words, $\min_\varepsilon(a)$ is a smoothening of $\varepsilon^{-1} \wedge a$.

**Definition 3.16.** *Let $L \in \mathbb{N}_\infty$, $\varepsilon \in 2^{-\mathbb{N}_\infty}$, $\delta \in (0,1)$ and let $\varphi_\delta \in C_c^\infty([0,\infty))$ be a $\delta^{-1}$-Gevrey function[3.1] such that $\mathrm{supp}(\varphi_\delta) = [0,2]$ and $\varphi_\delta([0,1]) = 1$. We define $G_t^{L,\varepsilon}$ as the operator on $\mathfrak{h}$ with kernel*

$$G_t^{L,\varepsilon}(y;x) \coloneqq \frac{U}{\Gamma(d/2-\gamma)} \Big( t \mathbf{1}_{0\leqslant t<1} \int_{2^{-2}}^\infty + \mathbf{1}_{t\geqslant 1} \int_{1/\min_\varepsilon(2^t)^2}^\infty \Big)$$
$$L^{-d} \sum_{k\in(\mathbb{T}_{L,\varepsilon}^d)^*} e^{ik(y-x)} \zeta^{-\gamma+d/2-1} e^{-\zeta(k^2+1)} \varphi_\delta(\varepsilon|k|) \,\mathrm{d}\zeta, \quad (3.25)$$

*for $y,x \in (-L/2, L/2]^d$ and $G_t^{L,\varepsilon}(y;x) = 0$ otherwise, where the sum should be intended as integral if $L = \infty$.*

**Remark 3.17.** Note the dependence of $G_t^{L,\varepsilon}$ on $\delta$ when $\varepsilon > 0$. We avoid writing this dependence explicitly to keep our notation as light as possible.

We abridge $G_t^\varepsilon \coloneqq G_t^{\infty,\varepsilon}$ and note that, because of translation invariance, the Poisson summation formula implies $G_t^{L,\varepsilon}(y;x) = \sum_{m\in\mathbb{Z}^d} G_t^\varepsilon(y+mL;x)$, for $y,x \in (-L/2, L/2]^d$. We also have

$$\dot{G}_t^{L,\varepsilon}(y;x) = \mathbf{1}_{0\leqslant t<1} G_1^{L,\varepsilon}(y;x)$$
$$+ c_{d,\gamma}\, \chi_{1,-\log_2 \varepsilon}(t)\, U \frac{2^{(2\gamma-d)s}}{L^d} \sum_{k\in\frac{2\pi}{L}\mathbb{Z}^d} e^{ik(y-x)} e^{-2^{-2t}(k^2+1)} \varphi(\varepsilon|k|) \quad (3.26)$$

where $\chi_{a,b}(t)$ is a positive smooth approximation to $\mathbf{1}_{a\leqslant t\leqslant b}$, such that $\chi_{a,b}(t) = 1$ for $a \leqslant t \leqslant b$ and $\chi_{a,b}(t) = 0$ for $t < a$ or $t > b + \eta$ for some small fixed $\eta > 0$ and where $c_{d,\gamma}$ is some constant depending on $d$ and $\gamma$. We can write $\dot{G}_t^{L,\varepsilon} = (\mathfrak{C}_t^{L,\varepsilon})^2 U$ and read off the explicit expression for $\mathfrak{C}_t^{L,\varepsilon}$ from (3.26). Note the commutation relations $[\mathfrak{C}_t^{L,\varepsilon}, U] = [\mathfrak{C}_t^{L,\varepsilon}, \Theta] = 0$. We can finally introduce the following $(L,\varepsilon)$-dependent family of GBMs.

**Definition 3.18.** *We let $(X_t^{L,\varepsilon})_t^{L\in\mathbb{N}_\infty, \varepsilon\in 2^{-\mathbb{N}_\infty}}$ be the anticommuting family of norm-compatible GBMs, see Definitions 2.15 and 2.18, such that*

$$\omega\big(X_t^{L,\varepsilon}(f) X_t^{L',\varepsilon'}(g)\big) = \Big\langle\!\!\Big\langle f, \int_0^t \mathfrak{C}_s^{L,\varepsilon} \mathfrak{C}_s^{L',\varepsilon'} U \,\mathrm{d}s\, g \Big\rangle\!\!\Big\rangle_\Theta \quad \forall f,g \in \mathfrak{h}.$$

**Remark 3.19.** In Appendix A, we prove that there exists a sufficiently large FNPS $(\mathcal{M}, \omega, (\mathcal{M}_t)_t)$ where such a family can be constructed. Note that we are requiring that $(X_t^{L,\varepsilon})_t$ at different $L$ and $\varepsilon$ are correlated and this is actually natural if they originate from colouring the same "Grassmann white noise", see Remark A.10. In particular, note that the Gaussian field $X_t^{L,\varepsilon} - X_t^{L',\varepsilon'}$ has covariance $\int_0^t \big(\mathfrak{C}_s^{L,\varepsilon} - \mathfrak{C}_s^{L',\varepsilon'}\big)^* \big(\mathfrak{C}_s^{L,\varepsilon} - \mathfrak{C}_s^{L',\varepsilon'}\big) U \,\mathrm{d}s$.

---

3.1. We thus have $\sup_{x\in\mathbb{R}^+} |\partial_\nu^n \varphi_\delta(x)| \leqslant C^{1+n}(n!)^{\delta^{-1}}$, see [GMR21, Rod93].

Even though $X_t^{L,\varepsilon}(f)$ can be defined for $f \in \mathfrak{h}$, in view of the application to the stochastic quantisation theorem, see Theorem 3.12, we think of them as indexed by $\mathfrak{h}_{L,\varepsilon} \subset \mathfrak{h}$, when $L \in \mathbb{N}$ and $\varepsilon > 0$. In other words, for fixed $L \in \mathbb{N}$ and $\varepsilon > 0$, we consider only the image of $\bigwedge \mathfrak{h}_{L,\varepsilon}$ under $X_t^{L,\varepsilon}$. In this regard, it is particularly convenient to switch to the Kronecker delta basis and therefore obtain what we call a field on $\mathbb{T}_{L,\varepsilon}^d$, that is, a map of the form $\psi \colon \mathbb{T}_{L,\varepsilon}^d \to \mathscr{M}^4$.

Define $\delta_{x;\uparrow}^{L,\varepsilon} \coloneqq (\delta_x^{L,\varepsilon}, 0) \in L^2(\mathbb{T}_L^d; \mathbb{C}^2)$ and $\delta_{x;\downarrow}^{L,\varepsilon} \coloneqq (0, \delta_x^{L,\varepsilon}) \in L^2(\mathbb{T}_L^d; \mathbb{C}^2)$, where $\delta_x^{L,\varepsilon}$ is the extension of the Kronecker delta to $\mathbb{T}_L^d$ via Fourier transform, i.e. $\delta_x^{L,\varepsilon} \coloneqq \mathcal{F}_{L,0}^{-1} \mathcal{F}_{L,\varepsilon} \delta_x$. The set $(\delta_{x;\sigma}^{L,\varepsilon})_{x \in \mathbb{T}_{L,\varepsilon}^d, \sigma=\uparrow,\downarrow}$ is an orthogonal basis for $L^2(\mathbb{T}_{L,\varepsilon}^d; \mathbb{C}^2) \hookrightarrow L^2(\mathbb{T}_L^d; \mathbb{C}^2)$ (it is not orthonormal though).

**Definition 3.20.** *If $(X_t^{L,\varepsilon})_t^{L \in \mathbb{N}_\infty, \varepsilon \in 2^{-\mathbb{N}_\infty}}$ are the GBMs of Definition 3.18, by abuse of notation we let $X_t^{L,\varepsilon} \colon \mathbb{R}^d \to \mathscr{M}^4$ denote the corresponding field on $\mathbb{R}^d$, $x \mapsto X_{t,x}^{L,\varepsilon} \coloneqq (X_{t,x,\mu}^{L,\varepsilon})_{\mu \in \{\uparrow,\downarrow\} \times \{\pm\}}$ where*

$$X_{t,x,(\sigma,+)}^{L,\varepsilon} \coloneqq X_t^{L,\varepsilon}(\delta_{x;\sigma}^{L,\varepsilon} \oplus 0), \qquad X_{t,x,(\sigma,-)}^{L,\varepsilon} \coloneqq X_t^{L,\varepsilon}(0 \oplus \delta_{x;\sigma}^{L,\varepsilon}), \qquad (3.27)$$

*for $x \in \mathbb{T}_L^d$ and $X_{t,x+mL}^{L,\varepsilon} \coloneqq X_{t,x}^{L,\varepsilon}$ for any $m \in \mathbb{Z}^d$.*

We introduce the following point-wise symmetric products for Grassmann fields. If $\psi^{(1)}, \psi^{(2)}$ and $\psi^{(3)}$ are fields we set:

$$(\psi^{(1)} \cdot \psi^{(2)})_x \coloneqq \frac{1}{2} \sum_{\sigma=\uparrow,\downarrow} \sum_{\pi \in S_2} \psi_{x,(\sigma,+)}^{(\pi(1))} \psi_{x,(\sigma,-)}^{(\pi(2))},$$

$$(\psi^{(1)} \cdot \psi^{(2)} \cdot \psi^{(3)})_{x,\mu} \coloneqq -\frac{1}{6} \operatorname{sign}(\mu) \sum_{\pi \in S_3} \psi_{x,\mu}^{(\pi(1))} (\psi^{(\pi(2))} \cdot \psi^{(\pi(3))})_x,$$

where $S_n$ denotes the set of permutations of $n$ elements, and where if $\mu = (\sigma, \rho) \in \{\uparrow, \downarrow\} \times \{\pm\}$ we set

$$\operatorname{sign}(\mu) \coloneqq \rho. \qquad (3.28)$$

We also abridge the notation by setting $(\psi_x)^2 \coloneqq (\psi \cdot \psi)_x$ and $(\psi_{x,\mu})^3 \coloneqq (\psi \cdot \psi \cdot \psi)_{x,\mu}$.

We are finally in a position to apply the FBSDE to the construction of the interacting measure. In the FNPS $(\mathscr{M}, \omega, (\mathscr{M}_t)_t)$ we are interested in computing the expectation of observables[3.2] of the form $P(\psi)$, for some $P \in \bigwedge B_{1,1,0-}^{\gamma^+}$, see (1.12) and Remark 2.2, so that $P(\psi) \in \mathscr{A}$. This expectation is obtained as the weak limit

$$\omega^V(P(\psi)) = \lim_{L \to \infty} \lim_{\varepsilon \to 0} \frac{\omega\left(P(X^{L,\varepsilon}) e^{V_\infty^{L,\varepsilon}(X^{L,\varepsilon})}\right)}{\omega\left(e^{V_\infty^{L,\varepsilon}(X^{L,\varepsilon})}\right)}$$

$$= \lim_{L \to \infty} \lim_{\varepsilon \to 0} \omega(P(\Psi^{L,\varepsilon}))$$

$$= \omega(P(\Psi)),$$

with $X^{L,\varepsilon} \coloneqq X_\infty^{L,\varepsilon}$, $\Psi^{L,\varepsilon} \coloneqq \Psi_\infty^{L,\varepsilon}$, $P(\Psi) = \lim_{L \to \infty} \lim_{\varepsilon \to 0} P(\Psi^{L,\varepsilon})$ in $\mathscr{M}$ (see Remark 4.26), where $X_t^{L,\varepsilon}$ are the GBMs introduced above, where the potential $V_T^{L,\varepsilon}$ reads

$$V_T^{L,\varepsilon}(X^{L,\varepsilon}) \coloneqq \int_{\mathbb{T}_{L,\varepsilon}^d} \left[\frac{\lambda}{2}[(X_x^{L,\varepsilon})^2]^2 + \mu_T^\varepsilon (X_x^{L,\varepsilon})^2\right] dx, \qquad (3.29)$$

for suitable constant $\mu_T^\varepsilon = \mu_T^\varepsilon(\lambda) \in \mathbb{R}$, and where the lattice field $\Psi_t^{L,\varepsilon}$ solves for $T = \infty$

$$d\Psi_t^{L,\varepsilon} = -\dot{G}_t^{L,\varepsilon} \omega_t(DV_T^{L,\varepsilon}(\Psi_T^{L,\varepsilon})) dt + dX_t^{L,\varepsilon}, \qquad \Psi_0^{L,\varepsilon} = 0. \qquad (3.30)$$

---

3.2. See also Lemma 4.29 for the actual applications.

Note that by (3.27), the drift term is

$$\mathrm{D}V_T^{L,\varepsilon}(\Psi_T^{L,\varepsilon})(y;\mu) = \lambda\,(\Psi_{T,y,\mu}^{L,\varepsilon})^3 - \mathrm{sign}(\mu)\,\mu_T^\varepsilon\,\Psi_{T,y,-\mu}^{L,\varepsilon}.$$

**Remark 3.21.** Let $\tilde{G}_t^{L,\varepsilon}$ denote the periodisation of $G_t^{L,\varepsilon}$, so that $\tilde{G}_t^{L,\varepsilon}(y+mL;x+m'L) = \tilde{G}_t^{L,\varepsilon}(y;x)$ for any $m,m' \in \mathbb{Z}^d$. Then, we can equivalently solve

$$\Psi_{t,x,\mu}^{L,\varepsilon} = -\sum_{\mu'}\int_0^t\int_{\mathbb{T}_L^d}\bigl(\dot{\tilde{G}}_s^{L,\varepsilon}(x;y)\bigr)_{\mu,\mu'}\,\omega_s(\mathrm{D}V_T^\varepsilon(\Psi_T^{L,\varepsilon})(y;-\mu'))\mathrm{d}y\mathrm{d}s + X_{t,x,\mu}^{L,\varepsilon}, \tag{3.31}$$

where $V_T^\varepsilon \equiv V_T^{\infty,\varepsilon}$ is as in (3.29) but with integral over the whole $\mathbb{R}^d$. This defines an extension of the solution to the whole $\mathbb{R}^d$, satisfying $\Psi_{t,x+mL}^{L,\varepsilon} = \Psi_{t,x}^{L,\varepsilon}$ for any $m \in \mathbb{Z}^d$.

Using Lemma 3.14, we solve the FBSDE by means of an interpolation scheme, for suitable $F_t^{L,\varepsilon}$. As it turns out, we can actually take it already in the infinite-volume case and likewise replace $\dot{\tilde{G}}_t^{L,\varepsilon}$ by $\dot{G}_t^\varepsilon$ in (3.31): in fact, because of translation invariance and because the field $X_t^{L,\varepsilon}$ is periodic, so is the solution $\Psi_t^{L,\varepsilon}$ and $F_t^\varepsilon(\Psi_t^{L,\varepsilon})$. We make this precise in the following lemma, which will be proven in Section 4.

**Lemma 3.22.** *Let* $(\Psi_t^{L,\varepsilon}, R_t^{L,\varepsilon})$ *solve the system*

$$\begin{aligned}\Psi_t^{L,\varepsilon} &= -\int_0^t \dot{G}_s^\varepsilon(F_s^\varepsilon(\Psi_s^{L,\varepsilon}) + R_s^{L,\varepsilon})\mathrm{d}s + X_t^{L,\varepsilon},\\ R_t^{L,\varepsilon} &= \int_t^\infty \omega_t(\mathcal{H}_s[F_s^\varepsilon](\Psi_s))\mathrm{d}s + \int_t^\infty \omega_t(\mathrm{D}F_s^\varepsilon(\Psi_s^{L,\varepsilon})\cdot\dot{G}_s^\varepsilon R_s^{L,\varepsilon})\mathrm{d}s\end{aligned}$$

*with* $\mathcal{H}_s[F_s^\varepsilon] := \partial_s F_s^\varepsilon + \tfrac{1}{2}\mathrm{D}_{\dot{G}_s^\varepsilon}^2 F_s^\varepsilon + \mathrm{D}F_s^\varepsilon\cdot\dot{G}_s^\varepsilon F_s^\varepsilon$. *Then,* $\Psi_t^{L,\varepsilon}$ *is periodic and solves the FBSDE (3.30).*

## 3.4 Sobolev spaces and covariance estimates

Before delving into the details of the solution, we conclude the section by introducing suitable norms on fields on $\mathbb{R}^d$[3.3]. The crucial property of (regularised) Grassmann field is that they are (smooth) bounded objects, unlike standard (commutative) Gaussian random variables. This motivates working with the space $L^\infty(\mathbb{R}^d;\mathcal{M}^4)$, which is introduced in the usual way as the $L^\infty$ space of functions taking values in a Banach space. This space is equipped with the norm

$$\|\varphi\|_{L^\infty(\mathbb{R}^d;\mathcal{M}^4)} := \sup_{\mu\in\{\uparrow,\downarrow\}\times\{\pm\}}\sup_{x\in\mathbb{R}^d}\|\varphi_{x,\mu}\|.$$

We will also need to control the size of the derivatives of Grassmann fields. If $n \in \mathbb{N}_0$ and $\varphi\colon \mathbb{R}^d \to \mathcal{M}^4$ is sufficiently smooth, we introduce the notation

$$\partial^0\varphi := \varphi \qquad \partial^n\varphi := \partial^n\varphi := (\partial_\nu^n\varphi)_{\nu\in\{1,\ldots,d\}^n},$$

to denote the tensor of the fields with $n$ derivatives, where if $\nu = (\nu_i)_{i=1}^n \in \{1,\ldots,d\}^n$

$$\partial_\nu^n\varphi(x) := \frac{\partial^n}{\partial^\nu x}\varphi(x) = \frac{\partial^n\varphi(x)}{\partial x_{\nu_1}\cdots\partial x_{\nu_n}}.$$

We introduce the following weighted Sobolev spaces

$$\mathscr{C}_t^n(\mathbb{R}^d;\mathcal{M}^4) := \{\varphi\colon \mathbb{R}^d \to \mathcal{M}^4 \mid \|\varphi\|_{\mathscr{C}_t^n(\mathbb{R}^d;\mathcal{M}^4)} < \infty\}$$

---

3.3. If fields are on the torus $\mathbb{T}_L^d$, the definitions below carry over by simply identifying $\mathbb{T}_L^d$ with the box $(-L/2,L/2]^d \subset \mathbb{R}^d$.

where
$$\|\varphi\|_{\mathscr{C}_t^n(\mathbb{R}^d;\mathscr{M}^4)} := \sum_{\nu \in \{1,\ldots,d\}^n} 2^{-nt} \|\partial_\nu^n \varphi\|_{L^\infty(\mathbb{R}^d;\mathscr{M}^4)}. \tag{3.32}$$

Note that the seminorms (3.32) controls only the size of the derivatives and is therefore a norm only after quotienting, which we left implicit in the definition of $\mathscr{C}_t^n(\mathbb{R}^d;\mathscr{M}^4)$. For simplicity, we abridge the notation to $\|\cdot\|_{L^\infty}$, $\|\cdot\|_{\mathscr{C}_t^n}$. For controlling the infinite-volume limit, we shall also consider the following weighted seminorms

$$\begin{aligned}
\|\varphi\|_{L^\infty(\eta)} &:= \sup_{\mu \in \{\uparrow,\downarrow\}\times\{\pm\}} \sup_{x \in \mathbb{R}^d} \varrho_\eta(x) \|\varphi_{x,\mu}\|, \\
\|\varphi\|_{\mathscr{C}_t^n(\eta)} &:= \sum_{\nu \in \{1,\ldots,d\}^n} 2^{-nt} \|\partial_\nu^n \varphi\|_{L^\infty(\eta)},
\end{aligned} \tag{3.33}$$

where $\varrho_\eta(x) := (1+|x|^2)^{-\eta/2}$ satisfies the compatibility condition $\varrho_\eta(x) \lesssim_\eta \frac{\varrho_\eta(y)}{\varrho_\eta(y-x)}$, as follows straightforwardly by the triangular inequality.

In Section 4.1, we study multilinear operators in the fields. To this end, we also consider fields on $\mathbb{R}^{dk}$ for some $k \in \mathbb{N}$, as maps $\varphi \colon \mathbb{R}^{dk} \to \mathscr{M}^{4^k}$, since in fact, if $\varphi$ is a field on $\mathbb{R}^2$, the tensor product $\psi^{\otimes k}$ is a field on $\mathbb{R}^{dk}$. Given $\varphi \colon \mathbb{R}^{2q} \to \mathscr{M}^{4^k}$ and a tuple $m = (m_i)_{i=1}^k$, we let $\partial^m \psi := (\partial^{m_1} \otimes \cdots \otimes \partial^{m_q}) \psi$. For the sake of brevity, we abridge our notation to

$$\psi^{(m)} := \partial^m \psi.$$

If $\psi$ is a field on $\mathbb{R}^2$ and $m = (m_i)_{i=1}^k$ is a tuple, we will further abridge our notation to

$$\psi^{\otimes(m)} := (\psi^{\otimes k})^{(m)} = (\partial^{m_1}\psi) \otimes \cdots \otimes (\partial^{m_k}\psi).$$

On fields on $\mathbb{R}^{dk}$, we consider the seminorms

$$\|\varphi\|_{\mathscr{C}_t^m(\mathbb{R}^{dk};\mathscr{M}^{4^k})} := \sum_{\substack{\nu_1,\ldots,\nu_k \\ |\nu_i|=m_i}} 2^{-[m]t} \|(\partial^{m_1} \otimes \cdots \otimes \partial^{m_k})\varphi\|_{L^\infty(\mathbb{R}^{dk};\mathscr{M}^{4^k})}. \tag{3.34}$$

where $[m] := \sum_i m_i$, and let

$$\mathscr{C}_t^m(\mathbb{R}^{dk};\mathscr{M}^{4^k}) := \Big\{\varphi \colon \mathbb{R}^{dk} \to \mathscr{M}^{4^k} \Big| \|\varphi\|_{\mathscr{C}_t^m(\mathbb{R}^{dk};\mathscr{M}^{4^k})} < \infty\Big\},$$

where the quotienting is again implicit. We conclude this section with some estimates on the covariance and on the GBM of Definition 3.18.

**Lemma 3.23.** *Let $\varepsilon \in 2^{-\mathbb{N}_\infty}$ and let $\dot{G}_t^\varepsilon$ be as in Definition 3.16. Then,*

$$\sup_{\mu,\mu'} \sup_{\nu \in \{1,\ldots,d\}^n} \left|\frac{\partial^n}{\partial^\nu x}(\dot{G}_s^\varepsilon(x;y))_{\mu,\mu'}\right| \lesssim_n 2^{(2\gamma+n)s} e^{-\bar{c}(2^s|x-y|)^\delta}, \tag{3.35}$$

*for some universal $\bar{c}$, with $\delta = 1$ if $\varepsilon = 0$ and otherwise as in Definition 3.16. If $\varepsilon > \varepsilon' \in 2^{-\mathbb{N}_\infty}$ and $\theta > 0$, then*

$$\sup_{\mu,\mu'} \sup_{\nu \in \{1,\ldots,d\}^n} \left|\frac{\partial^n}{\partial^\nu x}(\dot{G}_s^\varepsilon(x;y) - \dot{G}_s^{\varepsilon'}(x;y))_{\mu,\mu'}\right| \lesssim_{n,\theta} \varepsilon^\theta 2^{(2\gamma+n+\theta)s} e^{-\bar{c}(2^s|x-y|)^\delta}. \tag{3.36}$$

Because of the translational invariance of $\dot{G}_s^\varepsilon$, we can consider it to be a $\mathbb{C}^{4\times 4}$-valued function on $\mathbb{R}^d$. We let $L^p_{s;c}(\mathbb{R}^d;\mathbb{C}^{4\times 4})$ be the $L^p$ space restricted on $\mathbb{R}^d$ with respect to the measure $w_{s;c}(x)dx$, where

$$w_{s;c}(x) := e^{c(2^s|x|)^\delta} \tag{3.37}$$

for $c \geqslant 0$ and with $\delta$ as in Definition 3.16. The following corollary is a simple consequence of Lemma 3.23.

**Corollary 3.24.** *Under the same assumptions of Lemma 3.23, we have*

$$\sup_{\nu \in \{1,\ldots,d\}^n} \|\partial_\nu^n \dot{G}_s^\varepsilon\|_{L^p_{s;c}} \lesssim_n 2^{(2\gamma+n-d/p)s}, \tag{3.38}$$

$$\sup_{\nu \in \{1,\ldots,d\}^n} \|\partial_\nu^n (\dot{G}_s^\varepsilon - \dot{G}_s^{\varepsilon'})\|_{L^p_{s;c}} \lesssim_{n,\theta} 2^{(2\gamma+n+\theta-d/p)s} \varepsilon^\theta,$$

*for $c < \bar{c}$, where $\lesssim$ depends on the weight as well.*

**Proof of Lemma 3.23.** We focus on the case $s \in [1, -\log_2 \varepsilon]$, since for $s \in [0,1)$ the proof is similar. We note that $(\partial_\nu^n \dot{G}_s^\varepsilon)_{\mu,\mu'}$ is proportional to the integral (up to universal factors)

$$\left(\frac{\partial^n}{\partial^\nu x} \dot{G}_s^\varepsilon\right)_{\mu,\mu'}(y;x) \sim 2^{(2\gamma+n)s} \int_{\mathbb{R}^d} e^{ik2^s(y-x)} e^{-(k^2+2^{-2s})} \left[\prod_{i=1}^n -ik_{\nu_i}\right] \varphi_\delta(2^s \varepsilon |k|) dk. \tag{3.39}$$

We want to prove a bound of the form

$$|2^{ms}|y-x|^m (\partial_\nu^n \dot{G}_s^\varepsilon)_{\mu,\mu'}(y;x)| \lesssim_n C^m (m!)^{\delta-1} \quad \forall m \in \mathbb{N} \tag{3.40}$$

from which (3.35) for $\varepsilon > 0$ follows by optimising over $m \in \mathbb{N}$. To prove (3.40), we use integration by parts on (3.39)

$$|2^{ms}|y-x|^m (\partial_\nu^n \dot{G}_s^\varepsilon)_{\mu,\mu'}(y;x)|$$
$$\lesssim 2^{(2\gamma+n)s} \tilde{C}^n \sup_{\nu'} \int_{\mathbb{R}^d} \left|\frac{\partial^m}{\partial^{\nu'} k} e^{-(k^2+2^{-2s})} \left[\prod_{i=1}^n -ik_{\nu_i}\right] \varphi_\delta(2^s \varepsilon |k|)\right| dk.$$

By using the Gevrey condition, we have

$$\sup_{k \in \mathbb{R}^d} \left|\frac{\partial^m}{\partial^{\nu'} k} \varphi_\delta(2^s \varepsilon |k|)\right| \leqslant C^{1+m} (2^s \varepsilon)^m (m!)^{\delta-1}$$

whereas, by using the analyticity in a strip of size one around the real axis we have $\left|\frac{\partial^n}{\partial^{\nu'} k} e^{-(k^2+2^{-2s})}\right| \leqslant C^{1+m} m! e^{-k^2}$. Since $2^s \varepsilon \leqslant 1$, (3.40) follows.

Repeating the same strategy for the difference, we obtain

$$\left|2^{ms}|y-x|^m \left[\left(\frac{\partial^n}{\partial^\nu x} \dot{G}_s^\varepsilon\right)_{\mu,\mu'}(y;x) - \left(\frac{\partial^n}{\partial^\nu x} \dot{G}_s^{\varepsilon'}\right)_{\mu,\mu'}(y;x)\right]\right|$$
$$\lesssim_n 2^{(2\gamma+n)s} \tilde{C}^n \sup_{\nu'} \int_{\mathbb{R}^d} \left|\frac{\partial^m}{\partial^{\nu'} k} e^{-(k^2+2^{-2s})} \left[\prod_{i=1}^n -ik_{\nu_i}\right] (\varphi_\delta(2^s \varepsilon |k|) - \varphi_\delta(2^s \varepsilon' |k|))\right| dk. \tag{3.41}$$

We note that $\varphi_\delta(2^s \varepsilon |k|) - \varphi_\delta(2^s \varepsilon' |k|)$ and its derivatives are vanishing unless $|k| \gtrsim \varepsilon^{-1} 2^{-s}$, therefore

$$\left|2^{ms}|y-x|^m \left[\left(\frac{\partial^n}{\partial^\nu x} \dot{G}_s^\varepsilon\right)_{\mu,\mu'}(y;x) - \left(\frac{\partial^n}{\partial^\nu x} \dot{G}_s^{\varepsilon'}\right)_{\mu,\mu'}(y;x)\right]\right|$$
$$\lesssim_n 2^{(2\gamma+n)s} C^m (m!)^{\delta-1} e^{-\frac{(\varepsilon 2^s)^{-2}}{2}}$$
$$\leqslant 2^{(2\gamma+n+\theta)s} C^m (m!)^{\delta-1} \sup_{s \in \mathbb{R}^+} 2^{-\theta s} e^{-\frac{(\varepsilon 2^s)^{-2}}{2}}$$
$$\lesssim_\theta 2^{(2\gamma+n+\theta)s} \varepsilon^\theta C^m (m!)^{\delta-1},$$

so that, optimising over $m \in \mathbb{N}$ gives the stretched exponential decay in (3.36).

When $\varepsilon = 0$, we have $\varphi_\delta(2^s \varepsilon |k|) \equiv 1$ and thus the integrand is analytic in a strip around the real axis in each of the variables $(k_i)_{i=1}^d$. Accordingly, we deform the integration contour by shifting $k_1 \to k_1 + i\operatorname{sign}(y_1-x_1)/2$ and obtain

$$\left| \left( \frac{\partial^n}{\partial^\nu x} \dot{G}_s \right)_{\mu,\mu'} (2^{-s}y; 2^{-s}x) \right| \lesssim 2^{(2\gamma+n)s} e^{-|y_1-x_1|/2} \int_{\mathbb{R}^d} e^{-\operatorname{Re}\left(k^2 - \frac{1}{2} + ik_1 \operatorname{sign}(y_1-x_1)\right)} |k|^n \mathrm{d}k$$
$$\lesssim_n 2^{(2\gamma+n)s} e^{-|y_1-x_1|/2}.$$

By rotation invariance this implies the bound (3.35). □

**Remark 3.25.** In a similar way and by using the Poisson summation formula to restrict to finite $L$, one can prove a similar decay estimate for $\mathfrak{C}_s^{L,\varepsilon}$, that is,

$$\sup_{\mu,\mu'} \sup_{\nu \in \{1,\ldots,d\}^n} \left| \frac{\partial^n}{\partial^\nu x} (\mathfrak{C}_s^{L,\varepsilon}(x;y))_{\mu,\mu'} \right| \lesssim 2^{(\gamma+d/2+n)s} \sum_{m \in \mathbb{Z}^d} e^{-\bar{c}(2^s|x-y-mL|)^\delta},$$

with $\delta = 1$ if $\varepsilon = 0$, and for the difference $\frac{\partial^n}{\partial^\nu x}\big(\mathfrak{C}_s^{L,\varepsilon}(x;y) - \mathfrak{C}_s^{L,\varepsilon'}(x;y)\big)_{\mu,\mu'}$. Note also the bounds

$$\sup_{\nu \in \{1,\ldots,d\}^n} \|\partial_\nu^n \mathfrak{C}_s^{L,\varepsilon}\|_{L^p_{s;c}(\mathbb{T}_L^d; \mathbb{C}^{4 \times 4})} \lesssim_{n,c} 2^{(\gamma+d/2+n-d/p)s}. \tag{3.42}$$

We also have the following estimates.

**Corollary 3.26.** *Let $X_t^{L,\varepsilon}$ be the field in Definition 3.20. Then, for any $L' \geqslant L \in \mathbb{N}_\infty$, $\varepsilon \geqslant \varepsilon' \in 2^{-\mathbb{N}_\infty}$, $\theta \geqslant 0$ and $\eta \geqslant 0$*

$$\|X_t^{L,\varepsilon}\|_{\mathscr{C}_t^n} \lesssim_n 2^{\gamma t}, \qquad \|X_t^{L,\varepsilon} - X_t^{L',\varepsilon'}\|_{\mathscr{C}_t^n(\eta)} \lesssim_{n,\theta,\eta} 2^{\gamma t}[2^{\theta t} \varepsilon^\theta \mathbf{1}_{\varepsilon > \varepsilon'} + L^{-\eta} \mathbf{1}_{L' > L}].$$

*Furthermore, for any $L' \geqslant L \in \mathbb{N}_\infty$, $\varepsilon \geqslant \varepsilon' \in 2^{-\mathbb{N}_\infty}$, $\theta \geqslant 0$, $\eta \geqslant 0$ and $\rho > 0$, uniformly in $t \in [0,\infty]$,*

$$\|X_t^{L,\varepsilon}\|_{B_{\infty,\infty}^{-\gamma-\rho}} \lesssim \rho^{-1/2}, \qquad \|X_t^{L,\varepsilon} - X_t^{L',\varepsilon'}\|_{B_{\infty,\infty,\eta}^{-\gamma-\theta-\rho}} \lesssim \rho^{-1/2}(\varepsilon^\theta \mathbf{1}_{\varepsilon > \varepsilon'} + L^{-\eta} \mathbf{1}_{L' > L}). \tag{3.43}$$

**Proof.** The proof is a consequence of the norm compatibility and the exponential estimates for $\mathfrak{C}_t^{L,\varepsilon}$, $\mathfrak{C}_t^{L,\varepsilon} - \mathfrak{C}_t^{L,\varepsilon'}$, see Remarks 3.25, and $\mathfrak{C}_t^\varepsilon - \mathfrak{C}_t^{L,\varepsilon}$. When we consider differences of fields, we shall use

$$\|X_t^{L,\varepsilon}(f) - X_t^{L',\varepsilon'}(f')\| \leqslant \|X_t^{L,\varepsilon}(f) - X_t^{L',\varepsilon'}(f)\| + \|X_t^{L',\varepsilon'}(f-f')\|,$$

with, e.g., $f = \delta_x^{L,\varepsilon} \oplus 0 \in \mathfrak{h}_L \cong L^2(\mathbb{T}_L^d; \mathbb{C}^4)$ and $f' = \delta_x^{L',\varepsilon'} \oplus 0 \in \mathfrak{h}_{L'}$ or their derivatives, see Definition 3.20. Accordingly, we have that

$$\|\partial_\nu^n X_{t,x}^{L,\varepsilon} - \partial_\nu^n X_{t,x}^{L,\varepsilon'}\| \lesssim \left( \int_0^t \|(\mathfrak{C}_s^{L,\varepsilon} - \mathfrak{C}_s^{L,\varepsilon'}) \partial_\nu^n \delta_x^{L,\varepsilon}\|_{\mathfrak{h}_L}^2 \mathrm{d}s \right)^{1/2}$$
$$+ \left( \int_0^t \|\mathfrak{C}_s^{L,\varepsilon'} \partial_\nu^n (\delta_x^{L,\varepsilon} - \delta_x^{L,\varepsilon'})\|_{\mathfrak{h}_L}^2 \mathrm{d}s \right)^{1/2}$$
$$\lesssim_{n,\theta} 2^{(\gamma+n+\theta)t} \varepsilon^\theta,$$

where in the first term we used that $|\varphi_\delta(\varepsilon|k|) - \varphi_\delta(\varepsilon'^{-1}|k|)| \lesssim (|k|\varepsilon)^\theta |\varphi_\delta(\varepsilon|k|) - \varphi_\delta(\varepsilon'^{-1}|k|)|$, compare with the definition of $\mathfrak{C}_s^{L,\varepsilon}$, whereas in the second that the Fourier series of $\delta_x^{L,\varepsilon} - \delta_x^{L,\varepsilon'}$ starts from $\gtrsim \varepsilon^{-1}$ and that $\mathbf{1}_{|k| \geqslant \varepsilon^{-1}} \lesssim \mathbf{1}_{|k| \geqslant \varepsilon^{-1}}(|k|\varepsilon)^\theta$. This proves the bound on $\|X_t^{L,\varepsilon} - X_t^{L,\varepsilon'}\|_{\mathscr{C}_t^n}$. The bound on $\|X_t^{L,\varepsilon}\|_{\mathscr{C}_t^n}$ follows in an easier way.

The claim on $\|X_t^{L,'\varepsilon} - X_t^{L,\varepsilon}\|_{\mathscr{C}_t^n(\eta)}$, follows by proving the same estimate for $\|X_t^\varepsilon - X_t^{L,\varepsilon}\|_{\mathscr{C}_t^n(\eta)}$. We have

$$\|\partial_\nu^n X_{t,x}^\varepsilon - \partial_\nu^n X_{t,x}^{L,\varepsilon}\| \lesssim \left(\int_0^t \|(\mathfrak{C}_s^\varepsilon - \mathfrak{C}_s^{L,\varepsilon})\partial_\nu^n \delta_x^\varepsilon\|_{\mathfrak{h}}^2 ds\right)^{1/2}$$
$$+ \left(\int_0^t \|\mathfrak{C}_s^\varepsilon \partial_\nu^n (\delta_x^\varepsilon - \delta_x^{L,\varepsilon})\|_{\mathfrak{h}}^2 ds\right)^{1/2}$$
$$=: a_x + b_x.$$

where $\mathfrak{h} \cong L^2(\mathbb{R}^d; \mathbb{C}^4)$. To control $\sup_{x \in \mathbb{R}^d} \varrho_\eta(x) a_x$ we split the integration coming from $\|\cdot\|_{\mathfrak{h}}^2$ as follows:

$$\sup_{x \in \mathbb{R}^d} \varrho_\eta(x) a_x \lesssim \sup_{x \in B_{L/2}} \left(\int_0^t \int_{B_L} \frac{\varrho_\eta(z)^2}{\varrho_\eta(z-x)^2} \left|\frac{\partial^n}{\partial^\nu x}(\mathfrak{C}_s^\varepsilon - \mathfrak{C}_s^{L,\varepsilon})(x;z)\right|^2 dz ds\right)^{1/2}$$
$$+ \sup_{x \in (B_{L/2})^c} \varrho_\eta(x) \left(\int_0^t \int_{B_L} \left|\frac{\partial^n}{\partial^\nu x}(\mathfrak{C}_s^\varepsilon - \mathfrak{C}_s^{L,\varepsilon})(x;z)\right|^2 dz ds\right)^{1/2} \quad (3.44)$$
$$+ \sup_{x \in \mathbb{R}^d} \left(\int_0^t \int_{(B_{L/2})^c} \frac{\varrho_\eta(z)^2}{\varrho_\eta(z-x)^2} \left|\frac{\partial^n}{\partial^\nu x} \mathfrak{C}_s^\varepsilon(x;z)\right|^2 dz ds\right)^{1/2}$$
$$=: \mathrm{I} + \mathrm{II} + \mathrm{III},$$

where we used the property $\varrho_\eta(x) \lesssim_\eta \varrho_\eta(z)/\varrho_\eta(z-x)$ and where $B_L := [-L/2, L/2]^d$. We bound I by noting that for $x, z \in B_L$ $(\mathfrak{C}_s^\varepsilon - \mathfrak{C}_s^{L,\varepsilon})(x;z) = -\sum_{m \in \mathbb{Z}^d \setminus \{0\}} \mathfrak{C}_s^\varepsilon(x; z + mL)$ and by using the estimate for $\partial^\nu \mathfrak{C}_s^\varepsilon$ for $x \in B_{L/2}$, which controls the growing term $\varrho_\eta^{-1}$. The fact that $|x - z| \le L/2$ and the term $m = 0$ is removed from the Poisson summation allows us to gain extra decay in $L$. Hence, we have

$$\mathrm{I} \lesssim \left(\int_0^t 2^{(2\gamma+2n)s} 2^{-\bar{c}(2^s L)^\delta/2} \sum_{m,m' \in \mathbb{Z}^d} e^{-\frac{|m|+|m'|}{2}L} ds\right)^{1/2} \lesssim 2^{-\bar{c}L^\delta/2}.$$

We bound II by using the decay of $\mathfrak{C}_s^\varepsilon$ and $\mathfrak{C}_s^{L,\varepsilon}$ separately to integrate over $B_L$, the decay in $L$ being obtained directly from the weight $\varrho_\eta$, since $x$ is away from the origin; thus, we have $\mathrm{II} \lesssim 2^{(\gamma+n)t} L^{-\eta}$. To bound III we use $\varrho_\eta(z) \lesssim L^{-\eta}$ and control the remaining integral using the decay of $\mathfrak{C}_s^\varepsilon$, implying the same bound as for I and II.

To control $\sup_{x \in \mathbb{R}^d} \varrho_\eta(x) b_x$ we split the integral as in (3.44) and follow the same strategy as above: in particular, if $z \in \mathbb{T}_L^d$ we use the Poisson summation formula $\delta_x^\varepsilon(z) - \delta_x^{L,\varepsilon}(z) = -\sum_{m \in \mathbb{Z}^d \setminus \{0\}} \delta_x^\varepsilon(z + mL)$ to write $(\mathfrak{C}_s^\varepsilon \partial_\nu^n (\delta_x^\varepsilon - \delta_x^{L,\varepsilon}))(z) = -\sum_{m \in \mathbb{Z}^d \setminus \{0\}} \frac{\partial^n}{\partial^\nu x} \mathfrak{C}_s^\varepsilon(x + mL; z)$, whereas for $z \in (B_L)^c$ we simply have $(\mathfrak{C}_s^\varepsilon \partial_\nu^n (\delta_x^\varepsilon - \delta_x^{L,\varepsilon}))(z) = \frac{\partial^n}{\partial^\nu x} \mathfrak{C}_s^\varepsilon(x; z)$ and the extra decay in $L$ is obtained from $\varrho_\eta$. Therefore $\sup_{x \in \mathbb{R}^d} \varrho_\eta(x) b_x \lesssim 2^{(\gamma+n)t} L^{-\eta}$ and thus the desired bound is proven.

We are left with proving (3.43). Denote by $(\Delta_j)_{j \ge -1}$ the Littlewood–Paley blocks. We have

$$\|X_t^{L,\varepsilon}(\Delta_j(\cdot - x)) - X_t^{L,\varepsilon'}(\Delta_j(\cdot - x))\| \lesssim \left(\int_0^t \|(\mathfrak{C}_s^{L,\varepsilon} - \mathfrak{C}_s^{L,\varepsilon'}) \mathbf{1}_{B_L} \Delta_j(\cdot - x)\|_{\mathfrak{h}_L}^2 ds\right)^{1/2}$$
$$\lesssim \left(\int_0^t 2^{-2\rho s} \varepsilon^{2\theta} \|(\mathbb{1} - \Delta)^{\gamma/2 + \theta/2 + \rho/2 - d/4} \Delta_j(\cdot - x)\|_{\mathfrak{h}}^2 ds\right)^{1/2}$$
$$\lesssim \rho^{-1/2} \varepsilon^\theta 2^{(\gamma+\theta+\rho)j},$$

$\Delta$ being the Laplacian, where we used that $|\varphi_\delta(\varepsilon|k|) - \varphi_\delta(\varepsilon^{-1}|k|)|^2 \lesssim (|k|\varepsilon)^{2\theta} |\varphi_\delta(\varepsilon|k|) - \varphi_\delta(\varepsilon^{-1}|k|)|^2$ and extracted sufficient decay in $s$ by taking $\rho > 0$. The bound on $\|X_t^{L,\varepsilon}(\Delta_j(\cdot - x))\|$ follows likewise. Finally, we prove the last bound following a strategy similar to (3.44). We split $\|\cdot\|_{\mathfrak{h}}^2$ as

$$\|(\mathfrak{C}_s^\varepsilon - \mathfrak{C}_s^{L,\varepsilon}) w^{-1} w \Delta_j(\cdot - x)\|_{\mathfrak{h}}^2$$
$$\lesssim \|\mathbf{1}_{B_L}(\mathfrak{C}_s^\varepsilon - \mathfrak{C}_s^{L,\varepsilon}) w^{-1} \mathbf{1}_{B_{L/2}} w \Delta_j(\cdot - x)\|_{\mathfrak{h}}^2 + \|\mathbf{1}_{(B_L)^c} \mathfrak{C}_s^\varepsilon w^{-1} \mathbf{1}_{B_{L/2}} w \Delta_j(\cdot - x)\|_{\mathfrak{h}}^2$$
$$+ \|(\mathfrak{C}_s^\varepsilon - \mathfrak{C}_s^{L,\varepsilon}) w^{-1} \mathbf{1}_{(B_{L/2})^c} w \Delta_j(\cdot - x)\|_{\mathfrak{h}}^2 =: \mathrm{I} + \mathrm{II} + \mathrm{III},$$

where $w = (\mathbb{1} - \Delta)^{\gamma/2+\rho/2-d/4}$. We bound I by using Young's inequality and the Poisson summation formula $(\mathfrak{C}_s^\varepsilon - \mathfrak{C}_s^{L,\varepsilon})(x;z) = -\sum_{m \in \mathbb{Z}^d \setminus \{0\}} \mathfrak{C}_s^\varepsilon(x; z + mL)$, so that $\text{I} \lesssim 2^{-2\rho s} 2^{-\bar{c}(2^s L)^\delta/2} \|w \Delta_j(\cdot - x)\|_{\mathfrak{h}}^2$. We bound I I by using the Cauchy–Schwarz inequality and by noting that for $x \in (B_L)^c$

$$\|\mathfrak{C}_s^\varepsilon(x;\cdot) w^{-1} \mathbf{1}_{B_{L/2}}\|_{\mathfrak{h}_L}^2 \lesssim 2^{(d-2\rho)s} e^{-\bar{c}(2^s|x|)^\delta/2},$$

which can be integrated in $dx$ and thus $\text{I I} \lesssim 2^{-2\rho s} 2^{-\bar{c}(2^s L)^\delta/2} \|w \Delta_j(\cdot - x)\|_{\mathfrak{h}}^2$. Finally, we estimate the last term by using separately the decay of $\mathfrak{C}_s^\varepsilon$ and $\mathfrak{C}_s^{L,\varepsilon}$ so that $\|(\mathfrak{C}_s^\varepsilon - \mathfrak{C}_s^{L,\varepsilon}) w^{-1}\|_{L^1(\mathbb{R}^d;\mathbb{C}^4)}^2 \lesssim 2^{-2\rho s}$ and $\text{I I I} \lesssim 2^{-2\rho s} \|\mathbf{1}_{(B_{L/2})^c} w \Delta_j(\cdot - x)\|_{\mathfrak{h}}^2$. The terms I and I I decay already in $L$. Concerning the term I I I, denote by $K_j$ the kernel of $w \Delta_j$ and note the faster than any power decay $|K_j(x;y)| \lesssim \frac{2^{(\gamma+\rho+d)j/2} C_N}{1+2^j|x-y|^N}$, $N \in \mathbb{N}$. Therefore

$$\sup_{x \in \mathbb{R}^d} \varrho_\eta(x) \left( \int_0^t \text{I I I} \, ds \right)^{1/2} \lesssim \rho^{-1/2} \sup_{x \in \mathbb{R}^d} \varrho_\eta(x) \|\mathbf{1}_{(B_{L/2})^c} w \Delta_j(\cdot - x)\|_{\mathfrak{h}}$$
$$\lesssim \rho^{-1/2} L^{-\eta} 2^{(\gamma+\rho)j}.$$

This concludes the proof. □

## 4 Solution in the subcritical regime

In Section 4.1 we introduce Banach spaces of Grassmann polynomials as suitable multilinear operators in the fields, possibly with derivatives. In Section 4.2 we truncate the Polchinski flow equation and study global polynomial solutions in such Banach spaces. Note that the closer to the critical regime (that is, the closer $\gamma$ is to $d/4$), the less we can truncate the flow equation. Furthermore, the less we truncate the equation, the more regular is the remainder produced by the truncation. In Section 4.3 we prove global existence of the pair $(\Psi, R)$ by a fixed-point argument, compare with Lemma 3.22, prove Theorem 1.1 and, as a byproduct, prove the existence of the weak limit $\omega^V$ stated in Theorem 1.3. Finally, in Section 4.4 we show the exponential decay of clustered correlations and the short-distance divergence of the two-point function and thus prove Theorem 1.3.

### 4.1 Spaces of Grassmann monomials

We solve the truncated flow equation in suitable spaces of multilinear operators on fields and their derivatives. For technical reasons, it will suffice to consider fields with either no or two derivatives. Accordingly, to keep track of the derivatives, we henceforth consider only tuples of the form $m = (m_i)_{i=1}^k \in \{0,2\}^k$ for some $k \in \mathbb{N}$. If $m = (m_i)_{i=1}^k$ for some $k \in \mathbb{N}$, we denote its length by $|m| := k$ and use $[m] := \sum_{i=1}^k m_i$ to count the total number of derivatives. We will also consider the empty tuple, denoted by $\emptyset$. We naturally set $|\emptyset| = [\emptyset] = 0$. Letting $m = (m_i)_{i=1}^k$, $m = (m_i')_{i=1}^{k'}$ and $m'' = (m_i'')_{i=1}^{k''}$ be possibly empty tuples:

- $m' \subseteq m$, if $m'$ can be formed by possibly removing elements from $m$; note that $m \subseteq m$ and $\emptyset \subseteq m$ for any possibly empty tuple $m$.
- $m \circ m'$ denotes the concatenated tuple obtained by juxtaposition, that is, $m \circ m' = ((m \circ m')_i)_{i=1}^{k+k'}$ with $(m \circ m')_i = m_i$ for $i = 1, \ldots, k$ and $(m \circ m')_i = m_i'$ for $i = k+1, \ldots, k+k'$.
- $S_m$ denotes the set of permutations $\pi$ of $(1, \ldots, k)$ such that $m_{\pi(i)} = m_i$.
- $S_{m',m''}^m$ denotes the set of permutations of $(1, \ldots, k)$ into $(I', I'')$ such that $I' = (i_1' \ldots i_{k'}')$ and $I'' = (i_1'' \ldots i_{k''}'')$ are ordered (that is, $i_1' < \ldots < i_{k'}'$ and $i_1'' < \ldots < i_{k''}''$) and such that $m' = (m_{i_j'})_{j=1}^{k'}$ and $m'' = (m_{i_j''})_{j=1}^{k''}$. In particular, if either $m'$ or $m''$ are empty, $S_{m',m''}^m$ consists of the trivial permutation only.

**Remark 4.1.** Note the following binomial expansion, for $\psi, \varphi$ odd variables,

$$(\psi + \varphi)^{\otimes(m)} = \sum_{m', m'' \subseteq m} \sum_{\pi \in S_{m',m''}^m} s(\pi) P_\pi(\psi^{\otimes(m')} \otimes \varphi^{\otimes(m'')}),$$

where, $s$ is a suitable sign depending on $\pi$, $P_\pi$ lifts the action of $\pi$ to $(P_\pi \psi)((x_1, \mu_1); \ldots; (x_k, \mu_k)) = \psi((x_{\pi(1)}, \mu_{\pi(1)}); \ldots; (x_{\pi(k)}, \mu_{\pi(k)}))$. Note of course that the right-hand side is vanishing if $|S_{m',m''}^m| = 0$.

Let $\mathcal{B}(V; W)$ denote the Banach space of bounded linear operators from V to W. We are interested in the following spaces of linear operators

$$\mathcal{L}_t^{(k,m)} := \mathcal{B}(\mathcal{C}_t^m(\mathbb{R}^{dk}; \mathcal{M}^{4^k}); L^\infty(\mathbb{R}^d; \mathcal{M}^4)),$$

where $k \in \mathbb{N}$ and $m$ is a tuple of length $k$. With any operator $K \in \mathcal{L}_t^{(k,m)}$, we can associate an integral kernel, which is generally a distribution. For example, for any tuple $m = (m_i)_{i=1}^k$ we can write

$$(K_m(\psi^{(m)}))(y; \mu)$$
$$= \sum_{\substack{\mu_1, \ldots, \mu_k \\ \mu_i \in \{\uparrow, \downarrow\} \times \{\pm\}}} \sum_{\substack{\nu_1, \ldots, \nu_k: \\ \nu_i \in \{1, \ldots, d\}^{m_i}}} \int_{\mathbb{R}^{dk}} K_m((y, \mu); (x_1, \mu_1, \nu_1); \cdots; (x_k, \mu_k, \nu_k)) \times$$
$$\times \partial_{\nu_1}^{m_1} \otimes \cdots \partial_{\nu_k}^{m_k} \otimes \psi_{x_1, \mu_1; \ldots; x_k, \mu_k} \mathrm{d}x_1 \cdots \mathrm{d}x_k,$$

In our setting, such kernels are a.e. antisymmetric and invariant under translations and reflections across coordinate planes, that is, for any permutation $\pi \in S_m$ and any $z \in \mathbb{R}^d$

$$K_m((y, \mu); (x_1, \mu_1, \nu_1); \cdots; (x_k, \mu_k, \nu_k))$$
$$= K_m((y-z, \mu); (x_1-z, \mu_1, \nu_1); \ldots; (x_k-z, \mu_k, \nu_k))$$
$$= K_m((-y, \mu); (-x_1, \mu_1, \nu_1); \cdots; (-x_k, \mu_k, \nu_k))$$
$$= \mathrm{sign}(\pi) K_m((y, \mu); (x_{\pi(1)}, \mu_{\pi(1)}, \nu_{\pi(1)}); \cdots; (x_{\pi(k)}, \mu_{\pi(k)}, \nu_{\pi(k)})).$$

The multilinear operators we consider exhibit scale-dependent spatial decay in the separation of points. To keep track of this decay, we introduce the following scale-dependent weight, which is a generalisation of (3.37),

$$w_{t;c}^{(k)}(x_1; \ldots; x_k) := e^{c(2^t \mathrm{St}(x_1; \ldots; x_k))^\delta},$$

for $c \geq 0$, $\delta$ as in Definition 3.16, where $\mathrm{St}(x_1; \ldots; x_k)$ denotes the Steiner diameter of the set $\{x_1, \ldots, x_k\}$, see, e.g., [GMR21]. If $K \in \mathcal{L}_t^{(k,m)}$, we let $K \cdot w_{t;c}^{(k+1)}$ be the operator whose kernel is a.e. the point-wise product of the kernel of $K$ with $w_{t;c}^{(k+1)}(x_1; \ldots; x_{k+1})$. Equivalently,

$$((K \cdot w_{t;c}^{(k+1)})(\psi^{(m)}))(y; \mu) := (K(w_{t;c}^{(k+1)}(y; \cdot) \psi^{(m)}))(y; \mu).$$

We introduce

$$\mathcal{L}_t^{(k,m)}(s; c) := \{K \in \mathcal{L}_t^{(k,m)} | \|K \cdot w_{s;c}^{(k+1)}\|_{\mathcal{L}_t^{(k,m)}} < \infty\},$$

and for $K \in \mathcal{L}_t^{(k,m)}(s; c)$ let $\|K\|_{\mathcal{L}_t^{(k,m)}(s;c)} := \|K \cdot w_{s;c}^{(k+1)}\|_{\mathcal{L}_t^{(k,m)}}$. We finally define the spaces of monomials with spatial decay in the separation of points.

**Definition 4.2.** *Let $k \in \mathbb{N}$ and $m = (m_i)_i \in \{0, 2\}^k$. The weighted space of monomials of degree $k \in \mathbb{N}$ at scale $t \geq 0$ is defined as the following subspace of the direct sum of $(\mathcal{L}_t^{(k,m)})_m$:*

$$\mathcal{L}_t^{(k)}(s; c) := \left\{ (K_m)_m \in \bigoplus_{m \in \{0,2\}^k} \mathcal{L}_t^{(k,m)}(s; c) \,\middle|\, \text{antisym., trans. \& refl. invariant} \right\}.$$

*The space $\mathscr{L}_t^{(k)}(s;c)$ is equipped with the norm*

$$\|K\|_{\mathscr{L}_t^{(k)}(s;c)} = \sup_{m \in \{0,2\}^k} 2^{[m](s-t)} \|K_m\|_{\mathscr{L}_t^{(k,m)}(s;c)}. \tag{4.1}$$

*Finally, we let $\mathscr{L}_t^{(k)} := \mathscr{L}_t^{(k)}(0;0)$ denote the space of monomials of degree k, at scale t.*

**Remark 4.3.** Note that for any $s \geqslant s'$ and $K \in \mathscr{L}_t^{(k)}(s;c)$

$$\|K\|_{\mathscr{L}_t^{(k)}(s';c)} \leqslant \|K\|_{\mathscr{L}_t^{(k)}(s;c)}$$

and thus $\mathscr{L}_t^{(k)}(s;c) \subset \mathscr{L}_t^{(k)}(s';c)$, which follows because $w_{s';c}^{(k)} \leqslant w_{s;c}^{(k)}$.

**Remark 4.4.** Loosely speaking, in the norm $\|\cdot\|_{\mathscr{L}_t^{(k)}(s;c)}$, $t$ is the scale of the fields, whereas $s$ is the scale of the monomial itself. For controlling the flow equation, it is crucial to take $t \leqslant s$ (except for the local monomial of degree one, see the next section). Since the field scale is not larger than the monomial one, the renormalisation operator gives a gain, as can be seen in eq. (4.5) below and motivates the choice of the scale-dependent weight in $\|\cdot\|_{\mathscr{L}_t^{(k)}(s;c)}$.

As long as cut-offs are present, functions in the Grassmann fields are polynomials and can be therefore written as $F_t(\psi_t) := \sum_k F_t^{(k)}(\psi_t)$, where the sum contains only finitely many terms, $\psi_t$ being a Grassmann field and $F^{(k)} \in \mathscr{L}_t^{(k)}$, so that

$$F_t^{(k)}(\psi_t) = \sum_{m \in \{0,2\}^k} F_t^{(k,m)}(\psi_t^{\otimes (m)}). \tag{4.2}$$

This allows us to lift the flow equation to the space of sequences $(F^{(k)})_{k \geqslant 1}$, see Section 4.2 for details. Note, however, that the expansion (4.2) is not unique, see also Remark 4.6 below, and we will exploit this flexibility to write the lifted flow equation in such a way that the correct boundary data can be obtained, that is, the correct local "renormalised counterterm" in the interacting measure.

We now introduce two operators on $\mathscr{L}_t^{(1)}$ that allow us to rearrange the expansion in (4.2), namely the so-called localisation and renormalisation operators. These operators are well-known in the study of the renormalisation group flow, see, e.g., [Mas08, BBS19] and are crucial for our purposes because, as explained before, allow us to obtain non-trivial local boundary terms in the flow equation.

**Definition 4.5.** *(Localisation & Renormalisation). The localisation and renormalisation operators are respectively the projections* Loc *and* Ren *onto $\mathscr{L}_t^{(1)}$ defined as follows: if $K = (K_0, K_2)$, we let*

$$\text{Loc}\,K := (\text{Loc}_0 K_0, 0), \qquad \text{Ren}\,K := (0, K_2 + (\mathbb{1} - \text{Loc}_0) K_0),$$

*where*

$$(\text{Loc}_0 K_0)((y;\mu);(x;\mu')) := \delta(y-x) \int_{\mathbb{R}^d} K((y,\mu);(z;\mu'))\mathrm{d}z$$

*with $\delta(y-x)$ denoting the Dirac delta distribution.*

**Remark 4.6.** Consider $K(\psi) := \sum_{m \in \{0,2\}^k} K_m(\psi^{(m)})$. Note that even though $\text{Loc} + \text{Ren} \neq \mathbb{1}$ on $\mathscr{L}_t^{(1)}$, we still have the identity $K(\psi) = \text{Loc}\,K(\psi) + \text{Ren}\,K(\psi)$. In other words, we have the equivalent relation $K \sim (\text{Loc} + \text{Ren}) K = (\text{Loc}_0 K_0, K_2 + (\mathbb{1} - \text{Loc}_0) K_0)$, which we will use to properly rewrite the lifted flow equation on $\mathscr{L}_t^{(k)}$.

In view of the application to the flow equation, we provide the following lemma.

**Lemma 4.7.** *For any $r,s,t \geqslant 0$ and any $K \in \mathscr{L}_t^{(1)}$ we have*

$$\|\operatorname{Loc} K\|_{\mathscr{L}_r^{(1)}(s;c)} = \|\operatorname{Loc} K\|_{\mathscr{L}_t^{(1)}(s;c)} \leqslant \|K\|_{\mathscr{L}_t^{(1)}}. \tag{4.3}$$

*If $K \in \mathscr{L}_t^{(1)}(s;c)$, then for any $c' < c$*

$$\|\operatorname{Ren} K\|_{\mathscr{L}_t^{(1)}(s;c')} \lesssim_{c-c'} \|K\|_{\mathscr{L}_t^{(1)}(s;c)}. \tag{4.4}$$

**Proof.** The identity and the inequality in (4.3) follow by noting that $(\operatorname{Loc} K)_2 = 0$ and that $w_{s;\delta;c}^{(2)}(y;y) = 1$. To prove (4.4), one should note that, for $\psi$ smooth enough,

$$\begin{aligned}
&((\mathbb{1} - \operatorname{Loc}_0) K_0(\psi))(y;\mu) \\
&= \sum_{\mu'} \int_{\mathbb{R}^d} K_0((y,\mu);(x;\mu'))[\psi(x;\mu') - \psi(y;\mu') - (x-y) \cdot \partial \psi(y;\mu')]\mathrm{d}x \\
&= \sum_{\mu'} \int_{\mathbb{R}^d} K_0((y,\mu);(x;\mu')) \int_0^1 (1-t)(x-y)^2 \cdot \partial^2 \psi(y+t(x-y);\mu')\mathrm{d}t\mathrm{d}x,
\end{aligned}$$

where in the second line we used that $K_0$ is translation and reflection invariant, so that $\int_{\mathbb{R}^d} K_0((y,\mu);(x;\mu'))(x-y)\mathrm{d}y = 0$, and where in the last line the tensors $(x-y)^2$ and $\partial^2$ are contracted in the internal indices. We have

$$\begin{aligned}
&\|(\mathbb{1} - \operatorname{Loc}_0) K_0\|_{\mathscr{L}_t^{(1,(2))}(s;c')} \\
&\leqslant 2^{2t} \sup_{y,\mu} \sum_{\mu'} \int_{\mathbb{R}^d} \int_0^1 (1-t)|K_0((y,\mu);(x;\mu'))(x-y)^2| w_{s;c'}^{(2)}(y;y+t(x-y))\mathrm{d}t\mathrm{d}x \\
&\leqslant 2^{2t} \sup_{y,\mu} \sum_{\mu'} \int_{\mathbb{R}^d} \int_0^1 (1-t)|K_0((y,\mu);(x;\mu'))(x-y)^2| w_{s;c'}^{(2)}(y;x)\mathrm{d}t\mathrm{d}x \\
&\lesssim 2^{2(t-s)} \sup_{y,\mu} \sum_{\mu'} \int_{\mathbb{R}^d} |K_0((y,\mu);(x;\mu'))| w_{s;c}^{(2)}(y;x)\mathrm{d}x \\
&= 2^{2(t-s)} \|K_0\|_{\mathscr{L}_t^{(1,(0))}(s;c)},
\end{aligned} \tag{4.5}$$

where $\lesssim$ is up to constants depending on $c - c'$, implying the claim. □

We conclude this section by introducing the Laplacian and the contraction operators. Loosely speaking, these operators implement the action of $\frac{1}{2}\mathrm{D}_{\dot{G}_s^\varepsilon}^2$ and of $\langle \dot{G}_s^\varepsilon \cdot, \mathrm{D} \cdot \rangle$ on the space of monomials described above. In the definition below, we think of $\dot{G}_s^\varepsilon$ as an element of $C^\infty(\mathbb{R}^d \times \mathbb{R}^d; \mathbb{C}^{4 \times 4})$ and thus write $\dot{G}_s^{\varepsilon(m)}$ for some $m \in \{0,2\}^2$.

**Definition 4.8.** *The Laplacian operator is the linear operator $\mathbb{D}_{\dot{G}_s^\varepsilon}^2 \colon \mathscr{L}_t^{(k+2)} \to \mathscr{L}_t^{(k)}$ for any $k \geqslant 1$ that acts as follows: if $K \in \mathscr{L}_t^{(k+2,m)}$ and $m' \subset m$, with $|m'| = k$*

$$(\mathbb{D}_{\dot{G}_s^\varepsilon}^2 K)_{m'}(\psi^{(m')}) := \sum_{m'' \subseteq m} \sum_{\pi \in S_{m',m''}^m} s(\pi) K_m\bigl(P_\pi \psi^{(m')} \otimes \dot{G}_s^{\varepsilon(m'')}\bigr).$$

*where, $s(\pi) \in \{\pm 1\}$ being a sign due to the anticommutativity of the fields. The contraction operator is the bilinear operator $\mathbb{C}_{\dot{G}_s^\varepsilon} \colon \mathscr{L}_t^{(k+1)} \times \mathscr{L}_t^{(k')} \to \mathscr{L}_t^{(k+k')}$, for any $k, k' \geqslant 1$ that acts as follows: if $K \in \mathscr{L}_t^{(k+1,m)}$, $K' \in \mathscr{L}_t^{(k',m')}$ and $m'' \subset m$, with $|m''| = k$*

$$\begin{aligned}
&(\mathbb{C}_{\dot{G}_s^\varepsilon}(K, K'))_{m'',m'}(\psi^{(m'' \circ m')}) \\
&:= \sum_{\bar{m} \subseteq m} \sum_{\pi \in S_{m'',\bar{m}}^m} s'(\pi) \sum_{\mu'} \int K_m\bigl(P_\pi \psi^{(m'')} \otimes \dot{G}_s^{\varepsilon(\bar{m})}((y';\mu'),\cdot)\bigr) \cdot K'_{m'}(\psi^{(m')})((y';\mu'))\mathrm{d}y'
\end{aligned}$$

*where $s'(\pi) \in \{\pm 1\}$ is another suitable sign related to the permutation $\pi$.*

In view of the application to the flow equation the following lemma is crucial.

**Lemma 4.9.** *If $K \in \mathscr{L}_t^{(k+2)}(s;c)$, then $\mathbb{D}_{\dot{G}_s^\varepsilon}^2(K) \in \mathscr{L}_t^{(k)}(s;c)$ and*

$$\|\mathbb{D}_{\dot{G}_s^\varepsilon}^2(K)\|_{\mathscr{L}_t^{(k)}(s;c)} \lesssim k^2 \|K\|_{\mathscr{L}_t^{(k+2)}(s;c)} \|\dot{G}_s^\varepsilon\|_{L^\infty}. \tag{4.6}$$

*If $K \in \mathscr{L}_t^{(k+1)}(s;c)$ and $K' \in \mathscr{L}_t^{(k')}(s;c)$, then $\mathbb{C}_{\dot{G}_s^\varepsilon}(K,K') \in \mathscr{L}_t^{(k+k')}(s;c)$*

$$\|\mathbb{C}_{\dot{G}_s^\varepsilon}(K,K')\|_{\mathscr{L}_t^{(k+k')}(s;c)} \lesssim k \|K\|_{\mathscr{L}_t^{(k+1)}(s;c)} \|K'\|_{\mathscr{L}_t^{(k')}(s;c)} \|\dot{G}_s^\varepsilon\|_{L_{s;c}^1}. \tag{4.7}$$

**Proof.** The bound on the Laplacian operator follows by noting that there are at most $(k+1)(k+2)$ elements in the sum and by noting that

$$\begin{aligned}
\|K_m(P_\pi(\cdot \dot{G}_s^{\varepsilon(m'')}))\|_{\mathscr{L}_t^{(k,m')}(s;c)} &= \||K_m(P_\pi(\cdot \dot{G}_s^{\varepsilon(m'')}))| \cdot w_{s;c}^{(k+1)}\|_{\mathscr{L}_t^{(k,m')}} \\
&\leqslant 2^{-[m'']t} \|\dot{G}_s^{\varepsilon(m'')}\|_{L^\infty} \||K_m| \cdot w_{s;c}^{(k+3)}\|_{\mathscr{L}_t^{(k+2,m)}} \\
&\lesssim 2^{[m''](s-t)} \|\dot{G}_s^\varepsilon\|_{L^\infty} \|K_m\|_{\mathscr{L}_t^{(k+2,m)}(s;c)},
\end{aligned}$$

where in the last line we used that $\|\dot{G}_s^{\varepsilon(m'')}\|_{L^\infty} \lesssim 2^{[m'']s} \|\dot{G}_s^\varepsilon\|_{L^\infty}$, see Corollary 3.24, and that the Steiner diameter is a non-decreasing function in the cardinality of the set. To prove the bound on the contraction operator, we observe that there are at most $k+1$ terms in the sum and that

$$\begin{aligned}
&\left\|\int K_m(P_\pi(\cdot \dot{G}_s^{\varepsilon(\bar{m})}((y,\mu),\cdot))) K'_{m'}((y,\mu)) dy\right\|_{\mathscr{L}_t^{(k+k',m''\circ m')}(s;c)} \\
&= \left\|\left|\int K_m(P_\pi(\cdot \dot{G}_s^{\varepsilon(\bar{m})}((y,\mu),\cdot))) K'_{m'}((y,\mu)) dy\right| \cdot w_{s;c}^{(k+k'+1)}\right\|_{\mathscr{L}_t^{(k+k',m''\circ m')}} \\
&\leqslant \left\|\int |K_m(P_\pi(\cdot \dot{G}_s^{\varepsilon(\bar{m})}((y,\mu),\cdot)))| dy \cdot w_{s;c}^{(k+2)}(y;\cdot)\right\|_{\mathscr{L}_t^{(k,m'')}} \|K'_{m'} \cdot w_{s;c}^{(k'+1)}\|_{\mathscr{L}_t^{(k',m')}} \\
&\leqslant 2^{-[\bar{m}]t} \|\dot{G}_s^{\varepsilon(\bar{m})}\|_{L_{s;c}^1} \|K_m \cdot w_{s;c}^{(k+2)}\|_{\mathscr{L}_t^{(k+1,m'')}} \|K'_{m'} \cdot w_{s;c}^{(k'+1)}\|_{\mathscr{L}_t^{(k',m')}} \\
&\lesssim 2^{[\bar{m}](s-t)} \|\dot{G}_s^\varepsilon\|_{L_{s;c}^1} \|K_m \cdot w_{s;c}^{(k+2)}\|_{\mathscr{L}_t^{(k+1,m)}} \|K'_{m'} \cdot w_{s;c}^{(k'+1)}\|_{\mathscr{L}_t^{(k',m')}},
\end{aligned}$$

where in the third line we used that $\mathrm{St}(\{x,x_1,\ldots,x_k\} \cup \{x'_1,\ldots,x'_{k'}\}) \leqslant \mathrm{St}(\{x,x_1,\ldots,x_k,y\}) + \mathrm{St}(\{y,x'_1,\ldots,x'_{k'}\})$ and in the fourth line that $\mathrm{St}(\{x,x_1,\ldots,x_k,y\}) \leqslant \mathrm{St}(\{x,x_1,\ldots,x_k\}) + |x_k-y|$. $\square$

**Remark 4.10.** Note that the scale dependent weight in the norm $\|\cdot\|_{\mathscr{L}_t^{(k+k')}(s;c)}$, see (4.1) are crucial for these bounds to hold. In fact, when the derivatives act on $\dot{G}_s^\varepsilon$ one loses $2^s$ for each derivative, but this factor is then cancelled by the gain in the kernel, i.e., due to the weight in (4.1). The following corollary is an immediate consequence of Lemma 4.9 and Lemma 4.7.

**Corollary 4.11.** *If $K \in \mathscr{L}_t^{(3)}(s;c)$, then for any $c' < c$, $\mathrm{Ren}\,\mathbb{D}_{\dot{G}_s^\varepsilon}^2(K) \in \mathscr{L}_t^{(1)}(s;c')$ and*

$$\|(\mathrm{Ren}\,\mathbb{D}_{\dot{G}_s^\varepsilon}^2(K))_2\|_{\mathscr{L}_t^{(1,(2))}(s;c')} \lesssim 2^{-2(s-t)} \|K\|_{\mathscr{L}_t^{(3)}(s;c)} \|\dot{G}_s^\varepsilon\|_{L_{s;c}^\infty}.$$

*If $K, K' \in \mathscr{L}_t^{(1)}(s;c)$, then for any $c' < c$, $\mathrm{Ren}\,\mathbb{C}_{\dot{G}_s^\varepsilon}(K,K') \in \mathscr{L}_t^{(1)}(s;c')$ and*

$$\|(\mathrm{Ren}\,\mathbb{C}_{\dot{G}_s^\varepsilon}^2(K,K'))_2\|_{\mathscr{L}_t^{(1,(2))}(s;c')} \lesssim_{c-c'} 2^{-2(s-t)} \|K\|_{\mathscr{L}_t^{(1)}(s;c)} \|K'\|_{\mathscr{L}_t^{(1)}(s;c)} \|\dot{G}_s^\varepsilon\|_{L_{s;c}^1}.$$

### 4.2 The truncated flow equation

For technical convenience, we study the flow equation by introducing a further grading with respect to the derivative in the flow parameter. Thus, we decompose the function $(F_s^\varepsilon)_s$ as follows

$$F_s^\varepsilon(\psi_t) = \sum_{\ell \geqslant 0} F_s^{\varepsilon[\ell]}(\psi_t) = \sum_{\ell \geqslant 0, k \geqslant 1} F_s^{\varepsilon[\ell](k)}(\psi_t), \tag{4.8}$$

where the label $\ell$ is the grading with respect to the operator $\partial_s$, whereas $k$ is the grading with respect to the functional derivative D, so that $F_s^{\varepsilon[\ell](k)} \in \mathscr{L}_s^{(k)}$ are monomials of degree $k$. This grading decomposition is natural because the flow equation involves both differential operators. We introduce the following truncated flow equation

$$\partial_s F_s^\varepsilon(\psi_s) + \frac{1}{2} \mathrm{D}_{\dot G_s^\varepsilon}^2 F_s^\varepsilon(\psi_s) + \Pi_{\leqslant n} \mathrm{D} F_s^\varepsilon(\psi_s) \cdot \dot G_s^\varepsilon F_s^\varepsilon(\psi_s) = 0, \qquad (4.9)$$

where $\Pi_{\leqslant n}$ is the projection on terms with $\ell \leqslant n$. For convenience, we also introduce the notation

$$H_{>n;t}^\varepsilon(\psi_t) := \Pi_{>n} \mathrm{D} F_t^\varepsilon(\psi_t) \cdot \dot G_t^\varepsilon F_t^\varepsilon(\psi_t). \qquad (4.10)$$

The truncation is convenient because it allows us to look for polynomial solutions. Note that we could truncate with respect to both gradings, with no major difference in the analysis. Expanding as in (4.8), the truncated flow equation reads, for $\ell = 0, \ldots, n$:

$$\partial_s F_s^{\varepsilon[\ell+1](k)} + \mathbb{D}_{\dot G_s^\varepsilon}^2 (F_s^{\varepsilon[\ell](k+2)}) + \sum_{\ell'=0}^{\ell} \sum_{k'=0}^{k-1} \mathbb{C}_{\dot G_s^\varepsilon} (F_s^{\varepsilon[\ell'](k-k')}, F_s^{\varepsilon[\ell-\ell'](k'+1)}) = 0, \qquad (4.11)$$

where the operators $\mathbb{D}_{\dot G_s^\varepsilon}^2$ and $\mathbb{C}_{\dot G_s^\varepsilon}$ are as in Definition 4.8. Note that eq. (4.11) is triangular in $\ell$ and this will bring technical simplification to the analysis.

We let

$$c_\ell := (1 - \ell n^{-1}) \bar c / 2, \qquad (4.12)$$

$\bar c$ being the constant in Lemma 3.23. We make the following Ansatz for the norm of the monomials $F_s^{\varepsilon[\ell](k)}$, for suitable constants $C > 0$, $\alpha, \beta, \kappa \geqslant 0$:

$$\|F_s^{\varepsilon[\ell](k)}\|_{\mathscr{L}_t^{(k)}(s;c_\ell)} \lesssim C^k 2^{(\alpha - \beta k - \kappa \ell)s}, \qquad 0 \leqslant t \leqslant s. \qquad (4.13)$$

If we interpret the exponent multiplying $s$ as the regularity, we have that $\beta$ quantifies the gain in regularity due to the grading in the fields $k$ whereas $\kappa$ the gain due to the perturbative order in $\ell$. Because the problem is super-renormalisable, we will have $\beta > \gamma$, $-\gamma$ being the regularity of the fields and furthermore $\kappa > 0$. Now note that because of the "initial condition" $F_s^{[0]}(\psi_s^\varepsilon) = \lambda(\psi_s^\varepsilon)^3$, the Ansatz (4.13) can hold true only if

$$\alpha - 3\beta \geqslant 0. \qquad (4.14)$$

Furthermore, as shall see shortly, these estimates can be propagated only if the following conditions are satisfied:

$$5\beta - \alpha - 2\gamma - \kappa > 0, \qquad d + 4\beta - 2\alpha - 2\gamma - \kappa > 0, \qquad (4.15)$$

and

$$3\beta + 2 - \alpha - 2\gamma - \kappa > 0, \qquad d + 2\beta + 2 - 2\alpha - 2\gamma - \kappa > 0 \qquad (4.16)$$

The constraints in (4.15) makes irrelevant the terms with $k \geqslant 3$, in the sense of integrability at $s = \infty$ of the Laplacian and contraction operators respectively. On the other hand, the constraints in (4.16) make irrelevant the renormalised term with $k = 1$.

**Remark 4.12.** Note that the constraints (4.14), (4.15) and (4.16) can be satisfied only provided that

$$\begin{aligned} \gamma < \min\{d/4, 1\}, \qquad \beta \in (\gamma, \min\{d/2 - \gamma, d/4 + 1/2 - \gamma/2\}), \\ \alpha \in [3\beta, d/2 + \min\{2\beta, \beta + 1\} - \gamma). \end{aligned} \qquad (4.17)$$

The condition $\gamma < d/4$ is what makes the model super-renormalisable. As hinted in the introduction, the constraint $\gamma < 1$ comes from the requirement that the counter-term is a mass term and does not have derivatives of the fields. Note that the super-renormalisability of the model allows for room in the choice of the parameters $\beta, \alpha$ and $\kappa$.

Given the Ansatz (4.13), the following norm topology is natural: for any graded sequence $F = (F_s^{[\ell](k)})_{s \geq 0}^{\ell \geq 0, k \geq 0}$, we introduce the norm

$$\|F\|_{C,\alpha,\beta,\kappa} := \sup_{\ell \geq 0, k \geq 0} \sup_{s,t: 0 \leq t \leq s} C^{-k} 2^{-(\alpha - \beta k - \kappa \ell)s} \|F_s^{[\ell](k)}\|_{\mathscr{L}_t^{(k)}(s;c_\ell)}, \qquad (4.18)$$

where for simplicity we dropped the dependence on $(c_\ell)_\ell$. The following theorem establishes the existence and uniqueness of a global solution of the truncated flow equations.

**Theorem 4.13.** *(Global solution). Let $\varepsilon \in 2^{-\mathbb{N}_\infty}$, $n \in \mathbb{N}$, $C > 0$, let $\alpha, \beta, \gamma, \kappa \geq 0$ satisfy the constraints in (4.14), (4.15), (4.16) and let $c_\ell$ be as in (4.12). The truncated flow equations (4.11) with boundary data*

$$F_0^{\varepsilon[0](k)} = K^{(k)} \mathbf{1}_{k \in \{1,3\}}, \qquad \sup_{k \geq 0} \sup_{0 \leq t \leq s} C^{-k} 2^{-(\alpha - \beta k)s} \|K^{(k)}\|_{\mathscr{L}_t^{(k)}(s;c_0)} < \infty,$$

*and for $\ell \geq 1$*

$$F_\infty^{\varepsilon[\ell](k)} = 0 \quad \forall k \neq 1, \qquad \operatorname{Ren} F_\infty^{\varepsilon[\ell](1)} = 0, \qquad \operatorname{Loc} F_0^{\varepsilon[\ell](1)} = 0$$

*have a unique global solution $F^\varepsilon = (F_s^{\varepsilon[\ell](k)} \mathbf{1}_{k \text{ odd}})_{s \geq 0}^{\ell \geq 0, k \geq 0}$, such that $F_s^{\varepsilon[\ell](k)} \in \mathscr{L}_t^{(k)}(s; c_\ell)$.*

**Remark 4.14.** Note that the monomials $K^{(k)}$ have to be local in order to satisfy the said bound.

**Proof.** The proof is by induction in $\ell$. We let $F^{\varepsilon[\leq \ell]}$ denote the graded sequence $\left(F_t^{\varepsilon[\ell'](k)} \mathbf{1}_{\ell' \leq \ell}\right)_{t \geq 0}^{\ell' \geq 0, k \geq 0}$ and show that $\|F^{\varepsilon[\leq \ell]}\|_{C,\alpha,\beta,\kappa}$ is finite.

The validity at $\ell = 0$ is obvious because $\partial_s F_s^{\varepsilon[0]} = 0$ and because of the assumptions. In particular, we have that $\|F^{\varepsilon[\leq 0]}\|_{C,\alpha,\beta,\kappa} < \infty$. To prove the induction step, we integrate the truncated flow equation (4.11) for $k \geq 3$ backwards from the final condition at $s = \infty$ with $F_\infty^{\varepsilon[\ell+1](k)} = 0$ for any $\ell \geq 0, k \geq 3$. We have:

$$F_t^{\varepsilon[\ell+1](k)} = \int_t^\infty \mathbb{D}_{\dot{G}_s^\varepsilon}(F_s^{\varepsilon[\ell](k+2)}) ds$$
$$+ \int_t^\infty \sum_{\ell'=0}^{\ell} \sum_{k'=0}^{k-1} \mathbb{C}_{\dot{G}_s^\varepsilon}(F_s^{\varepsilon[\ell'](k-k')}, F_s^{\varepsilon[\ell-\ell'](k'+1)}) ds$$

For the term with $k = 1$, we use the equivalence $F_t^{\varepsilon[\ell+1](1)} \sim \operatorname{Loc} F_t^{\varepsilon[\ell+1](1)} + \operatorname{Ren} F_t^{\varepsilon[\ell+1](1)}$, that is, $F_t^{\varepsilon[\ell+1](1)} \sim ((\operatorname{Loc} F_t^{\varepsilon[\ell+1](1)})_0, (\operatorname{Ren} F_t^{\varepsilon[\ell+1](1)})_2)$. We integrate the term with the renormalisation from $s = \infty$, with $(\operatorname{Ren} F_\infty^{\varepsilon[\ell+1](1)})_2 = 0$ for any $\ell \geq 0$

$$(\operatorname{Ren} F_t^{\varepsilon[\ell+1](1)})_2 = \int_t^\infty (\operatorname{Ren} \mathbb{D}_{\dot{G}_s^\varepsilon}^2(F_s^{\varepsilon[\ell](3)}))_2 ds$$
$$+ \int_t^\infty \sum_{\ell'=0}^{\ell} (\operatorname{Ren} \mathbb{C}_{\dot{G}_s^\varepsilon}(F_s^{\varepsilon[\ell'](1)}, F_s^{\varepsilon[\ell-\ell'](1)}))_2 ds,$$

and, finally, we integrate the term with the localisation from 0 to $t$, with $(\operatorname{Loc} F_t^{\varepsilon[\ell+1](1)})_0 = 0$ for any $\ell \geq 0$

$$(\operatorname{Loc} F_t^{\varepsilon[\ell+1](1)})_0 = -\int_0^t (\operatorname{Loc} \mathbb{D}_{\dot{G}_s^\varepsilon}^2(F_s^{\varepsilon[\ell](3)}))_0 ds$$
$$- \int_0^t \sum_{\ell'=0}^{\ell} (\operatorname{Loc} \mathbb{C}_{\dot{G}_s^\varepsilon}(F_s^{\varepsilon[\ell'](1)}, F_s^{\varepsilon[\ell-\ell'](1)}))_0 ds.$$

Now, recalling the kernel estimates for the propagator $\|\dot{G}_s^\varepsilon\|_{L^p_{s;c}} \lesssim 2^{-(d/p-2\gamma)s}$, for $c < \bar{c}$, see Corollary 3.24, Lemma 4.9 implies the following estimates:

$$\|\mathbb{D}^2_{\dot{G}_s^\varepsilon}(F_s^{\varepsilon[\ell](k+2)})\|_{\mathscr{L}_t^{(k)}(s;c_\ell)}$$
$$\lesssim k^2 C^{k+2} 2^{(\alpha-\beta k)s} 2^{2(\gamma-\beta)s} 2^{-\kappa \ell s} \|F^{\varepsilon[\leqslant \ell]}\|_{C,\alpha,\beta,\kappa}, \tag{4.19}$$

and for $\ell' = 0, \ldots, \ell$

$$\|\mathbb{C}_{\dot{G}_s^\varepsilon}(F_s^{\varepsilon[\ell'](k-k')}, F_s^{\varepsilon[\ell-\ell'](k'+1)})\|_{\mathscr{L}_t^{(k+k')}(s;c_\ell)}$$
$$\lesssim k C^{k+1} 2^{(\alpha-\beta k)s} 2^{-(d+\beta-\alpha-2\gamma)s} 2^{-\kappa \ell s} \|F^{\varepsilon[\leqslant \ell]}\|^2_{C,\alpha,\beta,\kappa}. \tag{4.20}$$

By using these bounds, the constraints (4.17) and Remark 4.3, we obtain that for $k \geqslant 3$ and $r \leqslant t$ $F_t^{\varepsilon[\ell+1](k)} \in \mathscr{L}_r^{(k)}(t;c_{\ell+1})$ and

$$\|F_t^{\varepsilon[\ell+1](k)}\|_{\mathscr{L}_r^{(k)}(t;c_{\ell+1})}$$
$$\leqslant \int_t^\infty \|\mathbb{D}^2_{\dot{G}_s^\varepsilon}(F_s^{\varepsilon[\ell](k+2)})\|_{\mathscr{L}_r^{(k)}(s;c_\ell)} ds$$
$$+ \int_t^\infty \sum_{\ell'=0}^\ell \sum_{k'=0}^{k-2} \|\mathbb{C}_{\dot{G}_s^\varepsilon}(F_s^{\varepsilon[\ell'](k-k')}, F_s^{\varepsilon[\ell-\ell'](k'+1)})\|_{\mathscr{L}_r^{(k)}(s;c_\ell)} ds$$
$$\lesssim k^2 C^{k+2} \|F^{\varepsilon[\leqslant \ell]}\|_{C,\alpha,\beta,\kappa} \int_t^\infty 2^{(\alpha-\beta k)s} 2^{2(\gamma-\beta)s} 2^{-\kappa \ell s} ds \tag{4.21}$$
$$+ (\ell+1) k C^{k+1} \|F^{\varepsilon[\leqslant \ell]}\|^2_{C,\alpha,\beta,\kappa} \int_t^\infty 2^{(\alpha-\beta k)s} 2^{-(d+\beta-\alpha-2\gamma)s} 2^{-\kappa \ell s} ds$$
$$\lesssim C^{k+1} 2^{(\alpha-\beta k)t} 2^{-\kappa(\ell+1)t} (\|F^{\varepsilon[\leqslant \ell]}\|_{C,\alpha,\beta,\kappa} + \|F^{\varepsilon[\leqslant \ell]}\|^2_{C,\alpha,\beta,\kappa}),$$

where in the last line $\lesssim$ is up to constants depending on $\ell, \alpha, \beta, \gamma, \kappa$. In a similar way, Corollary 3.24, Corollary 4.11 with $c' \coloneqq c_{\ell+1}$ and $c \coloneqq c_\ell$, and the constraints (4.17), we obtain that for $r \leqslant t$ $(\operatorname{Ren} F_t^{\varepsilon[\ell+1](1)})_2 \in \mathscr{L}_r^{(1,(2))}(t;c_{\ell+1})$ and

$$\|(\operatorname{Ren} F_t^{\varepsilon[\ell+1](1)})_2\|_{\mathscr{L}_r^{(1,(2))}(t;c_{\ell+1})}$$
$$\leqslant \int_t^\infty \|(\operatorname{Ren} \mathbb{D}^2_{\dot{G}_s^\varepsilon}(F_s^{\varepsilon[\ell](3)}))_2\|_{\mathscr{L}_r^{(1,(2))}(s;c_{\ell+1})} ds$$
$$+ \int_t^\infty \sum_{\ell'=0}^\ell \|(\operatorname{Ren} \mathbb{C}_{\dot{G}_s^\varepsilon}(F_s^{\varepsilon[\ell'](1)}, F_s^{\varepsilon[\ell-\ell'](1)}))_2\|_{\mathscr{L}_r^{(1,(2))}(s;c_{\ell+1})} ds \tag{4.22}$$
$$\lesssim C^3 \|F^{\varepsilon[\leqslant \ell]}\|_{C,\alpha,\beta,\kappa} 2^{2r} \int_t^\infty 2^{(\alpha-\beta-2)s} 2^{2(\gamma-\beta)s} 2^{-\kappa \ell s} ds$$
$$+ (\ell+1) C^2 \|F^{\varepsilon[\leqslant \ell]}\|^2_{C,\alpha,\beta,\kappa} 2^{2r} \int_t^\infty 2^{(\alpha-\beta-2)s} 2^{-(d+\beta-\alpha-2\gamma)s} 2^{-\kappa \ell s} ds$$
$$\lesssim C^2 2^{2(r-t)} 2^{(\alpha-\beta)t} 2^{-\kappa(\ell+1)t} (\|F^{\varepsilon[\leqslant \ell]}\|_{C,\alpha,\beta,\kappa} + \|F^{\varepsilon[\leqslant \ell]}\|^2_{C,\alpha,\beta,\kappa}).$$

Then, we control the term with the localisation operator, for $r \leqslant t$

$$\|(\operatorname{Loc} F_t^{\varepsilon[\ell+1](1)})_0\|_{\mathscr{L}_r^{(1,(0))}(t;c_{\ell+1})}$$
$$\leqslant \int_0^t \|\operatorname{Loc} \mathbb{D}^2_{\dot{G}_s^\varepsilon}(F_s^{\varepsilon[\ell](3)})\|_{\mathscr{L}_0^{(1)}(s;c_\ell)} ds$$
$$+ \sum_{\ell'=0}^\ell \int_0^t \|\operatorname{Loc} \mathbb{C}_{\dot{G}_s^\varepsilon}(F_s^{\varepsilon[\ell'](1)}, F_s^{\varepsilon[\ell-\ell'](1)})\|_{\mathscr{L}_0^{(1)}(s;c_\ell)} ds \tag{4.23}$$
$$\lesssim C^3 \|F^{\varepsilon[\leqslant \ell]}\|_{C,\alpha,\beta,\kappa} \int_0^t 2^{(\alpha-\beta)s} 2^{2(\gamma-\beta)s} 2^{-\kappa \ell s} ds$$
$$+ (\ell+1) C^2 \|F^{\varepsilon[\leqslant \ell]}\|^2_{C,\alpha,\beta,\kappa} \int_0^t 2^{(\alpha-\beta)s} 2^{-(d+\beta-\alpha-2\gamma)s} 2^{-\kappa \ell s} ds$$
$$\lesssim C^2 2^{(\alpha-\beta)t} 2^{-\kappa(\ell+1)t} (\|F^{\varepsilon[\leqslant \ell]}\|_{C,\alpha,\beta,\kappa} + \|F^{\varepsilon[\leqslant \ell]}\|^2_{C,\alpha,\beta,\kappa}),$$

where we used that $\|(K)_0\|_{\mathscr{L}_r^{(1,(0))}(t;c_{\ell+1})} \leqslant \|K\|_{\mathscr{L}_r^{(1)}(t;c_{\ell+1})}$, and the property (4.3) to set the scale of the field to zero. Putting together (4.22) and (4.23), we obtain for $r \leqslant t$

$$\|F_t^{\varepsilon[\ell+1](1)}\|_{\mathscr{L}_r^{(1)}(t;c_{\ell+1})} \lesssim C^2 2^{(\alpha-\beta)t} 2^{-\kappa(\ell+1)t}(\|F^{\varepsilon[\leqslant \ell]}\|_{C,\alpha,\beta,\kappa} + \|F^{\varepsilon[\leqslant \ell]}\|^2_{C,\alpha,\beta,\kappa}) \tag{4.24}$$

The bounds (4.21), (4.24) imply that $\|F^{\varepsilon[\leqslant \ell+1]}\|_{C,\alpha,\beta,\kappa}$ is finite. This concludes the proof. □

**Corollary 4.15.** *Let $F^\varepsilon$ be the unique solution to the flow equation truncated at n, with the same boundary data as in Theorem 4.13, but in particular $F_0^{\varepsilon[0]}(\psi_t) := \lambda(\psi_t)^3$, see eq. (3.27). Let $H_{>n}^\varepsilon$ be defined as in (4.10). Then, there exist constants $C_0 = C_0(\alpha,\beta,\gamma,\kappa,n)$ and $\lambda_0 = \lambda_0(C)$ independent of $\varepsilon$ and L such that, if $C \leqslant C_0$ and $\lambda \leqslant \lambda_0(C)$ then*

$$\|F^\varepsilon\|_{C,\alpha,\beta,\kappa} \lesssim \lambda, \qquad \|H_{>n}^\varepsilon\|_{C,\alpha,\beta,\kappa} \lesssim_n \lambda^2,$$

*uniformly in $\varepsilon$ and L, where the $\lesssim$ are up to universal constants.*

**Proof.** The proof is a corollary of the proof of Theorem 4.13. In fact, one can read off the bound for $\|F^\varepsilon\|_{C,\alpha,\beta,\kappa}$ directly from (4.21) (which holds for any $k \geqslant 0$), provided that $C$ and $\lambda$ are small enough, the latter depending on $C$. Finally, comparing with (4.11) we have $H_{>n,s}^{\varepsilon[\ell](k)} = \mathbf{1}_{n<\ell\leqslant 2n} \sum_{\ell'=0}^{\ell} \sum_{k'=0}^{k-1} \mathbb{C}_{\dot{G}_s^\varepsilon}(F_s^{\varepsilon[\ell'](k-k')}, F_s^{\varepsilon[\ell-\ell'](k'+1)})$, and therefore, by the proof of Theorem 4.13 we have the bound $\|H_{>n;s}^{\varepsilon[\ell](k)}\|_{\mathscr{L}_t^{(k)}(s;c_\ell)} \lesssim k(\ell+1)C^{k+1}2^{(\alpha-\beta k)s}2^{-(d+\beta-\alpha-2\gamma)s}2^{-\kappa \ell s}\|F\|^2_{C,\alpha,\beta,\kappa}$, which implies the claim. □

**Remark 4.16.** We first of all note that because of the symmetries in the internal component, one can show that for any $\mu,\mu',\bar{\mu} \in \{\uparrow,\downarrow\} \times \{\pm\}$

$$F_t^{\varepsilon[\ell](1)}((y;\mu);(x;\mu')) = \delta_{\mu,\mu'} F_t^{\varepsilon[\ell](1)}((y;\bar{\mu});(x;\bar{\mu})),$$

that is, $(F_t^{\varepsilon[\ell](1)}((y;\mu);(x;\mu')))_{\mu,\mu'}$ is diagonal and constant.

Given the global solution, and taking Remark 4.16 and eq. (3.30) into account, we can finally define the renormalised chemical potential.

**Definition 4.17.** *The scale-dependent chemical potential $\mu_t^\varepsilon$ associated with the global solution $F^\varepsilon$ of the truncated flow equations is defined as, for any $\mu \in \{\uparrow,\downarrow\} \times \{\pm\}$*

$$\mu_t^\varepsilon := -\text{sign}(\mu) \sum_{\ell=0}^{n} \int_{\mathbb{R}^d} (F_t^{\varepsilon[\ell](1)}((0;\mu);(x;\mu)))_0 \, dx,$$

*where $F_t^{\varepsilon[\ell](1)} = ((F_t^{\varepsilon[\ell](1)})_0, (F_t^{\varepsilon[\ell](1)})_2)$, $(F_t^{\varepsilon[\ell](1)})_0$ being the local part, and $\text{sign}(\mu)$ is defined in (3.28). The renormalised chemical potential is $\mu_\infty^\varepsilon := \lim_{T\to\infty} \mu_T^\varepsilon$.*

**Remark 4.18.** Note the bound

$$|\mu_T^\varepsilon| \leqslant \sum_{\ell=0}^{n} \|\text{Loc}\, F_T^{\varepsilon[\ell](1)}\|_{\mathscr{L}_T^{(1)}} \lesssim_n \|F^\varepsilon\|_{C,\alpha,\beta,\kappa} 2^{(\alpha-\beta-\kappa)T}$$

where by (4.14) and the second constraint in (4.15) $\alpha - \beta - \kappa > 2\gamma$.

Next, we prove the convergence $F^\varepsilon \to F$ in a weaker topology.

**Proposition 4.19.** *(Continuum limit). Let $\varepsilon > \varepsilon' \in 2^{-\mathbb{N}_\infty}$. Under the same assumptions of Theorem 4.13 but with $\alpha > 3\beta$, let $F^\varepsilon$ and $F^{\varepsilon'}$ be the global solutions of the truncated flow equation and let $\theta \in \left(0, \frac{\alpha-3\beta}{2}\right]$. Then*

$$\|F^\varepsilon - F^{\varepsilon'}\|_{C,\alpha,\beta,\kappa} \lesssim_{n,\theta} \varepsilon^\theta \sup_{\#=\varepsilon,\varepsilon'} (\|F^\#\|_{C,\alpha-\theta,\beta,\kappa} + \|F^\#\|^{n+2}_{C,\alpha-\theta,\beta,\kappa}).$$

**Proof.** Set

$$\delta_{\varepsilon,\varepsilon'} F_s := F_s^\varepsilon - F_s^{\varepsilon'}, \qquad \delta_{\varepsilon,\varepsilon'} \dot{G}_s := \dot{G}_s^\varepsilon - \dot{G}_s^{\varepsilon'}. \qquad (4.25)$$

We can write the truncated flow equation for the difference as follows:

$$\partial_s \delta_{\varepsilon,\varepsilon'} F_s^{[\ell+1](k)} = \delta_{\varepsilon,\varepsilon'} \mathbb{D}^2_{\dot{G}_s}(F_s^{[\ell](k+2)}) + \sum_{\ell'=0}^{\ell}\sum_{k'=1}^{k} \delta_{\varepsilon,\varepsilon'} \mathbb{C}_{\dot{G}_s}\bigl(F_s^{[\ell'](k')}, F_s^{[\ell-\ell'](k+1-k')}\bigr) \qquad (4.26)$$

with identically null boundary data $\delta_{\varepsilon,\varepsilon'} F_0^{[0]} = 0$. Above, the difference operator $\delta_{\varepsilon,\varepsilon'}$ acts on $\dot{G}_s$ as well,

$$\delta_{\varepsilon,\varepsilon'} \mathbb{D}^2_{\dot{G}_s}(F_s^{[\ell](k)}) = \mathbb{D}^2_{\dot{G}_s^\varepsilon}(\delta_{\varepsilon,\varepsilon'} F_s^{[\ell](k)}) + \mathbb{D}^2_{\delta_{\varepsilon,\varepsilon'} \dot{G}_s}(F_s^{\varepsilon'[\ell](k)})$$

and similar identity holds for $\delta_{\varepsilon,\varepsilon'} \mathbb{C}_{\dot{G}_s}$

$$\delta_{\varepsilon,\varepsilon'} \mathbb{C}_{\dot{G}_s}\bigl(F_s^{[\ell](k)}, F_s^{[\ell'](k')}\bigr) = \mathbb{C}_{\dot{G}_s^\varepsilon}\bigl(\delta_{\varepsilon,\varepsilon'} F_s^{[\ell](k)}, F_s^{\varepsilon[\ell'](k')}\bigr) + \mathbb{C}_{\dot{G}_s^\varepsilon}\bigl(F_s^{\varepsilon'[\ell](k)}, \delta_{\varepsilon,\varepsilon'} F_s^{[\ell'](k')}\bigr)$$
$$+ \mathbb{C}_{\delta_{\varepsilon,\varepsilon'} \dot{G}_s}\bigl(F_s^{\varepsilon'[\ell](k)}, F_s^{\varepsilon'[\ell'](k')}\bigr)$$

To estimate $\|\delta_{\varepsilon,\varepsilon'} F\|_{C,\alpha,\beta-\theta,\kappa} := \|F^\varepsilon - F^{\varepsilon'}\|_{C,\alpha,\beta-\theta,\kappa}$, we proceed as in the proof of Theorem 4.13, the crucial difference being that $\delta_{\varepsilon,\varepsilon'} F_0^{[0]} = 0$ and the estimate $\|\delta_{\varepsilon,\varepsilon'} \dot{G}_t\|_{L^p_{t;c}} \lesssim 2^{-(d/p-2\gamma-\theta)t} \varepsilon^\theta$ see Corollary 3.24. In particular, one wants to prove inductively that

$$\|\delta_{\varepsilon,\varepsilon'} F^{[\leq\ell]}\|_{C,\alpha,\beta,\kappa} \lesssim_\ell \varepsilon^\theta \sup_{\#=\varepsilon,\varepsilon'} (\|F^{\#[\leq\ell]}\|_{C,\alpha-\theta,\beta,\kappa} + \|F^{\#[\leq\ell]}\|_{C,\alpha-\theta,\beta,\kappa}^{\ell+2}).$$

Repeating the analysis in the proof of Theorem 4.13, we obtain the estimates

$$\|\delta_{\varepsilon,\varepsilon'} F_t^{[\ell+1](k)}\|_{\mathscr{L}_r^{(k)}(t;c_{\ell+1})}$$
$$\lesssim \varepsilon^\theta C^{k+1} 2^{(\alpha-\beta k)t} 2^{-\kappa(\ell+1)t}(\|F^{\#[\leq\ell]}\|_{C,\alpha-\theta,\beta,\kappa} + \|F^{\#[\leq\ell]}\|_{C,\alpha-\theta,\beta,\kappa}^{\ell+3})$$

compare with (4.21), implying the claim. □

## 4.3 Global solution of the FBSDE without cutoffs

In this section, we study the global existence and uniqueness of the pair $(\Psi^{L,\varepsilon}, R^{L,\varepsilon})$ solving the coupled equations

$$\begin{aligned}\Psi_t^{L,\varepsilon} &= -\int_0^t \dot{G}_s^\varepsilon(F_s^\varepsilon(\Psi_s^{L,\varepsilon}) + R_s^\varepsilon)\mathrm{d}s + X_t^{L,\varepsilon},\\ R_t^{L,\varepsilon} &= \int_t^\infty \omega_t(H_{>n;s}^\varepsilon(\Psi_s^{L,\varepsilon}))\mathrm{d}s + \int_t^\infty \omega_t(\mathrm{D}F_s^\varepsilon(\Psi_s^{L,\varepsilon}) \cdot \dot{G}_s^\varepsilon R_s^{L,\varepsilon})\mathrm{d}s,\end{aligned} \qquad (4.27)$$

compare with Lemma 3.22, with input $F_t^\varepsilon$ and $H_{>n;s}^\varepsilon$ graded sequences satisfying suitable smallness assumptions, and then study the limit as the cutoffs $L \to \infty$, $\varepsilon \to 0$.

To control the global solution of eq. (4.27) and its convergence as $L \to \infty$, $\varepsilon \to 0$, we use the following topology: for $\sigma, \sigma' \in \mathbb{R}$, $a, a' > 0$ and $\eta \geq 0$ we define the norm

$$\|(\psi,\psi')\|_* := a^{-1} \sup_{t \in \mathbb{R}^+} 2^{-\sigma t} \|\psi_t\|_{\mathscr{C}_t^{\leq n}(\eta)} + (a')^{-1} \sup_{t \in \mathbb{R}^+} 2^{-\sigma' t} \|\varphi_t\|_{L^\infty(\eta)}, \qquad (4.28)$$

see (3.33) for the definition of $\|\cdot\|_{L^\infty(\eta)}$ and $\|\cdot\|_{\mathscr{C}_t^n(\eta)}$, where we set

$$\|\psi_t\|_{\mathscr{C}_t^{\leq n}(\eta)} := \sum_{n'=0}^n \|\psi_t\|_{\mathscr{C}_t^n(\eta)}.$$

Note that $\|\cdot\|_{\mathscr{C}_t^{\leq n}} \equiv \|\cdot\|_{\mathscr{C}_t^{\leq n}(0)}$ is of course implied.

We can state our result as follows.

**Proposition 4.20.** *Let $L \in \mathbb{N}_\infty$ and $\varepsilon \in 2^{-\mathbb{N}_\infty}$. Let $\alpha, \beta, \gamma, \kappa$ satisfy (4.14), (4.15) and (4.16), let $\theta \in [0, (\beta - \gamma)/2]$ and let $F^\varepsilon = (F_t^{\varepsilon[\ell](k)})_{t \geq 0}^{\ell \geq 0, k \geq 1}$ and $H_{>n}^\varepsilon = (H_{>n;t}^{\varepsilon[\ell](k)})_{t \geq 0}^{\ell > n, k \geq 0}$ be the graded sequences in (4.27), with $n \geq \lfloor \alpha \kappa^{-1} \rfloor + 1$. There exist constants $D, c_n, \lambda_0 = \lambda_0(d, \alpha, \beta, \gamma, n)$ such that, if $\|F^\varepsilon\|_{C,\alpha,\beta,\kappa}, \|H_{>n}^\varepsilon\|_{C,\alpha,\beta,\kappa} \leq \lambda$ for $C < D^{-1}/2$ and $\lambda \leq \lambda_0$, then (4.27) has a unique global solution $(\Psi^{L,\varepsilon}, R^{L,\varepsilon})$ satisfying*

$$\|\|(\Psi^{L,\varepsilon}, R^{L,\varepsilon})\|\|_* \leq 1 \tag{4.29}$$

*with $\sigma = \gamma + \theta$, $\sigma' = \lfloor \alpha \kappa^{-1} \rfloor + 1 - n$, $a = D$, $a' = c_n \lambda$ and $\eta = 0$, these constants appearing in the definition of the norm (4.28).*

**Remark 4.21.** From the fact that $(\Psi^{L,\varepsilon}, R^{L,\varepsilon})$ solves eq. (4.27), which is an equation depending polynomially on $(\Psi^{L,\varepsilon}, R^{L,\varepsilon})$ and on the coupling constant $\lambda$, we obtain that $\lambda \mapsto (\Psi^{L,\varepsilon}, R^{L,\varepsilon})$ is an analytic function in a neighbourhood of 0, uniformly in $L, \varepsilon$. Indeed, since $(\Psi^{L,\varepsilon}, R^{L,\varepsilon})$ solves a finite-dimensional system which is analytic in all the variables and parameters, the function $\lambda \mapsto (\Psi^{L,\varepsilon}, R^{L,\varepsilon})$ is analytic. Furthermore, defining the power series

$$S^{L,\varepsilon}(\lambda) = \sum_{k=0}^{+\infty} \frac{1}{k!} \|\|(\partial_\lambda^k \Psi^{L,\varepsilon}|_{\lambda=0}, \partial_\lambda^k R^{L,\varepsilon}|_{\lambda=0})\|\|_*,$$

which is convergent for $\lambda$ small enough, using the fact that $F_s^\varepsilon$ is a polynomial in $\Psi$ and $\lambda$ with norm uniformly bounded in $\varepsilon$ and $\lambda$ small enough, we get that $S^{L,\varepsilon}$ satisfies the differential inequality

$$\frac{dS^{L,\varepsilon}}{d\lambda} \leq P(S^{L,\varepsilon}, \lambda),$$

where $P$ is a polynomial with positive coefficients, independent of $\varepsilon$ and $L$, but depending on the degree and the bounds on the coefficients of $F_s^\varepsilon$. Following [ABDG20, Appendix C], we have that $S^{L,\varepsilon}$ has a convergence radius independent of $L, \varepsilon$, which implies that $\lambda \mapsto (\Psi^{L,\varepsilon}, R^{L,\varepsilon})$ has a convergent radius independent of $L, \varepsilon$.

By taking the limit as $\varepsilon \to 0$ and $L \to \infty$, we get that also $\lambda \mapsto (\Psi, R) = \lim_{\varepsilon \to 0, L \to \infty} (\Psi^{L,\varepsilon}, R^{L,\varepsilon})$ is analytic in a neighbourhood of 0. This also implies that the $n$-point functions (which are expectations of the form $\omega(\Psi_\infty(f_1) \cdots \Psi_\infty(f_n))$) are analytic for $\lambda$ in a neighbourhood of 0.

The following lemma will be used in the proof of Proposition 4.20.

**Lemma 4.22.** *Let $N, n, n' \in \mathbb{N}$, let $F = (F_t^{[\ell](k)})_{t \geq 0}^{N \geq \ell > n, k \geq n'}$ be a graded sequence and let $\psi_t, \varphi_t$ be Grassmann fields such that $\|\psi_t\|_{\mathscr{C}_t^{\leq 2}}, \|\varphi_t\|_{\mathscr{C}_t^{\leq 2}} \leq D 2^{\rho t}$ for some $D, \rho > 0$. Then, provided that $\beta > \rho$ and $C < D^{-1}/2$ the following bounds hold true:*

$$\|F_t(\psi_t) - F_t(\varphi_t)\|_{L^\infty} \lesssim_N \|F\|_{C,\alpha,\beta,\kappa} 2^{[\alpha - \rho + n'(\rho - \beta) - n\kappa]t} \|\psi_t - \varphi_t\|_{\mathscr{C}_t^{\leq 2}}$$

$$\|DF_t(\psi_t) - DF_t(\varphi_t)\|_{\mathscr{B}(L^\infty; L^\infty)} \lesssim_N \|F\|_{C,\alpha,\beta,\kappa} 2^{[\alpha - 2\rho + n'(\rho - \beta) - n\kappa]t} \|\psi_t - \varphi_t\|_{\mathscr{C}_t^{\leq 2}}$$

**Proof.** Note the bounds

$$\|F_t^{[\ell](k)}(\psi_t) - F_t^{[\ell](k)}(\varphi_t)\|_{L^\infty}$$
$$\leq k \|F\|_{C,\alpha,\beta,\kappa} C^k D^k 2^{(\alpha - \rho)t} 2^{(\rho - \beta)kt} 2^{-\kappa \ell t} \|\psi_t - \varphi_t\|_{\mathscr{C}_t^{\leq 2}},$$

and

$$\|DF_t^{[\ell](k)}(\psi_t) - DF_t^{[\ell](k)}(\varphi_t)\|_{\mathscr{B}(L^\infty; L^\infty)}$$
$$\leq k(k+1) \|F\|_{C,\alpha,\beta,\kappa} C^{k+1} D^k 2^{(\alpha - 2\rho)t} 2^{(\rho - \beta)(k+1)t} 2^{-\kappa \ell t} \|\psi_t - \varphi_t\|_{\mathscr{C}_t^{\leq 2}}.$$

Then, the summation over $k \geq n'$ is controlled by the condition $CD < 1/2$, whereas the summation over $\ell$ is finite uniformly in $t \geq 0$ but depending on $N$. □

**Proof of Proposition 4.20.** The proof is an application of the Banach fixed-point theorem, since we can interpret (4.27) as a fixed-point equation $(\Psi^{L,\varepsilon}, R^{L,\varepsilon}) = \Gamma[(\Psi^{L,\varepsilon}, R^{L,\varepsilon})]$ where the operator $\Gamma: (\Psi^{L,\varepsilon}, R^{L,\varepsilon}) \mapsto (\tilde{\Psi}^{L,\varepsilon}, \tilde{R}^{L,\varepsilon})$ is defined in an obvious way.

We first prove that $\Gamma$ maps the closed ball $|||(\Psi^{L,\varepsilon}, R^{L,\varepsilon})|||_* \leqslant 1$ into itself. To this end, by using the assumption $\|F^\varepsilon\|_{C,\alpha,\beta,\kappa} \leqslant \lambda$ and that $\|\Psi_s^{L,\varepsilon}\|_{\mathscr{C}_s^{\leqslant 2}} \leqslant D 2^{(\gamma+\theta)s}$ in the said ball, for $C$ small enough Lemma 4.22 (with $\varphi = 0$) gives

$$\|F_s^\varepsilon(\Psi_s^{L,\varepsilon})\|_{L^\infty} \leqslant \tilde{c}_n \lambda 2^{(\alpha-\beta+\gamma+\theta)s}$$
$$\|\mathrm{D}F_s^\varepsilon(\Psi_s^{L,\varepsilon})\|_{\mathscr{B}(L^\infty; L^\infty)} \leqslant \tilde{c}_n \lambda 2^{(\alpha-\beta)s} \qquad (4.30)$$
$$\|H_{>n;s}^\varepsilon(\Psi_s^{L,\varepsilon})\|_{L^\infty} \leqslant \tilde{c}_n \lambda 2^{(\alpha-n\kappa)s},$$

for some constant $\tilde{c}_n$. To control the equation for the remainder, we need to integrate (4.30) in $\mathrm{d}s$, which can be done provided that $n > \alpha \kappa^{-1}$; actually, for $n$ sufficiently large we can make the remainder as small as we want. In fact, we have:

$$\begin{aligned}\|\tilde{R}_t^{L,\varepsilon}\|_{L^\infty} &\leqslant \int_t^\infty \|\omega_t(H_{>n;s}^\varepsilon(\Psi_s^{L,\varepsilon}))\|_{L^\infty} \mathrm{d}s + \int_t^\infty \|\omega_t(\mathrm{D}F_s^\varepsilon(\Psi_s^{L,\varepsilon}) \cdot \dot{G}_s^\varepsilon R_s^{L,\varepsilon})\|_{L^\infty} \mathrm{d}s \\ &\leqslant \int_t^\infty \|H_{>n;s}^\varepsilon(\Psi_s^{L,\varepsilon})\| \mathrm{d}s + \int_t^\infty \|\mathrm{D}F_s^\varepsilon(\Psi_s^{L,\varepsilon})\|_{\mathscr{B}(L^\infty; L^\infty)} \cdot \|\dot{G}_s^\varepsilon\|_{L^1} \|R_s^{L,\varepsilon}\|_{L^\infty} \mathrm{d}s \\ &\lesssim \tilde{c}_n \lambda 2^{-\sigma' t} + \tilde{c}_n \lambda \int_t^\infty 2^{-(d-2\gamma-\alpha+\beta)s} \|R_s^{L,\varepsilon}\|_{L^\infty} \mathrm{d}s,\end{aligned} \qquad (4.31)$$

where we used (4.30) and that $\|\dot{G}_s^\varepsilon\|_{L^1} \lesssim 2^{-(d-2\gamma)s}$ and that $n\kappa - \alpha \geqslant 1$. Since in the ball $|||(\Psi^{L,\varepsilon}, R^{L,\varepsilon})|||_* \leqslant 1$ and since $d - 2\gamma - \alpha + \beta > 0$, we have $\int_t^\infty 2^{-(d-2\gamma-\alpha+\beta)s} \|R_s^{L,\varepsilon}\|_{L^\infty} \mathrm{d}s \leqslant c_n \lambda 2^{-\sigma' t} (d - 2\gamma - \alpha + \beta)^{-1}$, hence, for $c_n$ large enough and for $\lambda$ small enough

$$(c_n \lambda)^{-1} \sup_{t \in \mathbb{R}^+} 2^{\sigma' t} \|\tilde{R}_t^{L,\varepsilon}\|_{L^\infty} \mathrm{d}s \leqslant \frac{1}{2} + \tilde{c}_n \lambda (d - 2\gamma - \alpha + \beta)^{-1} < 1. \qquad (4.32)$$

Next, by Lemma 3.26 we have the estimate $\|X_t^{L,\varepsilon}\|_{\mathscr{C}_t^{\leqslant 2}} \lesssim 2^{\gamma t}$, uniformly in $L, \varepsilon$. Therefore, we write:

$$\begin{aligned}\|\tilde{\Psi}_s^{L,\varepsilon}\|_{\mathscr{C}_t^{\leqslant 2}} &\leqslant \int_0^t \sum_{0 \leqslant n' \leqslant 2} \sum_\nu \|(2^{-n't} \partial_\nu^{n'} \dot{G}_s^\varepsilon)(F_s^\varepsilon(\Psi_s^{L,\varepsilon}) + R_s^{L,\varepsilon})\|_{L^\infty} \mathrm{d}s + \|X_t^{L,\varepsilon}\|_{\mathscr{C}_t^{\leqslant 2}} \\ &\lesssim \int_0^t 2^{-(d-2\gamma)s} (\|F_s^\varepsilon(\Psi_s^{L,\varepsilon})\|_{L^\infty} + \|R_s^{L,\varepsilon}\|_{L^\infty}) \mathrm{d}s + 2^{\gamma t} \\ &\lesssim \tilde{c}_n \lambda \int_0^t 2^{-(d-3\gamma-\alpha+\beta-\theta)s} \mathrm{d}s + 2^{\gamma t} \\ &\lesssim \tilde{c}_n \lambda 2^{(\gamma+\theta)t} + 2^{\gamma t},\end{aligned} \qquad (4.33)$$

where we used (4.30), the propagator estimates in Corollary 3.24, the assumption on $R^{L,\varepsilon}$ and again that $d - 2\gamma - \alpha + \beta > 0$. Therefore, putting together (4.32) and (4.33), for $D$ large enough and $\lambda$ small enough we have that the closed unit ball is left invariant by $\Gamma$, that is $|||(\tilde{\Psi}^{L,\varepsilon}, \tilde{R}^{L,\varepsilon})|||_* \leqslant 1$.

We now prove that $\Gamma$ is a contraction. To this end, assume that $(\tilde{\Psi}^{L,\varepsilon}, \tilde{R}^{L,\varepsilon})$ is another solution in the said ball and write the difference of the two solutions

$$\begin{aligned}\Psi_t^{L,\varepsilon} - \tilde{\Psi}_t^{L,\varepsilon} &= -\int_0^t \dot{G}_s^\varepsilon (F_s^\varepsilon(\Psi_s^{L,\varepsilon}) - F_s^\varepsilon(\tilde{\Psi}_s^{L,\varepsilon}) + R_s^{L,\varepsilon} - \tilde{R}_s^{L,\varepsilon}) \mathrm{d}s, \\ R_t^{L,\varepsilon} - \tilde{R}_t^{L,\varepsilon} &= \int_t^\infty \omega_t(H_{>n;s}^\varepsilon(\Psi_s^{L,\varepsilon}) - H_{>n;s}^\varepsilon(\tilde{\Psi}_s^{L,\varepsilon})) \mathrm{d}s \\ &\quad + \int_t^\infty \omega_t(\mathrm{D}F_s^\varepsilon(\Psi_s^{L,\varepsilon}) \cdot \dot{G}_s^\varepsilon(R_s^{L,\varepsilon} - \tilde{R}_s^{L,\varepsilon})) \mathrm{d}s \\ &\quad + \int_t^\infty \omega_t((\mathrm{D}F_s^\varepsilon(\Psi_s^{L,\varepsilon}) - \mathrm{D}F_s^\varepsilon(\tilde{\Psi}_s^{L,\varepsilon})) \cdot \dot{G}_s^\varepsilon \tilde{R}_s^{L,\varepsilon}) \mathrm{d}s.\end{aligned} \qquad (4.34)$$

By using the assumption $\|F^\varepsilon\|_{C,\alpha,\beta,\kappa} \leqslant \lambda$ as well as $\|\Psi_s^{L,\varepsilon}\|_{\mathscr{C}_s^{\leqslant 2}}, \|\tilde{\Psi}_s^{L,\varepsilon}\|_{\mathscr{C}_s^{\leqslant 2}} \leqslant D 2^{(\gamma+\theta)s}$, Lemma 4.22 implies

$$\begin{aligned}
\|F_s^\varepsilon(\Psi_s^{L,\varepsilon}) - F_s^\varepsilon(\tilde{\Psi}_s^{L,\varepsilon})\|_{L^\infty} &\leqslant \tilde{c}_n \lambda 2^{(\alpha-\beta)s} \|\Psi_s^{L,\varepsilon} - \tilde{\Psi}_s^{L,\varepsilon}\|_{\mathscr{C}_s^{\leqslant 2}}, \\
\|\mathrm{D}F_s^\varepsilon(\Psi_s^{L,\varepsilon}) - \mathrm{D}F_s^\varepsilon(\tilde{\Psi}_s^{L,\varepsilon})\|_{\mathscr{B}(L^\infty;L^\infty)} &\leqslant \tilde{c}_n \lambda 2^{(\alpha-\beta-\gamma-\theta)s} \|\Psi_s^{L,\varepsilon} - \tilde{\Psi}_s^{L,\varepsilon}\|_{\mathscr{C}_s^{\leqslant 2}}, \\
\|H_{>n;s}^\varepsilon(\Psi_s^{L,\varepsilon}) - H_{>n;s}^\varepsilon(\tilde{\Psi}_s^{L,\varepsilon})\|_{L^\infty} &\leqslant \tilde{c}_n \lambda 2^{(\alpha-\gamma-\theta-n\kappa)s} \|\Psi_s^{L,\varepsilon} - \tilde{\Psi}_s^{L,\varepsilon}\|_{\mathscr{C}_s^{\leqslant 2}}.
\end{aligned} \quad (4.35)$$

Therefore, plugging these estimates together with (4.30) into (4.34) we obtain:

$$\begin{aligned}
\|\Psi_t^{L,\varepsilon} &- \tilde{\Psi}_t^{L,\varepsilon}\|_{\mathscr{C}_t^{\leqslant 2}} \\
&\lesssim \tilde{c}_n \lambda 2^{(\gamma+\theta)t} (d - 2\gamma - \alpha + \beta)^{-1} \sup_{s \in \mathbb{R}^+} 2^{-(\gamma+\theta)s} \|\Psi_s^{L,\varepsilon} - \tilde{\Psi}_s^{L,\varepsilon}\|_{\mathscr{C}_s^{\leqslant 2}} \\
&\quad + (d - 2\gamma + \sigma')^{-1} \sup_{s \in \mathbb{R}^+} 2^{\sigma's} \|R_s^{L,\varepsilon} - \tilde{R}_s^{L,\varepsilon}\|_{L^\infty},
\end{aligned} \quad (4.36)$$

and

$$\begin{aligned}
\|R_t^{L,\varepsilon} &- \tilde{R}_t^{L,\varepsilon}\|_{L^\infty} \\
&\lesssim \tilde{c}_n \lambda 2^{-\sigma't} \sup_{s \in \mathbb{R}^+} 2^{-(\gamma+\theta)s} \|\Psi_s^{L,\varepsilon} - \tilde{\Psi}_s^{L,\varepsilon}\|_{\mathscr{C}_s^{\leqslant 2}} \\
&\quad + \tilde{c}_n \lambda 2^{-\sigma't}(d - 2\gamma - \alpha + \beta + \sigma')^{-1} \sup_{s \in \mathbb{R}^+} 2^{\sigma's} \|R_s^{L,\varepsilon} - \tilde{R}_s^{L,\varepsilon}\|_{L^\infty} \\
&\quad + \tilde{c}_n c_n \lambda^2 (d - 2\gamma - \alpha + \beta + \sigma') 2^{-\sigma's} \sup_{s \in \mathbb{R}^+} 2^{-(\gamma+\theta)s} \|\Psi_s^{L,\varepsilon} - \tilde{\Psi}_s^{L,\varepsilon}\|_{\mathscr{C}_s^{\leqslant 2}},
\end{aligned} \quad (4.37)$$

which for $c_n$ large enough and $\lambda$ small enough, imply the claim. $\square$

**Proof of Lemma 3.22.** The previous proof implies that Picard's iteration converges to the solution $\Psi_t^{L,\varepsilon}$. Since $X_t^{L,\varepsilon}(x + m L) = X_t^{L,\varepsilon}(x)$, see Definition 3.20, one can check that every $n$-th Picard's iterate $\Psi_t^{(n)L,\varepsilon}$ is periodic. In fact, if $\Psi_t^{(n)L,\varepsilon}$ is periodic, so are $F^\varepsilon(\Psi_t^{(n)L,\varepsilon})$ and $\dot{\mathfrak{g}}_t^\varepsilon F^\varepsilon(\Psi_t^{(n)L,\varepsilon})$ because of translation invariance, implying that also $\Psi_t^{(n+1)L,\varepsilon}$ is periodic. $\square$

As a consequence of Proposition 4.20, we have that the sequence $(\Psi^{L,\varepsilon})_{L,\varepsilon}$ is bounded in the topology $\|\cdot\|_{\gamma,0}$, see (1.9), introduced before Theorem 1.1. We now address the convergence of the sequence $(\Psi^{L,\varepsilon}, R^{L,\varepsilon})$ as $L \to \infty$ and $\varepsilon \to 0$. We first prove the convergence to the continuum limit, $(\Psi^{L,\varepsilon}, R^{L,\varepsilon}) \to (\Psi^L, R^L)$.

**Proposition 4.23.** *(Continuum limit). Let $L \in \mathbb{N}_\infty$ and $\varepsilon > \varepsilon' \in 2^{-\mathbb{N}_\infty}$. Under the same assumptions of Proposition 4.20, but with $\alpha > 3\beta$, $0 < \theta \leqslant \min\{\alpha - 3\beta; \beta - \gamma\}/2$, there exist $D_\theta$, $c_n$ and $\lambda = \lambda(d, \alpha, \beta, \gamma, n)$ such that, if $\|F^\varepsilon\|_{C,\alpha-\theta,\beta,\kappa}, \|F^{\varepsilon'}\|_{C,\alpha-\theta,\beta,\kappa}, \|H_{>n}^\varepsilon\|_{C,\alpha-\theta,\beta,\kappa}, \|H_{>n}^{\varepsilon'}\|_{C,\alpha-\theta,\beta,\kappa} \leqslant \lambda$ for $C < D_\theta^{-1}/2$, then uniformly in $L$*

$$\|\!|(\Psi^{L,\varepsilon} - \Psi^{L,\varepsilon'}, R^{L,\varepsilon} - R^{L,\varepsilon'})|\!\|_* \lesssim \varepsilon^\theta,$$

*with $\sigma = \gamma + \theta$, $\sigma' = \lfloor \alpha \kappa^{-1} \rfloor + 1 - n$, $a = D_\theta$, $a' = c_n \lambda$ and $\eta = 0$, these constants appearing in the definition of the norm (4.28).*

**Proof.** The equation for the difference is similar to eq. (4.34), but involves the differences $\delta_{\varepsilon,\varepsilon'} F_s := F_s^\varepsilon - F_s^{\varepsilon'}$, $\delta_{\varepsilon,\varepsilon'} \dot{G}_s := \dot{G}_s^\varepsilon - \dot{G}_s^{\varepsilon'}$ and $\delta_{\varepsilon,\varepsilon'} X_t^L := X_t^{L,\varepsilon} - X_t^{L,\varepsilon'}$ in an obvious way. We control such differences by Corollary 3.24, Corollary 3.26, and Proposition 4.19:

$$\begin{aligned}
\|\delta_{\varepsilon,\varepsilon'} \dot{G}_t\|_{L_t^p} &\lesssim 2^{-(d/p - 2\gamma - \theta)t} \varepsilon^\theta, \\
\|\delta_{\varepsilon,\varepsilon'} F\|_{C,\alpha,\beta,\kappa} &\lesssim \tilde{c}_n \varepsilon^\theta, \\
\|\delta_{\varepsilon,\varepsilon'} X_t^L\|_{\mathscr{C}_t^{\leqslant 2}} &\lesssim 2^{(\gamma+\theta)t} \varepsilon^\theta,
\end{aligned} \quad (4.38)$$

for some constant $\tilde{c}_n$. Therefore, we can follow the analysis after (4.34) to obtain the same estimates but with the extra gain $\varepsilon^\theta$ and thus prove the claim. $\square$

To prove convergence to the infinite-volume limit $(\Psi^L, R^L) \to (\Psi, R)$, we need to take $\eta > 0$ in $\|\!|\cdot|\!\|_*$.

**Proposition 4.24.** *(Infinite-volume limit). Let $L < L' \in \mathbb{N}_\infty$ and $\varepsilon \in 2^{-\mathbb{N}_\infty}$. Under the same assumptions of Proposition 4.20 and for any $\eta > 0$, there exist constants $D_\eta$, $c_n$ and $\lambda = \lambda(d, \alpha, \beta, \gamma, n)$ such that, if $\|F^\varepsilon\|_{C,\alpha,\beta,\kappa}, \|H^\varepsilon_{>n}\|_{C,\alpha,\beta,\kappa} \leqslant \lambda$ for $C < D_\eta^{-1}/2$, then uniformly in $\varepsilon$*

$$\|\!|(\Psi^{L',\varepsilon} - \Psi^{L,\varepsilon}, R^{L',\varepsilon} - R^{L,\varepsilon})|\!\|_* \lesssim L^{-\eta},$$

*with $\sigma = \gamma + \theta$, $\sigma' = \lfloor \alpha \kappa^{-1} \rfloor + 1 - n$, $a = D_\eta$, $a' = c_n \lambda$, these constants appearing in the definition of the norm (4.28).*

We skip the proof of this proposition because the analysis of the next section carries out a similar argument in a more complex setting. See also [ABDG20] for similar arguments in the context of the infinite-volume limit of Grassmann SPDEs.

So far we have established the existence of the field $\Psi^{L,\varepsilon}$ and its convergence in a suitable topology for distributional fields on $\mathbb{R}^+ \times \mathbb{R}^d$ as the cutoffs are removed, $L \to \infty$ and $\varepsilon \to 0$. We shall now consider the properties of the sequence $(\Psi^{L,\varepsilon}_t)$ for fixed $t$ in suitable Besov spaces.

**Corollary 4.25.** *With the same setting of Proposition 4.20, for any $L' \geqslant L \in \mathbb{N}_\infty$, $\varepsilon \geqslant \varepsilon' \in 2^{-\mathbb{N}_\infty}$, $\rho > 0$ and $t \in [0, \infty]$*

$$\|\Psi^{L,\varepsilon}_t\|_{B^{-\gamma-\rho}_{\infty,\infty}} \lesssim \rho^{-1/2}, \qquad \|\Psi^{L,\varepsilon}_t - \Psi^{L',\varepsilon'}_t\|_{B^{-\gamma-\theta-\rho}_{\infty,\infty,\eta}} \lesssim \rho^{-1/2}(\varepsilon^\theta \mathbf{1}_{\varepsilon>\varepsilon'} + L^{-\eta} \mathbf{1}_{L'>L}). \qquad (4.39)$$

**Proof.** Denote by $(\Delta_j)_{j \geqslant -1}$ the Littlewood–Paley blocks. It suffices to bound $\|\Psi^{L,\varepsilon}_t(\Delta_j)\|$ in the fixed-point equation (4.27) by noting that the function $\dot{G}^\varepsilon_s \Delta_j$ satisfies the estimate

$$\|\dot{G}^\varepsilon_s \Delta_j(\cdot - x)\|_{L^1} \lesssim 2^{-(d-\gamma)s} 2^{\gamma j} \|\Delta_j\|_{L^1} = 2^{-(d-\gamma)s} 2^{\gamma j},$$

and by recalling that $\|X^{L,\varepsilon}_t(\Delta_j(\cdot - x))\| \lesssim 2^{(\gamma+\rho)j}$, for any $\rho > 0$, see Corollary 3.26. Accordingly, for any $L \in \mathbb{N}_\infty$, $\varepsilon \in 2^{-\mathbb{N}_\infty}$ and $t \geqslant 0$ we have

$$\begin{aligned}
\|\Psi^{L,\varepsilon}_t(\Delta_j(\cdot - x))\| &\lesssim \int_0^t \|\dot{G}^\varepsilon_s \Delta_j(\cdot - x)\|_{L^1} \|(F^\varepsilon_s(\Psi^{L,\varepsilon}_s) + R^{L,\varepsilon}_s)\|_{L^\infty} ds + \|X^{L,\varepsilon}_t(\Delta_j(\cdot - x))\| \\
&\lesssim 2^{\gamma j} \int_0^t 2^{-(d-\gamma)s} \|(F^\varepsilon_s(\Psi^{L,\varepsilon}_s) + R^{L,\varepsilon}_s)\|_{L^\infty} ds + 2^{(\gamma+\rho)j} \\
&\lesssim_\rho 2^{(\gamma+\rho)j},
\end{aligned} \qquad (4.40)$$

where we used that $d - 2\gamma - \alpha + \beta > 0$. To prove the other bound, we first consider the equation for the difference

$$\begin{aligned}
&\Psi^{L,\varepsilon}_t(\Delta_j(\cdot - x)) - \Psi^{L,\varepsilon'}_t(\Delta_j(\cdot - x)) \\
&= \int_0^t \left[ \Delta_j(\cdot - x) \dot{G}^\varepsilon_s (F^\varepsilon_s(\Psi^{L,\varepsilon}_s) + R^{L,\varepsilon}_s) - \Delta_j(\cdot - x) \dot{G}^{\varepsilon'}_s \left(F^{\varepsilon'}_s\left(\Psi^{L,\varepsilon'}_s\right) + R^{L,\varepsilon'}_s\right) \right] ds \\
&\quad + X^{L,\varepsilon}_t(\Delta_j(\cdot - x)) - X^{L,\varepsilon'}_t(\Delta_j(\cdot - x)).
\end{aligned}$$

and write the integrand in terms of the differences $\delta_{\varepsilon,\varepsilon'} F_s$, $\delta_{\varepsilon,\varepsilon'} \dot{G}_s$, $\Psi^{L,\varepsilon} - \Psi^{L,\varepsilon'}$ and $R^{L,\varepsilon} - R^{L,\varepsilon'}$. The estimates for $\delta_{\varepsilon,\varepsilon'} F_s$ and $\delta_{\varepsilon,\varepsilon'} \dot{G}_s$ are provided in (4.38), the estimates for $\|\Psi^{L,\varepsilon} - \Psi^{L,\varepsilon'}\|_{\mathscr{C}_s^{\leq 2}}$ and for $\|R^{L,\varepsilon} - R^{L,\varepsilon'}\|_{L^\infty}$ in Proposition 4.23 whereas the estimate $\|X_t^{L,\varepsilon}(\Delta_j(\cdot-x)) - X_t^{L,\varepsilon'}(\Delta_j(\cdot-x))\| \lesssim_\rho \varepsilon^\theta 2^{(\gamma+\theta+\rho)j}$ is provided in Corollary 3.26. Carrying out the steps as in (4.40) we gain the extra $\varepsilon^\theta$ in the drift and thus prove that

$$\|\Psi_t^{L,\varepsilon} - \Psi_t^{L,\varepsilon'}\|_{B_{\infty,\infty}^{-\gamma-\theta-\rho}} \lesssim_\rho \varepsilon^\theta. \tag{4.41}$$

We then consider the equation for the other difference $\Psi_t^{L',\varepsilon}(\Delta_j(\cdot-x)) - \Psi_t^{L,\varepsilon}(\Delta_j(\cdot-x))$ and expand it in the differences $\Psi^{L',\varepsilon} - \Psi^{L,\varepsilon}$ and $R^{L',\varepsilon} - R^{L,\varepsilon}$ which are estimated by Proposition 4.24 in the $\|\cdot\|_{\mathscr{C}_t^{\leq 2}(\eta)}$ and $\|\cdot\|_{L^\infty(\eta)}$ norms respectively. Since $\|X_t^{L',\varepsilon}(\Delta_j(\cdot-x)) - X_t^{L,\varepsilon}(\Delta_j(\cdot-x))\| \lesssim_\rho L^{-\eta} 2^{(\gamma+\rho)j}$, see Corollary 3.26, following the same strategy used above, we have

$$\|\Psi_t^{L,\varepsilon} - \Psi_t^{L,\varepsilon'}\|_{B_{\infty,\infty,\eta}^{-\gamma-\rho}} \lesssim_\rho L^{-\eta},$$

which, together with (4.41) implies (4.39). □

**Proof Theorem 1.1.** The global existence of $\Psi^{L,\varepsilon}$ for $L \in \mathbb{N}$ and $\varepsilon \in 2^{-\mathbb{N}}$ follows by Proposition 4.20. In fact, having fixed $\gamma$, choose some $\alpha, \beta, \kappa$ such the constraints in (4.14), (4.15), (4.16) are satisfied. Choose $n = \lfloor \alpha \kappa^{-1} \rfloor + 2$ (so that $\gamma < \sigma'$) and $C \leqslant C_0(\alpha, \beta, \gamma, \kappa, n) \wedge D^{-1}/2$ and (by abuse of notation) $\lambda_0(d,\gamma) := 1 \wedge \lambda_0(C) \wedge \lambda_0(d,\alpha,\beta,\gamma)$, where $C_0$ and $\lambda_0(C)$ are as in Corollary 4.15, whereas $D$ and $\lambda_0(d,\alpha,\beta,\gamma)$ from Proposition 4.20. By Theorem 4.13 there exists a unique global solution of the $n$-truncated equation with boundary datum $F_0^{\varepsilon[0]}(\psi_t) := \lambda(\psi_t)^3$ for $\lambda \leqslant \lambda_0(d,\gamma)$. Denote by $\mu_\infty^\varepsilon(\lambda) \equiv \mu_\infty^\varepsilon(F_0^{\varepsilon[0]})$ the function defined in Definition 4.17 and note that $F_\infty^\varepsilon = DV_\infty^\varepsilon$ provided that $\mu^\varepsilon(\lambda) = \mu_\infty^\varepsilon(\lambda)$. By Corollary 4.15, $F^\varepsilon$ and $H_{>n}^\varepsilon$ satisfy the assumptions of Proposition 4.20, which for $\theta = 0$ and for any $L \in \mathbb{N}$ and $\varepsilon \in 2^{-\mathbb{N}}$ implies the existence of a unique pair $(\Psi^{L,\varepsilon}, R^{L,\varepsilon})$ solving (4.27) and by Lemma 3.14 also the FBSDE which was the goal.

By combining Proposition 4.23 with Proposition 4.24 we have the existence of the limit $\Psi = \lim_{L\to\infty} \lim_{\varepsilon\to 0} \Psi^{L,\varepsilon}$ and the bound (1.10). Finally, (1.11) is a direct consequence of Corollary 4.25. □

**Remark 4.26.** Note that the existence of the weak limit $\omega^V$ stated in Theorem 1.3 follows by Theorem 1.1, by duality $\|(\Psi_\infty - \Psi_\infty^{L,\varepsilon})(f)\| \lesssim \|f\|_{B_{1,1,-\eta}^{\gamma+\theta}} \|\Psi_\infty - \Psi_\infty^{L,\varepsilon}\|_{B_{\infty,\infty,\eta}^{-\gamma-\theta}}$ and by the stochastic quantisation formula of Theorem 3.12.

## 4.4 Exponential decay and short-distance divergence of correlations

In this section, we prove the exponential decay of correlation functions via a coupling method and the short-distance divergence of the two-point correlation function by using a generating function. We learned the coupling method from the work [Fun91] and more recently explored in [GHR22] in the context of stochastic quantisation. For the sake of simplicity, we shall restrict ourselves to the two-cluster correlation function. The basic idea is to approximate the solution of the FBSDE in suitably disjoint spatial domains by two independent Grassmann processes that solve some other FBSDEs with uncorrelated GBMs as sources. Because the latter are related to a massive operator, one expects that the approximation can be made "exponentially accurate", depending on the distance of the regions. This technique can be viewed as a sensible replacement of the well-known cluster expansions [Bry84, GJ87] of statistical mechanics, applied in the context of stochastic quantisation in [Dim90].

Let us discuss the coupling method in more detail. Consider $((f^{(i,k)})_{k=1}^{m_i})_{i=1}^2 \subset B_{1,1}^{\gamma^+}(\mathbb{R}^d; \mathbb{C}^4)$ compactly supported, define $D_i := \bigcup_{k=1}^{m_i} \mathrm{supp}(f^{(i,k)})$, and assume that $D_1 \cap D_2 = \emptyset$. We let $\mathcal{D}_i \supset D_i$ be such that $\mathcal{D}_1 \cap \mathcal{D}_2 = \emptyset$ and such that $\sup_i \mathrm{dist}(\mathcal{D}_i^c, D_i)$ is as large as possible, where $\mathcal{D}_i^c := \mathbb{R}^d \setminus \mathcal{D}_i$. For brevity, we will here consider the case $L = \infty$ and $\varepsilon = 0$ directly. Recall that $\dot{G}_t = \mathfrak{C}_t^2 U$, with $[\mathfrak{C}_t, U] = [\mathfrak{C}_t, \Theta] = 0$.

**Definition 4.27.** *We let $(X_t^{(i)})_t^{i=0,1,2}$ be the family of anticommuting norm-compatible GBMs such that, setting $\mathcal{D}_0 := \mathbb{R}^d$,*

$$\omega(X_t^{(i)}(f) X_t^{(j)}(g)) = \left\langle\!\!\left\langle f, \int_0^t \mathfrak{C}_t \mathbf{1}_{\mathcal{D}_i} \mathbf{1}_{\mathcal{D}_j} \mathfrak{C}_t U g \mathrm{d}s \right\rangle\!\!\right\rangle_\Theta \qquad \forall f, g \in \mathfrak{h}. \tag{4.42}$$

As anticipated, the process $(X_t^{(0)})_t := (X_t)_t$ is the diffusion of the FBSDEs considered so far, whereas $(X_t^{(i)})_t$, $i = 1, 2$, are independent local approximations of the latter in the region $D_i$. The independence is read out from (4.42), which, by the choice of $\mathcal{D}_i$ implies $\omega(X_t^{(1)}(f) X_t^{(2)}(g)) = 0$, for any $f, g \in \mathfrak{h}$.

If $\psi$ is a field on $\mathbb{R}^d$, that is, $\mathbb{R}^d \ni x \mapsto (\psi_\mu(x))_\mu$, we introduce the weighted topologies

$$\|\psi\|_{L_t^\infty(\mathcal{D};B)} := \sup_{x \in \mathbb{R}^d} e^{-\xi 2^t \mathrm{dist}(x,B)} e^{\xi \mathrm{dist}(x,\mathcal{D})} \|\psi_\mu(x)\|, \tag{4.43}$$

$$\|\psi\|_{\mathscr{C}_t^n(\mathcal{D};B)} := 2^{-nt} \sup_\nu \|\partial_\nu^n \psi\|_{L_t^\infty(\mathcal{D};B)}, \tag{4.44}$$

where $t \geq 0$, where $\xi > 0$ is small enough, and where $B, \mathcal{D} \subset \mathbb{R}^d$ are measurable domains. We also set $\|\cdot\|_{\mathscr{C}_t^{\leq n}(\mathcal{D};B)} := \sum_{0 \leq n' \leq n} \|\cdot\|_{\mathscr{C}_t^{n'}(\mathcal{D};B)}$. If $B = \mathbb{R}^d$, then we abridge our notation to $\|\cdot\|_{\mathscr{C}_t^n(\mathcal{D})}$, $\|\cdot\|_{\mathscr{C}_t^{\leq n}(\mathcal{D})}$ and $\|\cdot\|_{L^\infty(\mathcal{D})}$, the right-hand side of the latter being independent of $t$.

**Lemma 4.28.** *Let $(X_t^{(i)})_t^i$, $(D_i)_i$ and $(\mathcal{D}_i)_i$ be as above. Then, for $\xi > 0$ small enough, uniformly in $t \geq 0$*

$$\|X_t - X_t^{(i)}\|_{\mathscr{C}_t^{\leq n}(\mathcal{D}_i^c; D_i)} \lesssim \mathrm{dist}(D_i, \mathcal{D}_i^c)^{-\gamma} e^{-\xi \mathrm{dist}(D_i, \mathcal{D}_i^c)}. \tag{4.45}$$

**Proof.** We note that the derivative in $t$ of the covariance of $X_t - X_t^{(i)}$ is $\mathfrak{C}_t \mathbf{1}_{\mathcal{D}_i^c} \mathfrak{C}_t U = (\mathbf{1}_{\mathcal{D}_i^c} \mathfrak{C}_t)^* \mathbf{1}_{\mathcal{D}_i^c} \mathfrak{C}_t U$. Accordingly, by the norm compatibility we can write

$$\|(\partial_\nu^n X_t - \partial_\nu^n X_t^{(i)})(x)\|^2 \lesssim \int_{\mathbb{R}^d} \int_0^t |\partial_\nu^n \mathfrak{C}_s(z;x)|^2 \mathbf{1}_{\mathcal{D}_i^c}(z) \mathrm{d}z \, \mathrm{d}s,$$

where $\lesssim$ is up to some universal constant. Therefore, by using the estimate $|\partial_x^\gamma \mathfrak{C}_s(z;x)| \lesssim 2^{(\gamma+d/2+n)s} e^{-\bar{c} 2^s |z-x|}$, see Remark 3.25, we have

$$\begin{aligned}
&2^{-nt} e^{-\xi 2^t \mathrm{dist}(x, D_i)} e^{\xi \mathrm{dist}(x, \mathcal{D}_i^c)} \|(\partial_\nu^n X_t - \partial_\nu^n X_t^{(i)})(x)\| \\
&\lesssim 2^{-nt} e^{-\xi 2^t \mathrm{dist}(x, D_i)} e^{\xi \mathrm{dist}(x, \mathcal{D}_i^c)} \left(\int_{\mathbb{R}^d} \int_0^t |\partial_\nu^n \mathfrak{C}_s(z;x)|^2 \mathbf{1}_{\mathcal{D}_i^c}(z) \mathrm{d}s \mathrm{d}z\right)^{1/2} \\
&\lesssim 2^{-nt} \left[\int_{\mathbb{R}^d} \int_0^t 2^{(2\gamma+d+2n)s} e^{-\xi 2^s \mathrm{dist}(x, D_i)} e^{-(\bar{c} 2^s - \xi)|z-x|} \mathbf{1}_{\mathcal{D}_i^c}(z) \mathrm{d}s \mathrm{d}z\right]^{1/2} \\
&\lesssim \left[\int_0^t 2^{(2\gamma+d)s} e^{-\xi 2^s \mathrm{dist}(x, D_i)} e^{-\bar{c} 2^s \mathrm{dist}(x, \mathcal{D}_i^c)/3} \int_{\mathbb{R}^d} e^{-\bar{c} 2^s |z-x|/3} \mathrm{d}z \, \mathrm{d}s\right]^{1/2} \\
&\lesssim \left[\int_0^t 2^{2\gamma s} e^{-\xi 2^s \mathrm{dist}(D_i, \mathcal{D}_i^c)} \mathrm{d}s\right]^{1/2}, \\
&\lesssim_\xi \mathrm{dist}(D_i, \mathcal{D}_i^c)^{-\gamma} e^{-\xi \mathrm{dist}(D_i, \mathcal{D}_i^c)},
\end{aligned} \tag{4.46}$$

provided that $\xi \leqslant \bar{c}/3$, hence the claim. □

Let $(\Psi^{(i)})_{i=0,1,2}$ be the Grassmann processes that solve the FBSDEs with $(X^{(i)})_{i=0,1,2}$ as sources. More precisely, one should think of the pairs $(\Psi_t^{(i)}, R_t^{(i)})$, see Lemma 3.22, which solve a system of equations that are always well-defined, unlike the FBSDE.

Clearly, also $(\Psi^{(i)})_{i=1,2}$ are independent and one should think of $\Psi^{(i)}$ as approximating $\Psi := \Psi^{(0)}$ in the regions $D_i$ respectively. This allows for the representation.

**Lemma 4.29.** *Let $F, F' \in \bigoplus_{j=0}^n \left(B_{1,1,0\text{-}}^{\gamma^+}(\mathbb{R}^d; \mathbb{C}^4)\right)^{\wedge j}$, for some $n \in \mathbb{N}$, and let $X$, $\Psi$, $\Psi^{(i)}$ be as above. Then*

$$\text{Cov}_{\omega^\nu}(F(X_\infty); F'(X_\infty)) = \text{Cov}_\omega(F(\Psi_\infty) - F(\Psi_\infty^{(1)}); F'(\Psi_\infty)) + \text{Cov}_\omega(F(\Psi_\infty^{(1)}); F'(\Psi_\infty) - F'(\Psi_\infty^{(2)}))$$

**Proof.** First of all, note that $F(X_\infty), F'(X_\infty) \in \mathcal{A}$, where $\mathcal{A}$ was defined in (1.12), since in fact $X_\infty \equiv \psi$. By the stochastic quantisation formula in Theorem 3.12 and by Corollary 4.25

$$\begin{aligned}\text{Cov}_{\omega^\nu}(F(X_\infty); F'(X_\infty)) &= \lim_{L \to \infty} \lim_{\varepsilon \to 0} \text{Cov}_{\omega^\nu}(F(X_\infty^{L,\varepsilon}); F'(X_\infty^{L,\varepsilon})) \\ &= \lim_{L \to \infty} \lim_{\varepsilon \to 0} \text{Cov}_\omega(F(\Psi_\infty^{L,\varepsilon}); F'(\Psi_\infty^{L,\varepsilon})) \\ &= \text{Cov}_\omega(F(\Psi_\infty); F'(\Psi_\infty)).\end{aligned}$$

Note in particular that $F(X_\infty^{L,\varepsilon})$ and $F'(X_\infty^{L,\varepsilon})$ are elements of a finite-dimensional Grassmann algebra. Since $\Psi^{(1)}$ and $\Psi^{(2)}$ are independent, $\text{Cov}_\omega(F(\Psi_\infty^{(1)}); F'(\Psi_\infty^{(2)})) = 0$ and the claim follows by simple manipulations. □

This lemma reduces the problem to controlling, e.g., $|\text{Cov}_\omega(F(\Psi_\infty) - F(\Psi_\infty^{(1)}); F'(\Psi_\infty))|$. We also need another technical lemma, which is a generalisation of Lemma 4.22. Note that because we work without cut-offs, in particular with $\varepsilon = 0$, we work with exponential weights, compare with (4.44). In particular, the spaces $\mathcal{L}_t^{(k)}(s;c)$ should be intended with the exponential weights $w_{t;c}^{(k)}$ with $\delta = 1$. For simplicity, we denote by $\mathcal{B}_\xi(L^\infty; L^\infty)$ the weighted space $\mathcal{L}_0^{(1,0)}(0;\xi)$.

**Lemma 4.30.** *In the same setting of Lemma 4.22 the following bounds hold true, for $\xi$ small enough*

$$\begin{aligned}\|DF_t(\psi_t)\|_{\mathcal{B}_\xi(L^\infty;L^\infty)} &\lesssim_N \|F\|_{C,\alpha,\beta,\kappa} 2^{[\alpha-\rho+n'(\rho-\beta)-n\kappa]t} \\ \|F_t(\psi_t) - F_t(\varphi_t)\|_{L^\infty(\mathcal{D};B)} &\lesssim_N \|F\|_{C,\alpha,\beta,\kappa} 2^{[\alpha-\rho+n'(\rho-\beta)-n\kappa]t} \|\psi_t - \varphi_t\|_{\mathscr{C}_t^{\leqslant 2}(\mathcal{D};B)} \\ \|DF_t(\psi_t) - DF_t(\varphi_t)\|_{\mathcal{B}(L^\infty;L^\infty(\mathcal{D};B))} &\lesssim_N \|F\|_{C,\alpha,\beta,\kappa} 2^{[\alpha-2\rho+n'(\rho-\beta)-n\kappa]t} \|\psi_t - \varphi_t\|_{\mathscr{C}_t^{\leqslant 2}(\mathcal{D};B)}.\end{aligned}$$

**Proof.** The first bound follows by the same argument of Lemma 4.22. Regarding the other two bounds, we have to suitably distribute the weight. For example, we have

$$\begin{aligned}&\|F_t^{[\ell](k)}(\psi_t) - F_t^{[\ell](k)}(\varphi_t)\|_{L^\infty(\mathcal{D};B)} \\ &\leqslant k\|F\|_{C,\alpha,\beta,\kappa} C^k D^k 2^{(\alpha-\rho)t} 2^{(\rho-\beta)kt} 2^{-\kappa \ell t} \|\psi_t - \varphi_t\|_{\mathscr{C}_t^{\leqslant 2}(\mathcal{D};B)},\end{aligned}$$

compare with the proof of Lemma 4.22, where we used the triangle inequality to distribute the weight $e^{-\xi 2^t \text{dist}(y,B)} e^{\xi \text{dist}(y,\mathcal{D})}$, for $\xi$ small enough, on the Steiner tree for any $j=1,\ldots,k$:

$$e^{-\xi 2^t \text{dist}(y,B)} e^{\xi \text{dist}(y,\mathcal{D})} |F_t^{[\ell](k)}((y,\mu);(x_1,\mu_1)\ldots;(x_k,\mu_k))|$$
$$\leqslant |F_t^{[\ell](k)}((y,\mu);(x_1,\mu_1)\ldots;(x_k,\mu_k))| w_{t;c}^{(k+1)}(y;x_1;\ldots;x_k) e^{-\xi 2^t \text{dist}(x_j,B)} e^{\xi \text{dist}(x_j,\mathcal{D})}. \qquad \square$$

**Proof of Part 2 in Theorem 1.3.** By Corollary 4.25 and by duality we have

$$\|\Psi_\infty(f^{(j,k)})\| \lesssim \|f^{(j,k)}\|_{B_{1,1,0^-}^{\gamma^+}}, \qquad \|\Psi_\infty^{(i)}(f^{(j,k)})\| \lesssim \|f^{(j,k)}\|_{B_{1,1,0^-}^{\gamma^+}}.$$

Then, by Lemma 4.29 we can write

$$\left|\text{Cov}_{\omega^V}\left(\prod_{k=1}^{m_1} X_\infty(f^{(1,k)});\prod_{k=1}^{m_2} X_\infty(f^{(2,k)})\right)\right|$$
$$\lesssim \sum_{i=1,2,j=1,2,j\neq i} \prod_{k'=1}^{m_j} \|f^{(j,k')}\|_{B_{1,1,0^-}^{\gamma^+}} \left(\sum_k^{m_i} \|(\Psi_\infty - \Psi_\infty^{(i)})(f^{(i,k)})\| \prod_{\ell \neq k} \|f^{(i,\ell)}\|_{B_{1,1,0}^{\gamma^+}}\right) \quad (4.47)$$

If $\text{dist}(D_1,D_2) \leqslant 1$, we write $\|(\Psi_\infty - \Psi_\infty^{(i)})(f^{(i,k)})\| \leqslant \|\Psi_\infty(f^{(i,k)})\| + \|\Psi_\infty^{(i)}(f^{(i,k)})\| \lesssim \|f^{(i,k)}\|_{B_{1,1,0^-}^{\gamma^+}}$. Otherwise, we will prove that $\|(\Psi_\infty - \Psi_\infty^{(i)})(f^{(i,k)})\|$ decays exponentially in $\text{dist}(\mathcal{D}_i^c, D_i)$, and thus in $\text{dist}(D_1,D_2)$, since we can choose $\mathcal{D}_i$ so that $\min_i \text{dist}(\mathcal{D}_i^c, D_i) \sim \text{dist}(D_1,D_2)$.

First of all, we note that $\|f\|_{L^1} \lesssim \|f\|_{B_{1,1,0^-}^a}$ and hence, since $\text{supp}(f^{(i,k)}) \subseteq D_i$

$$\|e^{-\xi 2^t \text{dist}(\cdot,D_i)} e^{\xi \text{dist}(\cdot,\mathcal{D}_i^c)} f^{(i,k)}\|_{L^1} \lesssim \|f\|_{B_{1,1,0^-}^{\gamma^+}} e^{-\xi \text{dist}(D_i,\mathcal{D}_i^c)}.$$

Therefore, by Hölder's inequality we obtain

$$\|(\Psi_t - \Psi_t^{(i)})(f^{(i,k)})\| \lesssim \|f\|_{B_{1,1,0^-}^{\gamma^+}} e^{-\xi \text{dist}(D_i,\mathcal{D}_i^c)} \|\Psi_t - \Psi_t^{(i)}\|_{L_t^\infty(\mathcal{D}_i^c;D_i)}, \quad (4.48)$$

We first of all control $\|\Psi_t - \Psi_t^{(i)}\|_{\mathscr{C}_t^{\leqslant 2}(\mathcal{D}_i^c)} \equiv \|\Psi_t - \Psi_t^{(i)}\|_{\mathscr{C}_t^{\leqslant 2}(\mathcal{D}_i^c;\mathbb{R}^d)}$ and $\|R_t - R_t^{(i)}\|_{L^\infty(\mathcal{D}_i^c)} \equiv \|R_t - R_t^{(i)}\|_{L^\infty(\mathcal{D}_i^c;\mathbb{R}^d)}$. Note that by Lemma 4.30, the following bounds hold for some constant $\tilde{c}_n$:

$$\|DF_s(\Psi_s)\|_{\mathscr{B}_\xi(L^\infty;L^\infty)} \leqslant \tilde{c}_n \lambda 2^{(\alpha-\beta)s}$$
$$\|F_s(\Psi_s) - F_s(\Psi_s^{(i)})\|_{L^\infty(\mathcal{D}_i^c)} \leqslant \tilde{c}_n \lambda 2^{(\alpha-\beta)s} \|\Psi_s - \Psi_s^{(i)}\|_{\mathscr{C}_s^{\leqslant 2}(\mathcal{D}_i^c)},$$
$$\|DF_s(\Psi_s) - DF_s(\Psi_s^{(i)})\|_{\mathscr{B}(L^\infty;L^\infty(\mathcal{D}_i^c))} \leqslant \tilde{c}_n \lambda 2^{(\alpha-\beta-\gamma)s} \|\Psi_s - \Psi_s^{(i)}\|_{\mathscr{C}_s^{\leqslant 2}(\mathcal{D}_i^c)},$$
$$\|H_{>n;s}(\Psi_s) - H_{>n;s}(\Psi_s^{(i)})\|_{L^\infty(\mathcal{D}_i^c)} \leqslant \tilde{c}_n \lambda 2^{(\alpha-n\kappa-\gamma)s} \|\Psi_s - \Psi_s^{(i)}\|_{\mathscr{C}_s^{\leqslant 2}(\mathcal{D}_i^c)}.$$

We plug these bounds in the equations for $\Psi_t - \Psi_t^{(i)}$ and $R_t - R_t^{(i)}$, see Lemma 3.22 to obtain, for $\xi$ small enough

$$\|\Psi_t - \Psi_t^{(i)}\|_{\mathscr{C}_t^{\leqslant 2}(\mathcal{D}_i^c)}$$
$$\leqslant \int_0^t \|e^{\xi|\cdot|} \dot{G}_s\|_{L^1} \left[\|F_s(\Psi_s) - F_s(\Psi_s^{(i)})\|_{L^\infty(\mathcal{D}_i^c)} + \|R_s - R_s^{(i)}\|_{L^\infty(\mathcal{D}_i^c)}\right] ds$$
$$\quad + \|X_t - X_t^{(i)}\|_{\mathscr{C}_t^{\leqslant 2}(\mathcal{D}_i^c)}$$
$$\lesssim \tilde{c}_n \int_0^t 2^{-(d-2\gamma)s} \left[\lambda 2^{(\alpha-\beta)s} \|\Psi_s - \Psi_s^{(i)}\|_{\mathscr{C}_s^{\leqslant 2}(\mathcal{D}_i^c)} + \|R_s - R_s^{(i)}\|_{L^\infty(\mathcal{D}_i^c)}\right] ds + 2^{\gamma t} \quad (4.49)$$
$$\lesssim \tilde{c}_n \lambda (d - 2\gamma - \alpha + \beta)^{-1} 2^{\gamma t} \sup_{s \in \mathbb{R}^+} 2^{-\gamma s} \|\Psi_s - \Psi_s^{(i)}\|_{\mathscr{C}_s^{\leqslant 2}(\mathcal{D}_i^c)}$$
$$\quad + \int_0^\infty 2^{-(d-2\gamma)s} \|R_s - R_s^{(i)}\|_{L^\infty(\mathcal{D}_i^c)} ds + 2^{\gamma t},$$

and

$$\begin{aligned}
\|R_t - R_t^{(i)}\|_{L^\infty(\mathcal{D}_i^c)} &\lesssim \int_t^\infty \|H_{>n;s}(\Psi_s) - H_{>n;s}(\Psi_s^{(i)})\|_{L^\infty(\mathcal{D}_i^c)} \mathrm{d}s \\
&+ \int_t^\infty \|DF_s(\Psi_s)\|_{\mathcal{B}_\xi(L^\infty;L^\infty)} 2^{-(d-2\gamma)s} \|R_s - R_s^{(i)}\|_{L^\infty(\mathcal{D}_i^c)} \mathrm{d}s \\
&+ \int_t^\infty \|DF_s(\Psi_s) - DF_s(\Psi_s^{(i)})\|_{\mathcal{B}(L^\infty;L^\infty(\mathcal{D}_i^c))} 2^{-(d-2\gamma)s} \|R_s^{(i)}\|_{L^\infty} \mathrm{d}s \\
&\lesssim \tilde{c}_n \lambda [1 + c_n \lambda] \sup_{s \in \mathbb{R}^+} 2^{-\gamma s} \|\Psi_s - \Psi_s^{(i)}\|_{\mathscr{C}_s^{\leq 2}(\mathcal{D}_i^c)} \\
&+ \lambda \int_t^\infty 2^{-(d-2\gamma-\alpha+\beta)s} \|R_s - R_s^{(i)}\|_{L^\infty(\mathcal{D}_i^c)} \mathrm{d}s
\end{aligned}$$ (4.50)

where we used that $\int_t^\infty 2^{-(d-2\gamma-\alpha+\beta)s} \|R_s^{(i)}\|_{L^\infty} \mathrm{d}s \lesssim c_n \lambda$. For $\lambda$ small enough, (4.49) and (4.50) imply $\sup_{t \in \mathbb{R}^+} 2^{-\gamma t} \|\Psi_t - \Psi_t^{(i)}\|_{\mathscr{C}_t^{\leq 2}(\mathcal{D}_i^c)} \lesssim 1$ and $\sup_{t \in \mathbb{R}^+} \|R_t - R_t^{(i)}\|_{L^\infty(\mathcal{D}_i^c)} \lesssim_n \lambda$.

We now prove that $\|\Psi_t - \Psi_t^{(i)}\|_{\mathscr{C}_t^{\leq 2}(\mathcal{D}_i^c; D_i)}$ is actually small, uniformly in $t \in [0, \infty]$. In fact, by Lemma 4.30, we also have the bound

$$\|F_s(\Psi_s) - F_s(\Psi_s^{(i)})\|_{L^\infty(\mathcal{D}_i^c; D_i)} \lesssim_n \lambda \, 2^{(\alpha-\beta)s} \|\Psi_s - \Psi_s^{(i)}\|_{\mathscr{C}_s^{\leq 2}(\mathcal{D}_i^c; D_i)},$$

hence

$$\begin{aligned}
\|\Psi_t - \Psi_t^{(i)}\|_{\mathscr{C}_t^{\leq 2}(\mathcal{D}_i^c; D_i)} &\lesssim \int_0^t \|e^{2\xi|\cdot|} \dot{G}_s\|_{L^1} \left( \|F_s(\Psi_s) - F_s(\Psi_s^{(i)})\|_{L_s^\infty(\mathcal{D}_i^c; D_i)} + \|R_t - R_t^{(i)}\|_{L^\infty(\mathcal{D}_i^c)} \right) \mathrm{d}s + 1 \\
&\lesssim_n \lambda \int_0^t 2^{-(d-2\gamma)s} \left( 2^{(\alpha-\beta)s} \|\Psi_s - \Psi_s^{(i)}\|_{\mathscr{C}_s^{\leq 2}(\mathcal{D}_i^c; D_i)} + 1 \right) \mathrm{d}s + 1.
\end{aligned}$$ (4.51)

Since $d - 2\gamma - \alpha + \beta > 0$, for $\lambda$ small enough (4.51) implies uniformly in $t \in [0, \infty]$

$$\|\Psi_t - \Psi_t^{(i)}\|_{L_t^\infty(\mathcal{D}_i^c; D_i)} \lesssim 1,$$

which plugged into (4.48) implies the claim, provided that $(\mathcal{D}_i)_i$ are chosen so that $\min_i \mathrm{dist}(\mathcal{D}_i^c, D_i) \sim \mathrm{dist}(D_1, D_2)$. □

**Remark 4.31.** It is easy to see that the proof extends trivially to the case of $n$-cluster correlation functions. In this case, one finds independent solutions $(\Psi^{(i)})_{i=1,\ldots,n}$ in each of the clusters, and hence, the exponential clustering follows by simply controlling $\|(\Psi - \Psi^{(i)})(f^{(i,k)})\|$ precisely as was done above.

**Proof of Part 3 in Theorem 1.3.** The proof is based on the representation of the two-point function via a generating function. In our setting, a generating function is meaningful provided one introduces some "external" Grassmann variable $\phi$ compatible with the Grassmann martingale $X$ (and in particular with $\psi = X_\infty$) and such that $\omega(\phi) = \omega_t(\phi) = \phi$. This setting requires a simple modification of the NPS discussed in Section 2, whose details we leave to the reader.

For such a Grassmann variable $\phi$, and any $A, B \in \mathscr{A}$ with $B$ anticommuting with $\psi$ we have

$$\omega(A e^{V(\psi) + \phi B}) = \omega(A e^{V(\psi)}) + \phi \, \omega(BA e^{V(\psi)}).$$

Denote $\mathscr{G}_\phi \cong \mathbb{C} \oplus \phi \mathbb{C}$ and note the identity $(a + b\phi)^{-1} = a^{-1}(1 - a^{-1} b \phi)$ for any $a \in \mathbb{C} \setminus \{0\}$ and $b \in \mathbb{C}$. Accordingly, for any $B \in \mathscr{A}$ anticommuting with $\psi$ we can define the continuous linear functional on $\mathscr{A}$

$$\frac{\omega(\cdot e^{V(\psi) + \phi B})}{\omega(e^{V(\psi) + \phi B})} = \omega^V(\cdot) + \phi \, \mathrm{Cov}_{\omega^V}(B; \cdot).$$ (4.52)

For the sake of brevity, denote by $\partial_\phi \colon \mathcal{G}_\phi \to \mathbb{C}$ the projection on the coefficient of $\phi$. Considering now $A = \psi(f_1)$ and $B = \psi(f_2)$, since $\omega(\psi(f_i)) = 0$, $i = 1, 2$, we have

$$\omega^V(\psi(f_1)\psi(f_2)) = -\operatorname{Cov}_{\omega^V}(\psi(f_2); \psi(f_1)) = -\partial_\phi \omega^{\phi, f_2}(\psi(f_1)),$$

where we have introduced the shorthand $\omega^{\phi, f_2}$ for the l.h.s. of (4.52) with $B = \psi(f_2)$. As a consequence of Theorem 3.12 the following stochastic quantisation formula holds for $\omega^{\phi, f}$:

$$\omega^{\phi, f}(A(\psi)) = \omega(A(\Psi^{\phi, f}_\infty)),$$

where $(\Psi^{\phi, f}_t)_{t \geqslant 0}$ solves the following FBSDE on $[0, \infty]$

$$\mathrm{d}\Psi^{\phi, f}_t = \phi \dot{G}_t f \, \mathrm{d}t + \dot{G}_t \omega_t(DV(\Psi^{\phi, f}_\infty)) \mathrm{d}t + \mathrm{d}X_t,$$

$$R^{\phi, f}_t = \int_t^\infty \omega_t(H_{>n;s}(\Psi^{\phi, f}_s)) \mathrm{d}s + \int_t^\infty \omega_t(\langle \dot{G}_s R^{\phi, f}_s + \phi \dot{G}_t f, DF_s(\Psi^{\phi, f}_s)\rangle) \mathrm{d}s,$$

By simple computations, we have

$$\omega^V(\psi(f_1)\psi(f_2)) - \omega(\psi(f_1)\psi(f_2))$$
$$= -\partial_\phi \omega(\Psi^{\phi, f_2}_\infty(f_1)) - \langle\!\langle f_1, G f_2\rangle\!\rangle_\Theta = -\omega(\partial_\phi \mathscr{S}^{\phi, f_2}_\infty(f_1)),$$

where we set

$$\mathscr{S}^{\phi, f_2}_t := \int_0^t \dot{G}_s \omega_s(DV(\Psi^{\phi, f_2}_\infty)) \mathrm{d}s.$$

We now write $\omega_t(DV(\Psi^{\phi, f}_\infty)) = F_t(\Psi^{\phi, f}_t) + R^{\phi, f}_t$, where it is crucial to note that $F$ is independent of $\phi$ and $f$, that is, it can be taken as the solution of the truncated Polchinski equation as before. Note in particular that when $\phi$ is formally set to 0 we have that $\Psi^{\phi=0, f} = \Psi$, and $R^{\phi=0, f} = R$ so that we can write

$$\partial_\phi \Psi^{\phi, f} = \Psi^{\phi, f} - \Psi, \qquad \partial_\phi R^{\phi, f} = R^{\phi, f} - R. \tag{4.53}$$

The pair $(\partial_\phi \Psi^{\phi, f}, \partial_\phi R^{\phi, f})$ solves the forward-backward system

$$\begin{aligned}\partial_\phi \Psi^{\phi, f}_t &= G_t f + \int_0^t \dot{G}_s (\partial_\phi F_s(\Psi^{\phi, f}_s) + \partial_\phi R^{\phi, f}_s) \mathrm{d}s \\ \partial_\phi R^{\phi, f}_t &= \int_t^\infty \omega_t(\partial_\phi H_{>n;s}(\Psi^{\phi, f}_s)) \mathrm{d}s + \int_t^\infty \omega_t(DF_s(\Psi_s) \cdot \dot{G}_s \partial_\phi R^{\phi, f}_s) \mathrm{d}s \\ &\quad + \int_t^\infty \omega_t(DF_s(\Psi_s) \cdot \dot{G}_s f) \mathrm{d}s + \int_t^\infty \omega_t(\partial_\phi DF_s(\Psi^{\phi, f}_s) \cdot \dot{G}_s R_s) \mathrm{d}s\end{aligned} \tag{4.54}$$

and we have the representation

$$\partial_\phi \mathscr{S}^{\phi, f_2}_t = \int_0^t \dot{G}_s (\partial_\phi F_s(\Psi^{\phi, f_2}_s) + \partial_\phi R^{\phi, f_2}_s) \mathrm{d}s. \tag{4.55}$$

Along the lines of the proof of Proposition 4.20, we can show that the system (4.54) with $f$ such that $\|f\|_{L^1} = 1$ has a unique global solution $(\partial_\phi \Psi^{\phi, f}, \partial_\phi R^{\phi, f})$ satisfying

$$|\!|\!|(\partial_\phi \Psi^{\phi, f}, \partial_\phi R^{\phi, f})|\!|\!|_* \lesssim 1, \tag{4.56}$$

with $\sigma = 2\gamma$, $\sigma' = 0$, $a = D$, $a' = c_n \lambda$ and $\eta = 0$, these constants appearing in the definition of the norm (4.28). First of all, noting that $\partial_\phi F_s(\Psi^{\phi, f}_s) = F_s(\Psi^{\phi, f}_s) - F_s(\Psi_s)$ and so on, and taking (4.53) into account, the following estimates are a consequence of Lemma 4.22

$$\begin{aligned}\|\partial_\phi F_s(\Psi^{\phi, f}_s)\|_{L^\infty} &\lesssim \lambda \, 2^{(\alpha-\beta)s} \|\partial_\phi \Psi^{\phi, f_2}_s\|_{\mathscr{C}^{\leqslant 2}_s}, \\ \|\partial_\phi H_{>n;s}(\Psi^{\phi, f}_s)\|_{L^\infty} &\lesssim \lambda \, 2^{(\alpha-n\kappa-\gamma)s} \|\partial_\phi \Psi^{\phi, f_2}_s\|_{\mathscr{C}^{\leqslant 2}_s}, \\ \|\partial_\phi DF_s(\Psi^{\phi, f}_s)\|_{L^\infty} &\lesssim \lambda \, 2^{(\alpha-\beta-\gamma)s} \|\partial_\phi \Psi^{\phi, f_2}_s\|_{\mathscr{C}^{\leqslant 2}_s}.\end{aligned}$$

Then, by Young's inequality for convolutions we have the bound $\|G_t f\|_{\mathscr{C}_t^{\leqslant 2}} \lesssim 2^{2\gamma t}$ whereas the drift is estimated as

$$\|\dot{G}_s(\partial_\phi F_s(\Psi_s) + \partial_\phi R_s^{\phi,f})\|_{\mathscr{C}_s^{\leqslant 2}}$$
$$\lesssim 2^{-(d-2\gamma)s}(\lambda 2^{(\alpha-\beta)s}\|\partial_\phi \Psi_s^{\phi,f_2}\|_{\mathscr{C}_s^{\leqslant 2}} + \|\partial_\phi R_s^{\phi,f}\|_{L^\infty}).$$

On the other hand, the terms in the equation for $\partial_\phi R_t^{\phi,f}$ are bounded as follows

$$\|\omega_t(\partial_\phi H_{>n;s}(\Psi_s^{\phi,f}))\|_{L^\infty} \lesssim \lambda 2^{(\alpha-n\kappa-\gamma)s}\|\partial_\phi \Psi_s^{\phi,f}\|_{L^\infty},$$

$$\|\omega_t(DF_s(\Psi_s) \cdot \dot{G}_s \partial_\phi R_s^{\phi,f})\|_{L^\infty} \lesssim \|DF_s(\Psi_s)\|_{\mathscr{B}(L^\infty;L^\infty)} \|\dot{G}_s\|_{L^1} \|\partial_\phi R_s^{\phi,f}\|_{L^\infty}$$
$$\lesssim 2^{-(d-2\gamma-\alpha+\beta)s}\|\partial_\phi R_s^{\phi,f}\|_{L^\infty},$$

$$\|\omega_t(DF_s(\Psi_s) \cdot \dot{G}_s f)\|_{L^\infty} \lesssim \|DF_s(\Psi_s)\|_{\mathscr{B}(L^\infty;L^\infty)} \|\dot{G}_s\|_{L^1} \|f\|_{L^1}$$
$$\lesssim 2^{-(d-2\gamma-\alpha+\beta)s},$$

$$\|\omega_t(\partial_\phi DF_s(\Psi_s^{\phi,f}) \cdot \dot{G}_s R_s)\|_{L^\infty} \lesssim \|\partial_\phi DF_s(\Psi_s^{\phi,f})\|_{L^\infty} \|\dot{G}_s\|_{L^1} \|R_s\|_{L^\infty}$$
$$\lesssim 2^{-(d-\gamma-\alpha+\beta)}\|\partial_\phi \Psi_s^{\phi,f_2}\|_{\mathscr{C}_s^{\leqslant 2}}.$$

From these bounds, following exactly the same steps carried out in the proof of Proposition 4.20 we have that the system (4.54) has a unique solution in the ball defined in (4.56).

We will henceforth fix $\varepsilon > 0$ and take $f_i := \chi_{x_i}^\varepsilon$ where $\chi_x^\varepsilon \in C_c^\infty$ is as in the statement of the theorem, so that $\|f_i\|_{L^1} = 1$. Set $B_i := \{x \in \mathbb{R}^d | |x - x_i| \leqslant \varepsilon\}$, $i = 1, 2$. For $\xi$ small enough, using $\|f_2\|_{L^1} = 1$ and the estimates in Corollary 3.24 we have

$$\|G_t f_2\|_{\mathscr{C}_t^n(\mathbb{R}^d;B_1)} \leqslant \int_0^t \|\dot{G}_s f_2\|_{\mathscr{C}_s^n(\mathbb{R}^d;B_1)} \, ds$$
$$\leqslant \int_0^t \sup_\nu 2^{-ns} \|\partial_\nu^n \dot{G}_s\|_{L^\infty_{s;\xi}} 2^{-\xi 2^s \text{dist}(B_2,B_1)} \|f_2\|_{L^1}$$
$$\lesssim \int_0^t 2^{-\xi 2^s \text{dist}(B_1,B_2)} 2^{2\gamma s} ds$$
$$\lesssim \text{dist}(B_1, B_2)^{-2\gamma}.$$

where the norm $\|\cdot\|_{\mathscr{C}_t^n(\mathcal{D};B)}$ for measurable sets $\mathcal{D}, B \subseteq \mathbb{R}^d$ was defined in (4.44). Furthermore, by the triangle inequality

$$\int_0^t \|\dot{G}_s \partial_\phi F_s(\Psi_s^{\phi,f_2})\|_{\mathscr{C}_s^{\leqslant n}(\mathbb{R}^d;B_1)} ds$$
$$\lesssim \int_0^t \sup_\nu 2^{-ns} \|\dot{G}_s\|_{L^1_{\xi,1,s}} \|\partial_\phi \Psi_s^{\phi,f_2}\|_{\mathscr{C}_s^{\leqslant 2}(\mathbb{R}^d;B_1)} \|DF_s(\Psi_s)\|_{\mathscr{B}_s(L^\infty;L^\infty)} ds$$
$$\lesssim \lambda \int_0^t 2^{-(d-2\gamma-\alpha+\beta)s} \|\partial_\phi \Psi_s^{\phi,f_2}\|_{\mathscr{C}_s^{\leqslant 2}(\mathbb{R}^d;B_1)} ds$$
$$\lesssim \lambda \sup_{s \in \mathbb{R}^+} \|\partial_\phi \Psi_s^{\phi,f_2}\|_{\mathscr{C}_s^{\leqslant 2}(\mathbb{R}^d;B_1)}.$$

Furthermore,

$$\int_0^t \|\dot{G}_s \partial_\phi R_s^{\phi,f_2}\|_{L_t^\infty(\mathbb{R}^d;B_1)} ds \leqslant \int_0^t \|\dot{G}_s\|_{L^1_{s;\xi}} \|\partial_\phi R_s^{\phi,f_2}\|_{L^\infty} ds \lesssim \lambda,$$

hence

$$\sup_{t \in \mathbb{R}^+} \|\partial_\phi \Psi_t^{\phi,f_2}\|_{\mathscr{C}_t^{\leqslant 2}(\mathbb{R}^d;B_1)} \lesssim \lambda + \text{dist}(B_1,B_2)^{-2\gamma} + \lambda \sup_{t \in \mathbb{R}^+} \|\partial_\phi \Psi_t^{\phi,f_2}\|_{\mathscr{C}_t^{\leqslant 2}(\mathbb{R}^d;B_1)}$$

implying, for $\text{dist}(B_1,B_2)$ small enough

$$\sup_{t \in \mathbb{R}^+} \|\partial_\phi \Psi_t^{\phi,f_2}\|_{\mathscr{C}_t^{\leqslant 2}(\mathbb{R}^d;B_1)} \lesssim \text{dist}(B_1,B_2)^{-2\gamma}.$$

We plug this bound into the equation for $\partial_\phi \mathscr{S}_t^{\phi,f_2}(f_1)$, see (4.55), and obtain for $\xi$ small enough

$$\begin{aligned}\|\partial_\phi \mathscr{S}_t^{\phi,f_2}\|_{L_t^\infty(\mathbb{R}^d;B_1)} &\lesssim \int_0^t 2^{-(d-2\gamma)s}(\mathrm{dist}(B_1,B_2)^{-2\gamma}\lambda\, 2^{(\alpha-\beta)s}+\lambda)\mathrm{d}s \\ &\lesssim \lambda\,\mathrm{dist}(B_1,B_2)^{-2\gamma},\end{aligned}$$

where we used once more that $d-2\gamma-\alpha+\beta>0$. Since $f_1$ is supported on $B_1$, we have $\|\partial_\phi \mathscr{S}_t^{\phi,f_2}(f_1)\|_{L_t^\infty(\mathbb{R}^d;B_1)} \lesssim \|\partial_\phi \mathscr{S}_t^{\phi,f_2}\|_{L_t^\infty(\mathbb{R}^d;B_1)}\|f_1\|_{L^1}$, hence the claim. □

# Appendix A  Proof of Theorem 2.20

In this section we prove the existence of GBM as described in Theorem 2.20 and in Definition 3.18. More precisely, we prove the existence of a NPS where our GBMs are realised as suitable operators acting on a Hilbert space. The construction of Grassmann fields as operators acting on a Hilbert space was already obtained in the early days of Euclidean quantum field theory [OS72, OS73, Ost73], see also [GJ87]. More recently, this algebraic approach has been further developed in [ABDG20] with applications to Grassmann SDE. In addition to what was already done there, here we construct a conditional expectation and discuss GBMs in the framework of FNPST.

The proof of Theorem 2.20 is by construction and divided into three steps. First of all, we associate a FNPS with the pair $(\mathcal{H},(P_t)_{t\geqslant 0})$ consisting of a Hilbert space $\mathcal{H}$, which in our case will be chosen to be $L^2(\mathbb{R}_+;\mathfrak{h})$, and a net of projections $(P_t)_{t\geqslant 0}$. This association will be called Fock construction for obvious reasons. Then, we construct a conditional expectation associated with the Fock construction. Finally, we provide an explicit representation of the GBM by using the Osterwalder–Schrader construction [OS72, OS73, Ost73].

**Definition A.1.** *(Fock construction). Let $\mathcal{H}$ be a separable Hilbert space and $(P_t)_{t\geqslant 0}$ an increasing orthogonal net of projections on it, that is, for any $s\leqslant t$ we have*

$$P_t P_s = P_{t\wedge s} = P_s P_t. \tag{A.1}$$

*The Fock construction associated with the pair $(\mathcal{H},(P_t)_{t\geqslant 0})$ is the FNPS $(\mathcal{M},\omega,(\mathcal{M}_t)_{t\geqslant 0})$, where $\omega(\cdot)\coloneqq \langle\Omega,\cdot\,\Omega\rangle$ with $\Omega$ being the Fock vacuum on $\Gamma_a(\mathcal{H})$, the antisymmetric Fock space generated by $\mathcal{H}$, and where*

$$\mathcal{M}\coloneqq \mathrm{CAR}(\mathcal{H}),\qquad \mathcal{M}_t\coloneqq \mathrm{CAR}(P_t\mathcal{H})\subseteq \mathcal{M}.$$

**Remark A.2.** Verifying that the triple $(\mathcal{M},\omega,(\mathcal{M}_t)_{t\geqslant 0})$ is a FNPS is a simple check, see the definitions in Section 2.

Although the Fock construction is very simple, there is no general theory of conditional expectations, because the vector $\Omega$, does not have the modular property [Tom67, Tak70, Tak72]. If this were the case, on the other hand, the existence of a conditional expectation of $\mathcal{M}$ onto $\mathcal{M}_t$ would be a consequence of a theorem by Takesaki [Tak72]. Nonetheless, a conditional expectation can be explicitly defined by means of the fermionic exponential law [Der06].

For any $s\leqslant t$, we write

$$\mathcal{H}_t\coloneqq P_t\mathcal{H},\quad \mathcal{H}_{s,t}\coloneqq (\mathbb{1}-P_s)P_t\mathcal{H}=(\mathbb{1}-P_s)\mathcal{H}_t,\quad \mathcal{H}_{>t}\coloneqq (\mathbb{1}-P_t)\mathcal{H}.$$

By the property (A.1), we have that for any $s \leqslant t$ the subspaces $\mathcal{H}_s$, $\mathcal{H}_{s,t}$ and $\mathcal{H}_{>t}$ are orthogonal and we get $\mathcal{H} = \mathcal{H}_s \oplus \mathcal{H}_{s,t} \oplus \mathcal{H}_{>t} = \mathcal{H}_t \oplus \mathcal{H}_{>t}$.

The exponential law for fermions is the unitary operator $U_t : \Gamma_a(\mathcal{H}_t \oplus \mathcal{H}_{>t}) \to \Gamma_a(\mathcal{H}_t) \otimes \Gamma_a(\mathcal{H}_{>t})$, such that

$$U_t a(f) U_t^* = a(P_t f) \otimes \mathbb{1} + \Xi_t \otimes a((\mathbb{1} - P_t) f), \qquad U_t \Omega = \Omega_t \otimes \Omega_{>t}, \qquad (A.2)$$

where $\Xi_t = (-1)^{N_t}$[A.1], where $N_t$ is the number operator on $\Gamma_a(\mathcal{H}_t)$, and $\Omega, \Omega_t, \Omega_{>t}$ are the ground states of $\Gamma_a(\mathcal{H})$, $\Gamma_a(\mathcal{H}_t)$ and $\Gamma_a(\mathcal{H}_{>t})$ respectively. Furthermore, for any $s < t$, we can define the unitary operator $U_{s,t} : \Gamma_a(\mathcal{H}_s \oplus \mathcal{H}_{s,t} \oplus \mathcal{H}_{>t}) \to \Gamma_a(\mathcal{H}_s) \otimes \Gamma_a(\mathcal{H}_{s,t}) \otimes \Gamma_a(\mathcal{H}_{>t})$ such that

$$U_{s,t} a(f) U_{s,t}^* = a(P_s f) \otimes \mathbb{1} \otimes \mathbb{1} + \Xi_s \otimes a(P_{s,t} f) \otimes \mathbb{1} + \Xi_s \otimes \Xi_{s,t} \otimes a((\mathbb{1} - P_t) f) \qquad (A.3)$$

where $P_{s,t} = P_t(\mathbb{1} - P_s) = (\mathbb{1} - P_s) P_t = P_t - P_s$, and

$$U_{s,t} \Omega = \Omega_s \otimes \Omega_{s,t} \otimes \Omega_{>t}$$

where $\Omega_{s,t}$ is the ground state of $\Gamma_a(\mathcal{H}_{s,t})$, and $\Xi_{s,t} = (-1)^{N_{s,t}}$ where $N_{s,t}$ is the number operator of $\Gamma_a(\mathcal{H}_{s,t})$. It is important to note the following relation between $U_{s,t}$, $U_s$ and $U_t$, namely, for any $D \in \mathrm{CAR}(\mathcal{H})$ we have

$$U_{s,t} D U_{s,t}^* = (U_s \otimes \mathbb{1})(U_t D U_t^*)(U_s^* \otimes \mathbb{1}) = (\mathbb{1} \otimes U_t)(U_s D U_s^*)(\mathbb{1} \otimes U_t^*). \qquad (A.4)$$

**Remark A.3.** Hereafter if $\mathcal{K} \subset \mathcal{H}$ we identify $\mathrm{CAR}(\mathcal{K})$ with a subspace of $\mathrm{CAR}(\mathcal{H})$ in the following way

$$\begin{aligned}\mathrm{CAR}(\mathcal{K}) \\ = \mathrm{span}\{a(f_1) \cdots a(f_n) a^*(g_1) \cdots a^*(g_m) | n, m \in \mathbb{N}, (f_i)_{i=1}^n, (g_i)_{i=1}^m \subset \mathcal{K}, \}\end{aligned} \qquad (A.5)$$

In this way if $A \in \mathcal{B}(\Gamma_a(\mathcal{K})) \cap \mathrm{CAR}(\mathcal{K})$ we can define $\omega(A)$ as $\omega(A) = \langle \Omega, A \Omega \rangle$, where $\Omega$ is the ground state of $\Gamma_a(\mathcal{H})$ and $A$ is identified as an operator of $\mathrm{CAR}(\mathcal{H}) \subset \mathcal{B}(\Gamma_a(\mathcal{H}))$ thanks to the relation (A.5). With the previous notation, for any $s < t$ and $A \in \mathrm{CAR}(\mathcal{H}_s)$, $B \in \mathrm{CAR}(\mathcal{H}_{s,t})$ and $C \in \mathrm{CAR}(\mathcal{H}_{>t})$, we have

$$\omega(A) = \langle \Omega_s, A \Omega_s \rangle, \quad \omega(B) = \langle \Omega_{s,t}, B \Omega_{s,t} \rangle, \qquad \omega(C) = \langle \Omega_{>t}, C \Omega_{>t} \rangle.$$

Using the fact that $U_t$ is an isometry which sends the ground state of $\Gamma_a(\mathcal{H})$ into the one of $\Gamma_a(\mathcal{H}_t) \otimes \Gamma_a(\mathcal{H}_{>t})$, and using the previously described notation we get the equality

$$\omega\left(U_t^* \left(\sum_{i=1}^n A_i \otimes B_i\right) U_t\right) = \sum_{i=1}^n \omega(A_i) \omega(B_i) \qquad (A.6)$$

for any $t \in \mathbb{R}_+$, let $A_i \in \mathrm{CAR}(\mathcal{H}_t)$ and $B_i \in \mathrm{CAR}(\mathcal{H}_{>t})$.

The exponential law gives another representation of the CAR algebra we consider, as a *-sub-algebra of the algebra of bounded operators on $\Gamma_a(\mathcal{H}_t) \otimes \Gamma_a(\mathcal{H}_{>t})$. On the tensor product $\mathcal{B}(\Gamma_a(\mathcal{H}_t)) \otimes \mathcal{B}(\Gamma_a(\mathcal{H}_{>t}))$, we introduce the projection map $\tilde{\Omega}_{>t} : \mathcal{B}(\Gamma_a(\mathcal{H}_t)) \otimes \mathcal{B}(\Gamma_a(\mathcal{H}_{>t})) \to \mathcal{B}(\Gamma_a(\mathcal{H}_t)) \otimes \mathcal{B}(\Gamma_a(\mathcal{H}_{>t}))$:

$$\tilde{\Omega}_{>t}(A \otimes B) := A \otimes \mathbb{1} \langle \Omega_{>t}, B \Omega_{>t} \rangle.$$

The following properties hold true.

**Lemma A.4.** *The map $\tilde{\Omega}_{>t}$ is positive and continuous in the operator norm.*

---

A.1. Note that $\Xi_t$ is sometimes called parity operator of $\mathcal{H}_t$.

**Proof.** Consider a generic element $M = \sum_j A_j \otimes B_j \in \mathscr{B}(\Gamma_a(\mathcal{H}_t)) \otimes \mathscr{B}(\Gamma_a(\mathcal{H}_{>t}))$, the sum being finite. Set $\tilde{A}_j := A_j \langle \Omega_{>t}, B_j \Omega_{>t} \rangle$ and $A := \sum_j \tilde{A}_j$. To prove positivity, we note that an operator of the form $A \otimes \mathbb{1}$ is positive provided that $A$ is positive. We have

$$\langle \psi, A\psi \rangle_{\Gamma_a(\mathcal{H}_t)} = \sum_j \langle \psi, A_j \psi \rangle_{\Gamma_a(\mathcal{H}_t)} \langle \Omega_{>t}, B_j \Omega_{>t} \rangle_{\Gamma_a(\mathcal{H}_{>t})}$$

$$= \langle \psi \otimes \Omega_{>t}, M \psi \otimes \Omega_{>t} \rangle_{\Gamma_a(\mathcal{H}_t) \otimes \Gamma_a(\mathcal{H}_{>t})} \geq 0.$$

To prove continuity, we note that $\|A \otimes \mathbb{1}\|_{\mathscr{B}(\Gamma_a(\mathcal{H}_t) \otimes \Gamma_a(\mathcal{H}_{>t}))} \leq \|A\|_{\mathscr{B}(\Gamma_a(\mathcal{H}_t))}$, hence we have the bound:

$$\|\tilde{\Omega}_{>t}(M)\|_{\mathscr{B}(\Gamma_a(\mathcal{H}_t) \otimes \Gamma_a(\mathcal{H}_{>t}))} = \|A \otimes \mathbb{1}\|_{\mathscr{B}(\Gamma_a(\mathcal{H}_t) \otimes \Gamma_a(\mathcal{H}_{>t}))}$$

$$\leq \|A\|_{\mathscr{B}(\Gamma_a(\mathcal{H}_t))} = \sup_{\psi: \|\psi\|_{\Gamma_a(\mathcal{H}_t)}=1} \left| \sum_j \langle \psi, A_j \psi \rangle \langle \Omega_{>t}, B_j \Omega_{>t} \rangle \right|$$

$$\leq \|M\|_{\mathscr{B}(\Gamma_a(\mathcal{H}_t) \otimes \Gamma_a(\mathcal{H}_{>t}))}. \qquad \square$$

As a consequence of Lemma A.4, we can extend $\tilde{\Omega}_{>t}$ to the whole CAR algebra of $\Gamma_a(\mathcal{H}_t) \otimes \Gamma_a(\mathcal{H}_{>t})$. We finally define the conditional expectation on $\mathcal{M}_t$.

**Definition A.5.** *(Conditional Expectation). Let $(\Omega, \mathcal{M}, (\mathcal{M}_t)_{t \geq 0})$ be the Fock construction associated with $(\mathcal{H}, (P_t)_{t \geq 0})$. Let $U_t$ be the fermionic exponential law and $\tilde{\Omega}_{>t}$ be the (extension of the) projection defined above. We define the map $\omega_t: \mathcal{M} \to \mathcal{M}_t$ by:*

$$\omega_t(\cdot) := U_t^* \tilde{\Omega}_{>t}(U_t \cdot U_t^*) U_t.$$

Using the notation of Remark A.3, for any $A_i \in \mathrm{CAR}(\mathcal{H}_t)$ and $B_i \in \mathrm{CAR}(\mathcal{H}_{>t})$ we have

$$\omega_t\left(U_t^*\left(\sum_{i=1}^n A_i \otimes B_i\right) U_t\right) = U_t^* \tilde{\Omega}_{>t}\left(\sum_{i=1}^n A_i \otimes B_i\right) U_t = \sum_{i=1}^n \omega(B_i) A_i. \tag{A.7}$$

**Proposition A.6.** *The map $\omega_t: \mathcal{M} \to \mathcal{M}_t$ defined above is a conditional expectation.*

**Proof.** The property $\omega_t(\mathbb{1}_{\mathcal{M}}) = \mathbb{1}_{\mathcal{M}_t}$ is a straightforward consequence of the definition. The positivity of $\omega_t(\cdot)$ follows because unitary conjugation preserves positivity and because the map $\tilde{\Omega}_{>t}$ is positive, see Lemma A.4. We are left with checking the property

$$\omega_t(N_1 M N_2) = N_1 \omega_t(M) N_2,$$

for any $M \in \mathcal{M}$ and $N_1, N_2 \in \mathcal{M}_t$. To this end, it suffices to note that operators such as $N_1$ and $N_2$ are mapped to operators of the form $U_t N_i U_t^* = \tilde{N}_i \otimes \mathbb{1}$, for some $\tilde{N}_i \in \mathscr{B}(\Gamma_a(\mathcal{H}_t))$. Therefore, if $U_t M U_t^* = \sum_j A_j \otimes B_j$, the sum being finite, we have

$$\begin{aligned}\omega_t(N_1 M N_2) &= U_t^* \tilde{\Omega}_{>t}(U_t N_1 M N_2 U_t^*) U_t \\ &= U_t^* (\tilde{N}_1 \otimes \mathbb{1}) \sum_j (A_j \langle \Omega_{>t}, B_j \Omega_{>t}\rangle \otimes \mathbb{1})(\tilde{N}_2 \otimes \mathbb{1}) U_t \\ &= N_1 \omega_t(M) N_2.\end{aligned}$$

By the continuity of $\omega_t(\cdot)$ this extends to any $M \in \mathcal{M}$.

The fact that $\omega(\omega_t(D)) = \omega(D)$ for any $D \in \mathrm{CAR}(\mathcal{H})$, follows from the continuity of $\omega_t$ and from (A.6). We have now to prove the tower property, i.e. for any $D \in \mathrm{CAR}(\mathcal{H})$ we have $\omega_s(\omega_t(K)) = \omega_s(K)$. By the continuity of $\omega_s, \omega_t$ it is enough to prove the tower property for $D$ of the form

$$D = U_{s,t}^*\left(\sum_{i=1}^n A_i \otimes B_i \otimes C_i\right) U_{s,t}$$

for some $n \in \mathbb{N}$ and $A_i \in \mathrm{CAR}(\mathcal{H}_s)$, $B_i \in \mathrm{CAR}(\mathcal{H}_{s,t})$ and $C_i \in \Gamma_a(\mathcal{H}_{>t})$. By (A.4) and (A.7) we have

$$\omega_t\left(U_{s,t}^*\left(\sum_{i=1}^n A_i \otimes B_i \otimes C_i\right) U_{s,t}\right) = \omega_t\left(U_t^*\left(\sum_{i=1}^n (U_s^*(A_i \otimes B_i) U_s) \otimes C_i\right) U_t\right)$$
$$= \sum_{i=1}^n \omega(C_i) U_s^*(A_i \otimes B_i) U_s. \tag{A.8}$$

Thus, again by (A.4), (A.7) and (A.8), we get

$$\omega_s\left(\omega_t\left(\sum_{i=1}^n U_{s,t}^*(A_i \otimes B_i \otimes C_i) U_{s,t}\right)\right) = \omega_s\left(\sum_{i=1}^n \omega(C_i) U_s^*(A_i \otimes B_i) U_s\right)$$
$$= \sum_{i=1}^n \omega(C_i) \omega_s(U_s^*(A_i \otimes B_i) U_s)$$
$$= \sum_{i=1}^n A_i \omega(C_i) \omega(B_i)$$
$$= \omega_s\left(\sum_{i=1}^n U_s^*[A_i \otimes \{U_t^*(B_i \otimes C_i) U_t\}] U_s\right)$$
$$= \omega_s\left(U_{s,t}^*\left(\sum_{i=1}^n A_i \otimes B_i \otimes C_i\right) U_{s,t}\right).$$

$\square$

Before providing the construction of the GBM, let us set some further notation. If the GBM is indexed by the Hilbert space $\mathfrak{h}$, with conjugation $\Theta$, we let $\mathcal{H} := L^2(\mathbb{R}; \mathfrak{h})$ and $P_t := \hat{P}_t \otimes \mathbb{1}$ be an increasing net of projections on it, $\hat{P}_t$ be the projection of $L^2(\mathbb{R})$ onto $L^2([0,t])$.

**Definition A.7.** *Let $\mathfrak{h}$ be a Hilbert space with conjugation $\Theta$, let $(\mathcal{H}, (P_t)_{t \geq 0})$ be as above and $(\mathcal{M}, \omega, (\mathcal{M}_t)_{t \geq 0})$ be its Fock construction. Let $(G_t)_t \subset \mathcal{B}(\mathfrak{h})$ be $\Theta$-antisymmetric and differentiable in t and let $\dot{G}_t = C_t^* C_t U_t$ be a decomposition with unitary $U_t$. Define the processes*

$$X_t(f) := a^*\left(\int_0^t \delta_s \otimes C_s U_s f \mathrm{d}s\right) + a\left(\int_0^t \delta_s \otimes C_s \Theta f \mathrm{d}s\right), \tag{A.9}$$

*indexed by $\mathfrak{h}$.*

**Remark A.8.** Note that if $f \in \mathfrak{h}$, the functions $g_t := \int_0^t \delta_s \otimes C_s U_s f \mathrm{d}s$ and $\tilde{g}_t := \int_0^t \delta_s \otimes C_s \Theta f \mathrm{d}s$ are in $\mathcal{H}$. In fact, $\|g_t\|_{\mathcal{H}}^2 = \int_0^t \|C_s U_s f\|_{\mathfrak{h}}^2 \mathrm{d}s$ and likewise for $\tilde{g}_t$.

The following proposition completes the proof of Theorem 2.20.

**Proposition A.9.** *The process $(X_t)_{t \geq 0}$ defined in (A.9) are a norm-compatible GBM with covariance $G_t$.*

**Proof.** That $(X_t)_{t \geq 0}$ is adapted follows by direct inspection. We prove that it is a martingale, that is, $\omega_s(X_t(f)) = X_s(f)$ for any $0 \leq s \leq t$. For the sake of brevity, we let $g_{r,t} := \int_r^t \delta_s \otimes C_s U_s f \mathrm{d}s$ and $\tilde{g}_{r,t} := \int_r^t \delta_s \otimes C_s \Theta f \mathrm{d}s$. We let $\mathcal{H} = \mathcal{H}_s \oplus \mathcal{H}_{>s}$ and let $U_s$ be the corresponding fermionic exponential law. Note that $P_s g_{0,t} = g_{0,s}$ and $P_s \tilde{g}_{0,t} = \tilde{g}_{0,s}$ We have:

$$U_s X_t(f) U_s^* = a^*(g_{0,s}) \otimes \mathbb{1} + \Xi \otimes a^*(g_{s,t}) + a(\tilde{g}_{0,s}) \otimes \mathbb{1} + \Xi \otimes a(\tilde{g}_{s,t}).$$

Therefore

$$\omega_s(X_t(f)) = U_s^* \tilde{\Omega}_{>s}(U_s X_s(f) U_s^*) U_s = U_s^*(a^*(g_{0,s}) \otimes \mathbb{1} + a(\tilde{g}_{0,s}) \otimes \mathbb{1}) U_s = X_s(f).$$

Next, we prove that it is a centred Grassmann Gaussian process. We begin by checking the anti-commutation relations:

$$\begin{aligned}\{X_t(f), X_{t'}(f')\} &= \left\langle \int_0^t \delta_s \otimes C_s \Theta f \, ds, \int_0^{t'} \delta_r \otimes C_r U_r f' \, dr \right\rangle_{\mathcal{H}} \\ &\quad + \left\langle \int_0^{t'} \delta_s \otimes C_s \Theta f' \, ds, \int_0^t \delta_r \otimes C_r U_r f \, dr \right\rangle_{\mathcal{H}} \\ &= \int_0^{t \wedge t'} (\langle C_s \Theta f, C_s U_s f' \rangle_{\mathfrak{h}} + \langle C_s \Theta f', C_s U_s f \rangle_{\mathfrak{h}}) \, ds \\ &= \langle\!\langle f, G_{t \wedge t'} f' \rangle\!\rangle_\Theta + \langle\!\langle f', G_{t \wedge t'} f \rangle\!\rangle_\Theta = 0, \end{aligned}$$

where in the last step we used that $G_s$ is $\Theta$-antisymmetric. Then, we show that $X_t$ is a centred Gaussian. We compute

$$\left\{ a\left( \int_0^t \delta_r \otimes C_r U_r f \, dr \right), X_s(f') \right\} = \int_0^{t \wedge s} \langle C_r \Theta f, C_r U_r f' \rangle_{\mathfrak{h}} \, dr = \langle\!\langle f, G_{t \wedge s} f' \rangle\!\rangle_\Theta,$$

hence, the claim:

$$\begin{aligned}\omega(X_t(f) X_{t_1}(f_1) \cdots X_{t_n}(f_n)) &= \omega\left(a\left(\int_0^t \delta_r \otimes C_r U_r f \, dr\right) X_{t_1}(f_1) \cdots X_{t_n}(f_n)\right) \\ &= \sum_{j=1}^n (-)^j \langle\!\langle f, G_{t \wedge t_j} f_j \rangle\!\rangle_\Theta \, \omega\left(X_{t_1}(f_1) \cdots X_{t_j}(\!\!\!/f_j) \cdots X_{t_n}(f_n)\right)\end{aligned}$$

Finally, the norm-compatibility is a direct consequence of $\|a(f)\|, \|a^*(f)\| \leq \|f\|_{\mathfrak{h}}$. This concludes the proof of the proposition. □

Note in particular that Proposition A.9 implies that Definitions 3.18 and 4.27 are meaningful and that we obtain $(X_t^{L,\varepsilon})_t$ and $(X_t^{(i)})_t$ simply by plugging $\mathfrak{C}_s^{L,\varepsilon}$ and $\mathbf{1}_{D_i} \mathfrak{C}_s$ in place of $C_s$ in eq. (A.9).

**Remark A.10.** Let us fix $\mathfrak{h} := L^2(\mathbb{R}^d; \mathbb{C}^n)$ and some unitary $U$ for concreteness. In this case, we could introduce Grassmann white noise (associated with $U$ and $\Theta$) as the operator-valued distribution

$$\xi(ds, dz) := a^*(\delta_s ds \otimes U(\cdot, z) dz) + a(\delta_s ds \otimes \Theta(\cdot, z) dz),$$

so that, if $C_t$ commutes with $U$ and $\Theta$, we can write

$$X_t(f) = \int_0^t \int_{\mathbb{R}^d} (C_s f)(z) \, \xi(ds, dz).$$

This actually holds in our setting, see Section 3.3 and clarifies the meaning of Definition 3.18.